\documentclass[leqno]{amsart}
\usepackage{egastyle}
\usepackage{mymacros}
\usepackage[all,tips,dvips]{xy}
\usepackage{url}
\usepackage{bm}

\CompileMatrices

\newcommand\rig{\mathrm{rig}}
\newcommand\Ab{\operatorname{\mathbf{Ab}}}
\newcommand\Nil{\operatorname{\mathbf{Nil}}}
\newcommand\Trop{\operatorname{Trop}}
\newcommand\ev{\operatorname{ev}}

\newcommand\BT{\mathrm{BT}}

\newcommand\B{\mathbf{B}}
\newcommand\D{\mathbf{D}}
\newcommand\bflog{\operatorname{\mathbf{log}}}

\newcommand\face{\operatorname{face}}

\newcommand\New{\operatorname{New}}
\newcommand\Mod{\operatorname{\mathbf{Mod}}}
\DeclareMathOperator\Lie{Lie}
\DeclareMathOperator\hexp{hexp}
\DeclareMathOperator{\inn}{in}
\DeclareMathOperator{\Inn}{In}
\DeclareMathOperator{\minset}{minset}
\DeclareMathOperator{\conv}{conv}
\DeclareMathOperator{\size}{size}
\DeclareMathOperator{\vertices}{vertices}
\DeclareMathOperator{\relint}{relint}

\newcommand\MU{\bm\mu}
\newcommand\dV{{\mathrm{d}V}}
\newcommand\T{\bT}

\title{Higher-level canonical subgroups for $p$-divisible groups}
\author{Joseph Rabinoff}

\begin{document}

\begin{abstract}
  Let $R$ be a complete rank-$1$ valuation ring of mixed characteristic
  $(0,p)$, and let $K$ be its field of fractions.  A $g$-dimensional
  truncated Barsotti-Tate group $G$ of level $n$ over $R$ is said to have
  a \emph{level-$n$ canonical subgroup} if there is a $K$-subgroup of 
  $G\tensor_R K$ with geometric structure $(\Z/p^n\Z)^g$ consisting of
  points ``closest to zero''.  We give a nontrivial condition on the Hasse
  invariant of $G$ that guarantees the existence of the canonical subgroup,
  analogous to a result of Katz and Lubin for elliptic curves.  The bound
  is independent of the height and dimension of $G$.
\end{abstract}

{\let\thefootnote\relax
\footnote{MSC subject classes: 11G10 (primary), 14L05 (secondary)}
\footnote{Keywords: abelian varieties, $p$-divisible groups, 
  canonical subgroups, automorphic forms, tropical geometry}
}

\maketitle

\begin{center}
  Version dated: \textbf{\today}.
\end{center}

\section{Introduction}

\paragraph[Motivation] Let $K$ be a field that is complete with respect to a nontrivial
non-Archimedean valuation, with residue characteristic $p$.  The
level-$n$ canonical subgroup of an abelian variety $A$ of dimension $g$
defined over $K$ is a certain distinguished $K$-subgroup of $A[p^n]$ with
geometric structure $(\Z/p^n\Z)^g$; it should exist when the reduction type of
$A$ is ``not too far from ordinary''.  The classical theory of the canonical
subgroup (for elliptic curves), as introduced by Katz and
Lubin~\cite{katz:padicmfs}, has seen important applications in the
study of overconvergent $p$-adic modular forms; see for instance the work
of Buzzard and Taylor~\cite{buzzard_taylor:companion_forms},
Buzzard~\cite{buzzard:analytic_continuation}, and
Kassaei~\cite{kassaei:gluing_lemma}.  In an effort to extend these
overconvergent methods to the study of $p$-adic automorphic forms, there
has been much work in the past several years towards a better
understanding of the canonical subgroup of a higher-dimensional abelian
variety.  

\paragraph In this paper we give an intrinsic geometric construction of the
level-$n$ canonical subgroup of an arbitrary truncated Barsotti-Tate group
$G$ of level $n$ defined over the valuation ring $R$ of $K$,
subject to a universal bound on its ``Hasse invariant'' that is
independent of the height and dimension of $G$.
Building on these results, our construction can be relativized and
extended to algebraic families of 
abelian varieties over $K$ with arbitrary fibral reduction types, still
under the same Hasse invariant bound; this is the subject of ongoing work
in collaboration with Brian Conrad.  (Working with truncated Barsotti-Tate
groups over $R$ effectively restricts our present focus to the case of
good reduction.) 

\paragraph
Two \label{par:why.important} important advantages of our approach are its 
generality, and the fact that our bounds are nearly as good as possible,
in the sense of~\parref{par:katz.bound}.  Abbes and 
Mokrane~\cite{abbes_mokrane:sous_groupes_canoniques} and Andreatta and
Gasbarri~\cite{andreatta_gasbarri:canonical_subgroups} have methods for
constructing level-$1$ canonical subgroups of abelian varieties, subject
to more restrictive Hasse invariant bounds than ours.
Tian~\cite{tian:canonical_subgroup} and
Fargues-Tian~\cite{fargues_tian:canonical_subgroup} extend Abbes and
Mokrane's work to construct the level-$n$ canonical subgroup of an
arbitrary level-$n$ truncated Barsotti-Tate group,
again subject to more restrictive bounds and over noetherian base rings.
Kisin and Lai~\cite{kisin_lai:hilbert_modular_forms} and
Kassaei~\cite{kassaei:p_adic_mfs} also have theories of the level-$1$
canonical subgroup for certain universal families of abelian varieties;
their methods  rely on the integral structure of certain Shimura
varieties.  Recently, in an extension of their ``subgroup-free'' methods
\cite{goren_kassaei:subgroup_free}, Goren and Kassaei have made a
construction of the canonical subgroup of any level of a certain class of
abelian varieties with real multiplication, with better bounds than ours
in terms of partial Hasse invariants; again, the advantage of our approach is
that it works for arbitrary abelian varieties (independent of any level
structure) while remaining quite explicit. 

\paragraph
We expect that our results, along with certain compatibility properties
with respect to isogenies and the Frobenius endomorphism that will appear
in a future paper, should be useful in the application of overconvergent
methods to the study of $p$-adic automorphic forms on quite general
modular varieties such as Hilbert modular varieties, Siegel modular
varieties, and more general Shimura varieties of PEL type.  We expect that
the explicit nature of our bounds will be important in such applications.
One might hope in particular to find classical modular forms in $p$-adic
families of $p$-adic modular forms as
in~\cite{kisin_lai:hilbert_modular_forms}, and in the long term perhaps to
prove a ``control theorem'' in certain settings generalizing results of
Coleman~\cite{coleman:classical_overconvergent_mfs,coleman:classical_overconvergent_mfs_higher_level}
and Kassaei~\cite{kassaei:gluing_lemma}.

Our methods are also entirely different from those used in the literature
mentioned in~\parref{par:why.important}.  They involve some fairly
explicit  calculations with Zink's displays, and they use the language of
of tropical algebraic geometry in a crucial way.  We expect that these
methods and some of the related intermediate results will prove 
useful in other contexts. 

\paragraph[Overview of results]
Let $R$ be a complete rank-$1$ valuation ring of mixed characteristic $(0,p)$, with
valuation $\ord$ normalized so that $\ord(p) = 1$ and absolute value
$|\cdot| = p^{-\ord(\cdot)}$.  We do not require $R$ to be
discretely-valued; this generality will be useful for avoiding perfectness
hypotheses on the residue field.  Let
$G$ be a $p$-divisible group of height $h$ and dimension
$g$ over $R$, and let $G^\circ \cong \Spf(R\ps{X_1,\ldots,X_g})$
be the connected component of $G$.  Let $G_0 = G\tensor_R(R/pR)$.
For $r\in R/pR$ we let $\td r$ denote any lift to $R$; observe that
$\min\{\ord(\td r), 1\}$ depends only on $r$.
Let $V:G_0^{(p)}\to G_0$ be the Verschiebung homomorphism over $R/pR$, and
let $\dV$ be the associated map on tangent spaces.  Choosing bases, we may
identify $\dV$ with a $g\times g$ matrix with entries in $R/pR$.  The 
\emph{Hasse invariant} of $G$ is defined to be
\[ H(G) = \min\{\ord(\det(\dV)^\sim),~ 1\}. \]
One can show (Remark~\ref{rem:after.hasse}(ii)) that $H(G) = 0$ if and
only if the Cartier dual of $G^\circ$ 
is \'etale, which is the case exactly when the height of $G^\circ$ is
equal to $g$.  

\paragraph
The $p$-power torsion levels $G^\circ[p^n]$ are truncated Barsotti-Tate
groups over $\Spec(R)$, in the sense of~\cite[Chapter~1]{messing:thesis}
or~\parref{par:BTn}. 
For $\xi\in G^\circ[p^n](\bar K)$ let
$|\xi| = \max\{|X_1(\xi)|,\ldots,|X_g(\xi)|\} < 1$; this is the distance
of $\xi$ from the origin, and is independent of the choice of
$R$-parameters $X_i$ for $G^\circ$.
For $\rho > 0$ let $G^\circ[p^n]_{\leq\rho}$ be the $K$-subgroup of
$G^\circ[p^n]$ whose geometric points are
\[ G^\circ[p^n]_{\leq\rho}(\bar K)
= \{\xi\in G^\circ[p^n](\bar K)~:~|\xi|\leq\rho\}, \]
and let 
$G^\circ[p^\infty]_{\leq\rho} = \bigcup_{n=1}^\infty G^\circ[p^n]_{\leq\rho}$.
Obviously $G^\circ[p^n]_{\leq 1}\cong(\Z/p^n\Z)^h$, and for small enough
$\rho$ we have $G^\circ[p^n]_{\leq\rho}(\bar K) = \{0\}$.  Also note that
if $\rho\leq\mu$ then 
$G^\circ[p^n]_{\leq\rho}\subset G^\circ[p^n]_{\leq\mu}$.
If $G^\circ[p^n]_{\leq\rho}(\bar K)\cong(\Z/p^n\Z)^g$ for some $\rho > 0$
then we call $G^\circ[p^n]_{\leq\rho}$ the 
\emph{level-$n$ canonical subgroup of $G$}, and we say that this canonical
subgroup \emph{admits the radius $\rho$}.  

\paragraph
Such a subgroup, if it exists, does not depend on the
choice of $\rho$: it is the unique subgroup with geometric group structure
$(\Z/p^n\Z)^g$ whose geometric points are closer to the origin than all
other points of $G^\circ[p^n](\bar K)$.  If $H(G) = 0$ then $h = g$, so
$G^\circ[p^n] = G^\circ[p^n]_{\leq 1}$ serves as the level-$n$ canonical
subgroup for any $n$.  Conversely, if the level-$n$ canonical subgroup
exists for all $n$ then $H(G)=0$ (Remark~\ref{rem:all.cs.ordinary}).
In general it is not true that even the level-$1$ canonical subgroup
exists, even when $h=2$ and $g=1$, although it is clear from the
definition that if there exists a level-$n$ canonical subgroup $G_n$ of
$G$ then for all $1\leq m\leq n$ the $K$-subgroup $G_n[p^m]$ is the
level-$m$ canonical subgroup of $G$.  Also note that the canonical
subgroup of $G$ depends only on the connected component $G^\circ$, and
that its existence and formation are insensitive to valued extensions of
$K$.  

\paragraph
We will see (Remark~\ref{rem:after.hasse}(i) and~\parref{par:size.pdivgp})
that the Hasse invariant and level-$n$ canonical subgroup of 
$G$ are notions that are intrinsic to the truncated Barsotti-Tate
group $G[p^n]$.  
The main goal of this paper is to prove the following result.

\newpage
\begin{thm}
  Let $G$ be a truncated Barsotti-Tate group of level $N\geq 1$ defined
  over a complete mixed characteristic $(0,p)$ valuation ring
  $R$, and let $H(G)$ be its Hasse invariant.  
  If $H(G) < (p-1)/p^N$ then the level-$N$
  canonical subgroup of $G$ exists, and it admits the radius 
  $\rho = p^{-r}$ where 
  \[ r = \frac 1{p^{N-1}(p-1)} - \frac {H(G)}{p-1}. \]
  \label{thm:mainthm}
\end{thm}

\subparagraph
The analogous statement for $p$-divisible groups $G$ over $R$ then
follows; this statement (with $R$ fixed) is slightly weaker than
Theorem~\ref{thm:mainthm} over $R$, since we do not know if a truncated
Barsotti-Tate group of level $N$ over 
$R$ is necessarily isomorphic to the $p^N$-torsion of a $p$-divisible
group over $R$ when the residue field of $R$ is not perfect.  
An elementary argument will reduce Theorem~\ref{thm:mainthm} to the
special case when $G$ is the $p^N$-torsion in a $p$-divisible group.

\paragraph Theorem~\ref{thm:mainthm} \label{par:katz.bound}
is a partial generalization of results of
Katz and Lubin~\cite{katz:padicmfs} on level-$1$ canonical subgroups of
elliptic curves (see in particular Theorem~3.10.7 of
\emph{loc.\ cit.}), as extended by
Buzzard~\cite[\S3]{buzzard:analytic_continuation} to higher levels
(see~\cite[Theorem~4.2.5]{conrad:modular_curves} for a proof that
Buzzard's higher-level canonical subgroups agree with ours).
Katz \emph{et al} prove in particular that if $E$ is an elliptic curve over $R$
and if $G = E[p^\infty]$, then the level-$N$ canonical subgroup
of $G$ exists if and only if $H(G)<1/p^{N-2}(p+1)$, with the radius as in
Theorem~\ref{thm:mainthm}.
This (sharp) bound is better than our bound by a 
factor of $p^2/(p^2-1)$; see Remark~\ref{rem:why.our.bound}.
Katz also analyzes the behavior of 
the canonical subgroup with respect to degree-$p$ isogenies of elliptic
curves, and proves that the level-$N$ canonical subgroup lifts the kernel
of the $N$-fold relative Frobenius on $G\tensor_R(R/p^\lambda R)$ for
an explicit value of $\lambda\in(0,1)$ depending on $p$, $N$, and $H(G)$.
In future work we will prove analogous results in our
situation.  

\paragraph
Here we give a brief overview of the proof of Theorem~\ref{thm:mainthm}.  
Zink~\cite{zink:display} has defined a category of semi-linear algebraic
objects over the ring of Witt vectors $W(R)$ called displays, as well as a
functor $\cP\rightsquigarrow\BT_\cP$ from the category of (nilpotent)
displays over $R$ to the category of connected $p$-divisible groups over
$R$, which is an equivalence of
categories~\cite[Theorem~1.1]{lau:displays}.  We will  
first prove Theorem~\ref{thm:mainthm} in the case when $G^\circ$ extends
to a $p$-divisible group over over $R$, in which case
$G^\circ\cong\BT_\cP[p^N]$ for some display $\cP$ over $R$, 
and then in \S\ref{sec:not.noetherian} we will reduce
Theorem~\ref{thm:mainthm} to this case.

If $\cP$ is a display over $R$, the $p$-divisible formal group $\BT_\cP$
is given in terms of its functor 
of points; \S\ref{sec:gfgp} makes explicit how to recover a
$g$-dimensional formal group law from such an object.  In
\S\S\ref{sec:params}--\ref{sec:logarithm}, we calculate
the logarithm of $\BT_\cP$ in terms of the structure coefficients of
$\cP$.  Then in 
\S\ref{sec:genericfiber} we review the fact that the kernel of the logarithm is 
the $p$-power torsion $\BT_\cP[p^\infty]$ of $\BT_\cP$ as a $K$-analytic
subspace of a $p$-adic open unit ball of dimension $g$.

By the above considerations, we need to
understand the valuations of the points of the kernel of the logarithm.
Thinking of the logarithm as a $g$-tuple of power series
$\log_1,\ldots,\log_g$ in $g$ variables, this amounts to finding the
valuations of the common zeros of the power series $\log_i$.  The theory
of tropical algebraic geometry is ideally suited for such a task; we
briefly recall the relevant parts of the theory in
\S\ref{sec:tropical}.  In 
\S\ref{sec:thepicture} we study the tropicalizations of the $\log_i$
in enough detail to justify making a rather drastic deformation of 
$\ker(\log)$ in \S\ref{sec:counting}, the end result of which is a
way of counting the number of common zeros of the $\log_i$ contained in
the ball of radius $\rho$.  In particular we show that under the Hasse
invariant bound in the above Theorem, the $K$-group
$\BT_\cP[p^\infty]_{\leq\rho}$ 
has order $p^{Ng}$ for $\rho$ as in the Theorem; in
\S\ref{sec:groupstruct} we show that  
$\BT_\cP[p^\infty]_{\leq\rho}(\bar K)\cong(\Z/p^N\Z)^g$.

\begin{rem} Our \label{rem:why.our.bound}
  methods depend on the fact that when $G$ is a
  connected $p$-divisible group over $R$ with $H(G) < (p-1)/p^N$ then 
  $G[p^\infty]_{\leq\rho} = G[p^N]_{\leq\rho}$ with
  $\rho = p^{-r}, ~r = 1/p^{N-1}(p-1) - H(G)/(p-1)$.  This is \emph{false}
  in general when $H(G) \geq (p-1)/p^N$: for instance, let $E$ be an
  elliptic curve over $R$ with ordinary reduction and let $E'$ be an
  elliptic curve with supersingular reduction and Hasse invariant
  $H(E')=(p-1)/p$.  By the work of Katz \emph{et al} in \emph{loc.\ cit.,} $E'$
  has a level-$1$ canonical subgroup with radius 
  $\rho = p^{-1/p(p-1)}$, and $E$ has a level-$2$ canonical subgroup with
  the same radius.  Therefore 
  $(E\times E')[p^\infty]_{\leq\rho}\neq(E\times E')[p]_{\leq\rho}$
  (although $(E\times E')[p]_{\leq\rho}$ is the level-$1$ canonical
  subgroup of $E\times E'$).  Thus our methods break down when 
  $H(G)\geq (p-1)/p^N$.  It is unclear to us whether this bound is sharp in
  general, or whether Katz-Lubin's bound holds in higher dimensions.
\end{rem}

\paragraph[General notation] 
In this paper, all rings are commutative, and contain a unit element unless
\label{par:notation}
noted otherwise; rings without a unit element will generally be denoted by
a calligraphic symbol.  If $R$ is a ring, $\cN(R)$ will denote its
nilradical. 
We fix a prime $p$ once and for all.

The symbol $\subset$ means subset inclusion; $\subsetneq$ means strict
inclusion. For $r\in\R$ let
\[ \R_{>r} = \{x\in\R~:~x > r\} \qquad
\R_{\geq r} = \{x\in\R~:~x\geq r\}, \]
and similarly for $\R_{<r}$, $\R_{\leq r}$,
$\Z_{\leq r}, \Z_{<r}, \Z_{>r},$ and $\Z_{\geq r}$.  We let
$\N=\Z_{\geq1}$.  

We will use other bold capital letters $\bX,\bY,\ldots$ to denote a
sequence of indeterminates, that is, $\bX = X_1,\ldots,X_g$,
$\bY = Y_1,\ldots,Y_g$, etc.

\subparagraph We \label{par:witt.notation}
use the following notation for Witt vectors.
If $R$ is any ring, 
\[ W(R) = \{ (x_0,x_1,x_2,\ldots)~:~x_i\in R\} \]
denotes the ring of $p$-Witt vectors.  For $x\in W(R)$, the 
\emph{Witt components} of $x$ will be denoted $x_0,x_1,\ldots$, i.e.,
$x = (x_0,x_1,x_2,\ldots)\in W(R)$.  For $n\in\Z_{\geq 0}$ the 
\emph{$n$th ghost component} map is the ring homomorphism $w_n: W(R)\to R$
defined by 
\[ w_n(x) = x_0^{p^n} + px_1^{p^{n-1}} + \cdots + p^{n-1}x_{n-1}^p +
p^nx_n. \]
We set 
\[ I_R = \ker(w_0) = \{(x_0,x_1,x_2,\ldots)\in W(R)~:~x_0 = 0\}, \]
and for $x\in R$ we let $[x]$ denote the Witt vector
$(x,0,0,\ldots)\in W(R)$.
The Frobenius and Verschiebung endomorphisms are denoted
${}^F(\cdot),{}^V(\cdot):W(R)\to W(R)$, respectively, and are defined by
the relations 
\[ w_n({}^F x) = w_{n+1}(x) \qquad w_n({}^V x) = 
\begin{cases}
  0 & \quad\text{ if } n = 0 \\
  p\,w_{n-1}(x) & \quad\text{ otherwise}
\end{cases} \]
and by functoriality in $R$.  More explicitly, 
${}^V(x_0,x_1,\ldots) = (0,x_0,x_1,\ldots)$ (so $I_R = {}^V(W(R))$), and
when $pR=0$ we have 
${}^F(x_0,x_1,\ldots) = (x_0^p,x_1^p,\ldots)$.  When $R$ is not an
$\F_p$-algebra then the formula for ${}^Fx$ is difficult to write
explicitly; for instance, the first two components of ${}^F x$ are given by
\[ {}^F(x_0,x_1,\ldots) = 
\left(x_0^p + px_1,~ x_1^p + p x_2 - \sum_{i=0}^{p-1}
  \binom pi p^{p-i-1} x_0^{ip} x_1^{p-i}, ~\ldots\right). \]
The reader who is unfamiliar with Witt vectors over
arbitrary rings may want to do the exercise in~\cite{lang:algebra} on this
topic (the exercise number varies by edition, but can be found under the
``Witt vectors'' entry in the index), with the caveat that Lang's
definition of the Frobenius endomorphism is different from ours.  A more
complete treatment can be found in~\S17 of~\cite{hazewinkel:formal_groups}.

\subparagraph
Let $R$ be a ring, let $G$ be an $R$-group scheme, and let $n\in\N$.  We
will denote the multiplication-by-$n$ map by $[n]_G:G\to G$, and its
kernel by $G[n]$.  If $G = \Spec(A)$ is affine then we will also write
$[n]_G:A\to A$ for the corresponding $R$-algebra homomorphism.  When no
confusion is likely we will drop the subscript and write
$[n]=[n]_G$.  If $\fG\cong\Spf(R\ps{X_1,\ldots,X_g})$ is a formal group 
over $R$ then we define $[n]=[n]_\fG$ and $\fG[n]$ likewise.
If $G = \Spec(A)$ is an affine $R$-group scheme
(resp. $G\cong\Spf(R\ps{X_1,\ldots,X_g})$ is a formal group over $R$) and
$I\subset A$ is the augmentation ideal then the 
\emph{cotangent space} to $G$ is the $R$-module $I/I^2$, and the
\emph{tangent space} is $\Lie(G) = \Hom_R(I/I^2,R)$.

\paragraph[Acknowledgements]
This work was completed as part of the author's doctoral degree program.
The author would like to express sincere gratitude to his advisor Brian
Conrad for suggesting this problem as a thesis topic and for all of his
guidance.  He is equally indebted to his
advisor Ravi Vakil for his consistent support and advice throughout
the previous five years.  The author also thanks Sam Payne and Mark Kisin
who provided key input and interesting conversations, and the referee for
pointing out some errors and inaccuracies.

\medskip
\tableofcontents

\section{Generalities concerning $p$-divisible groups}
\label{sec:pdivgps}

\paragraph
Here we fix our ideas concerning formal Lie
groups, $p$-divisible groups, and truncated Barsotti-Tate groups.  For
details see~\cite[Chapters~I and~II]{messing:thesis}, 
\cite{tate:pdivgps}, and~\cite{illusie:pdivgp.deformation}.

\paragraph
Let $R$ be a ring endowed with the discrete topology.
A \emph{formal Lie variety of dimension $g$} over $R$ is a pointed formal
scheme $\fX = \Spf(A)$ over $\Spf(R)$ such that $A\cong R\ps{\bX}$, where 
$\bX = (X_1,\ldots,X_g)$ is a set of $g$ indeterminates and $R\ps\bX$ is
given the $(\bX)$-adic topology.  We call such an isomorphism 
$A\cong R\ps{\bX}$ a 
\emph{choice of parameters} for $\fX$.  A \emph{formal Lie group} over
$R$ is a formal Lie variety $\fG$ over $R$ which is a
group object in the category of pointed formal schemes over $\Spf(R)$.
Given a formal Lie group $\fG$ and a choice of parameters 
$\fG\cong\Spf(R\ps\bX)$, the group structure on $\fG$ 
is given by a $g$-dimensional  formal group law $F$ on
$R\ps\bX$.  Note that the tangent space of a $g$-dimensional formal Lie
group $\fG\cong\Spf(R\ps\bX)$ is the rank-$g$ free $R$-module
$\Hom_R((\bX)/(\bX)^2,R)$.  In this paper, all formal group 
laws and all formal Lie groups will be assumed commutative.

\paragraph
A \emph{$p$-divisible group} or \emph{Barsotti-Tate group} over $R$ is a
directed system $G = \{G(n)\}_{n\in\N}$ of finite locally-free
commutative group schemes over $\Spec(R)$ such that for all $n\in\N$,
\begin{enum}
\item the map $G(n)\to G(n+1)$ is a closed immersion
  identifying $G(n)$ with $G(n+1)[p^n]$, and
\item the resulting map $[p]: G(n+1)\to G(n)$ is faithfully flat.
\end{enum}
If $G = \{G(n)\}_{n\in\N}$ is a $p$-divisible group then we will write
$G[p^n] = G(n)$.  Suppose that $R$ is a $p$-adically separated and
complete local ring.  For all $n\geq 1$ the
$R$-group scheme $G[p^n]$ is finite and flat%
\footnote{A finitely generated module over a local ring $R$ is flat if and
  only if it is free by~\cite[Theorem~7.10]{matsumura:crt}.  Hence 
  a finite $R$-scheme $X = \Spec(A)$ is flat if and only if $A$ is
  a free $R$-module.} 
with constant order $p^{nh}$
for some number $h\in\Z_{\geq 0}$ (not depending on $n$), called the
\emph{height} of $G$.  The \emph{connected component} of a $p$-divisible
group $G$ over $R$ is the $p$-divisible group $G^\circ = \{G[p^n]^\circ\}_{n\in\N}$, where
$G[p^n]^\circ$ is the connected component of the identity section in
$G[p^n]$.  If $G=G^\circ$ we say that $G$ is \emph{connected}, or that $G$ is a
\emph{$p$-divisible formal group}.  

\paragraph We \label{par:G.p.infty}
continue to assume $R$ is a $p$-adically separated and complete local ring.
By~\cite[Chapter~II]{messing:thesis}, a $p$-divisible formal group $G$
over $R$ is a formal Lie group in the sense that there exists a formal Lie
group $\fG$ over $\Spf(R)$, canonically determined by $G$, such that 
$\fG[p^n]$ is a finite flat group scheme over $\Spec(R)$ and
$G \cong \fG[p^\infty] \coloneq \{\fG[p^n]\}_{n\in\N}$ naturally.
The \emph{tangent space} to a $p$-divisible group $G$ is the finite free
$R$-module $\Lie(G) \coloneq \Lie(\fG)$, where $\fG$ is the formal Lie
group determined by the $p$-divisible formal group $G^\circ$ as above,
and the \emph{dimension} of $G$ is the $R$-rank of $\Lie(G)$.  

\begin{rem}
  Let $\fG\cong\Spf(R\ps\bX)$ be a formal Lie group over $R$, and suppose
  that $R$ is a complete local \emph{noetherian} ring with residue
  characteristic $p$.  By~\cite[Proposition~2.2.1]{tate:pdivgps},
  $G = \{\fG[p^n]\}_{n\in\N}$ is a $p$-divisible formal group over $R$ if
  and only if $[p]:\fG\to\fG$ is an \emph{isogeny}, i.e., 
  $[p]:R\ps\bX\to R\ps\bX$ makes $R\ps\bX$ into a finite free module over
  itself.  This condition can be checked after base change to the residue field 
  $k$, or to any field extension of $k$.
\end{rem}

\begin{rem}
  Suppose that $R$ is a local $\F_p$-algebra.  \label{rem:captures.Lie}
  Let $G$ be a $p$-divisible formal group over $R$,
  and let $\fG\cong\Spf(R\ps\bX)$ be the associated formal Lie group.  Let
  $I = (\bX)\subset R\ps\bX$ be the augmentation ideal and let
  $J\subset R\ps\bX$ be the ideal defining $\fG[p]$.  Since the relative
  Frobenius $F:\fG\to\fG^{(p)}$ factors through $[p]:\fG\to\fG$, we have
  $J\subset(X_1^p,\ldots,X_g^p)\subset I^2$.  Thus $G[p]$ ``captures''
  the tangent space $\Lie(G)$ in the sense that the canonical map
  $\Lie(G[p])\to\Lie(G)$ is an isomorphism.  It follows that 
  $\Lie(G[p^n])\to\Lie(G)$ is an isomorphism for all $n\geq 1$.

  More generally, if $p$ is nilpotent in $R$ then 
  $\Lie(G[p^n])\to\Lie(G)$ is an isomorphism for large enough $n$;
  see~\cite[\S II.3]{messing:thesis}. 
\end{rem}

\paragraph Let \label{par:BTn}
$R$ be any ring, let $G$ be a finite locally-free group scheme over
$R$, and let $n\geq 1$.  We say that $G$ is a 
\emph{truncated $p$-divisible group of level $n$}, or a
\emph{truncated Barsotti-Tate group of level $n$}, or a $\BT_n$, 
provided that $G$ is killed by $p^n$, and
for $0\leq i\leq n$ the map $p^i: G\to G[p^{n-i}]$ is faithfully flat.
If $n=1$ we also require the following condition:
\begin{deflist}
\item[] Let $G_0 = G\tensor_R (R/pR)$, and let
  $F:G_0\to G_0^{(p)}$ and $V:G_0^{(p)}\to G_0$ be the
  relative Frobenius and Verschiebung homomorphisms,
  respectively (see~\cite[Exposé~VII$_{\text{A}}$.4]{sga3} for the
  definition).  We require that $\ker(F),\ker(V)$ be finite
  locally-free $(R/pR)$-group schemes and that $F: G_0\to\ker(V)$ and 
  $V:G_0^{(p)}\to\ker(F)$ be faithfully flat.
\end{deflist}
Suppose again that $R$ is a $p$-adically separated and complete local ring
with residue field $k$.
If $G$ is a $\BT_n$ for some $n\geq 1$ and $G_k \coloneq G\tensor_R k$ then
$\Lie(G_k)$ is a finite-dimensional $k$-vector space; its $k$-dimension is
defined to be the 
\emph{dimension} of $G$.  The order of $G[p]$ is equal to $p^h$ for some
number $h\geq 0$, called the \emph{height} of $G$.
Note that if $G$ is a $\BT_n$ of dimension $g$ over $R$ and $G^\circ$ is
the connected component of the identity section in $G$ then $G^\circ$ is
also a $\BT_n$ of dimension $g$ over $R$.

\paragraph
If $G$ is a $\BT_n$ of dimension $g$ and height $h$ for $n\geq 2$
then $G[p^{n-1}]$ is a $\BT_{n-1}$ of dimension $g$ and height $h$
by~\cite[Proposition~II.3.3.11]{messing:thesis}, and 
if $G$ is a $p$-divisible group of dimension $g$ and height $h$ then
$G[p^n]$ is a $\BT_n$ of dimension $g$ and height $h$ for all $n\geq 1$.
Conversely, if $R$ is noetherian and its residue field is perfect
then by a theorem of
Grothendieck~\cite[Theorem~4.4(e)]{illusie:pdivgp.deformation}, any
$\BT_n$ over $R$ is isomorphic to the $p^n$-torsion subgroup of a
Barsotti-Tate group over $R$.

\section{Formal groups and formal group laws} \label{sec:gfgp}

\paragraph
In~\cite{zink:cartiertheorie}, Zink develops an algebraic theory of formal
groups which is well-adapted to calculations with displays.  Here we make
explicit the relation between Zink's formal groups and classical formal
group laws.  

In this section, $R$ is any ring.

\begin{defn}
  A \emph{nilpotent $R$-algebra} is an $R$-algebra $\cN$ (generally
  without unit) such that for some $r\geq 0$, a product of any $r$
  elements of $\cN$ is equal to zero.  In other words, if $\cN^r$ denotes
  the ideal generated by $\{x_1 x_2\cdots x_r~:~x_i\in\cN\}$, then 
  $\cN^r = 0$.  We let $\Nil_R$ denote the category of nilpotent
  $R$-algebras.

  A \emph{topologically nilpotent $R$-algebra} is a topological
  $R$-algebra $\cN$ such that (i) $\bigcap_{r=1}^\infty\cN^r = 0$, and (ii)
  the natural map $\cN\to\invlimm_r\cN/\cN^r$ is a topological isomorphism,
  where $\cN/\cN^r$ is given the discrete topology.
  Let $\Nil^\wedge_R$ denote the category whose objects are the
  topologically nilpotent $R$-algebras and whose morphisms are 
  $R$-algebra homomorphisms.
\end{defn}

\newpage
\begin{rem} \label{rem:after.nilr}
\begin{enum}
\item If $\cN$ is a nilpotent (resp. topologically nilpotent)
    $R$-algebra, we put  a ring structure on  
    $R\dsum\cN$ by setting $(r,n)\cdot(r',n') = (rr', rn'+r'n+nn')$.  Thus
    $\Nil_R$ is equivalent to the category of augmented $R$-algebras with
    nilpotent augmentation ideals, and similarly for $\Nil_R^\wedge$.

\item An $R$-module $M$ can be viewed as a
    nilpotent $R$-algebra by setting $M^2 = 0$; in this way we think of the
    category $\Mod_R$ of $R$-modules as a full subcategory of $\Nil_R$.  

\item The category $\Nil_R$ is naturally an abelian category.  It is a
    full subcategory of $\Nil_R^\wedge$.

\item If $\cN_1,\cN_2\in\Nil_R^\wedge$ and $f:\cN_1\to\cN_2$ is an $R$-algebra
    homomorphism then $\cN_1^r\subset f\inv(\cN_2^r)$ for all $r\geq 1$,
    so $f$ is continuous.

\item A topologically nilpotent $R$-algebra is by definition an inverse
    limit of nilpotent $R$-algebras.   
    Hence any functor $G:\Nil_R\to\Ab$ canonically extends to
    $\Nil_R^\wedge$ by setting 
    \[ G(\cN) = \invlimm_r G(\cN/\cN^r). \]

\end{enum}
\end{rem}

\begin{eg}
  Let $\bX = X_1,\ldots,X_g$ be indeterminates, and let
  $\cN$ be the ideal $(\bX)R\ps{\bX}$.  Then
  $\cN$ is a topologically nilpotent $R$-algebra, and for all $r\geq 1$ the
  quotient $\cN/\cN^r = (\bX)R\ps{\bX}/(\bX)^r$ is a nilpotent
  $R$-algebra.
\end{eg}

\paragraph
Let $\Lambda_R\in\Nil_R$ be the 
$R$-module $R$ regarded as a nilpotent $R$-algebra, so
$\Lambda_R \cong \epsilon R[\epsilon]/\epsilon^2$.  

\begin{defn}[{\cite[Definition~80]{zink:display}}] \label{def:fg}
  A \emph{finite-dimensional abelian formal group} over $R$ is a functor
  $G$ from $\Nil_R$ to the category $\Ab$ of abelian groups, satisfying: 
  \begin{enum}
  \item $G(0) = 0$,
  \item $G$ takes exact sequences of $R$-algebras to exact sequences of
    abelian groups,
  \item $G$ naturally commutes with arbitrary direct sums in $\Mod_R$, and
  \item the tangent space $t_G = G(\Lambda_R)$ is a finite free
    $R$-module.
  \end{enum}
  A \emph{morphism} of finite-dimensional abelian formal groups is a
  natural transformation of functors.
\end{defn}

\subparagraph
Property~(iv) requires some explanation.  Let $\cN$ be an $R$-algebra
such that $\cN^2 = 0$.  Multiplication by $r\in R$ induces an
$R$-algebra endomorphism of $\cN$, and hence an endomorphism of 
$G(\cN)$.  The other axioms imply that this is an $R$-module structure
on $G(\cN)$.  In particular, $t_G$ is naturally an $R$-module.  

\subparagraph
In this paper, a \emph{formal group over $R$} is a finite-dimensional
abelian formal group over $R$ in the above sense, unless specified
otherwise.  The \emph{dimension} of $G$ 
is the rank of $t_G$.  We will implicitly extend such $G$ to a functor on
$\Nil_R^\wedge$ as in Remark~\ref{rem:after.nilr}(v).

\begin{rem}
  Zink in fact only requires that $t_G$ be a (finitely generated)
  projective $R$-module.  In the sequel we will assume that $R$ is local,
  so all projective $R$-modules are free
  (\cite[Theorem~2.5]{matsumura:crt}). 
\end{rem}

\paragraph
There is a natural functor  $\cA$ from the category of formal Lie groups
over $R$ to the category of formal groups over $R$, defined as follows.  Let
$\fG = \Spf(A)$ be a formal Lie group, and let $G = \cA(\fG):\Nil_R\to\Ab$
be the functor 
\[ G(\cN) = \ker\big(\Hom_R(\Spec(R\dsum\cN), \fG)\To
\Hom_R(\Spec(R),\fG)\big) = \Hom_R(I,\cN), \]
where $R\dsum\cN$ is the augmented $R$-algebra defined
in Remark~\ref{rem:after.nilr}(i) and $I$ is the augmentation ideal of $A$.
It is easily checked that $G$ is a formal 
group over $R$, and moreover that there is a natural identification 
of $R$-modules $\Lie(\fG)\cong t_G$.
Given a choice of parameters $\bX = X_1,\ldots,X_g$ 
for $A \cong R\ps\bX$, let $F$ be the formal group law on $R\ps\bX$
defining the group structure on $\fG$.  Then an element of $G(\cN)$ is a
homomorphism $R\ps\bX\to R\dsum\cN$ such that the image of
each $X_i$ is contained in $\cN$; hence $G(\cN)\cong\cN^{\dsum g}$, with 
$G(\cN)$ imposing the law of composition $(\bx,\by)\mapsto F(\bx,\by)$ on
$\cN^{\dsum g}$ via this bijection.

\begin{rem}
  Let $G$ be a $p$-divisible formal group over a local ring $R$ that is
  $p$-adically separated and complete, and let $\fG$ be the canonically
  associated formal Lie group such that $G \cong \fG[p^\infty]$ as
  in~\parref{par:G.p.infty}.  Then $\cA(\fG)$ is 
  a formal group which is naturally associated with $G$; in later sections we will
  identify $G$ with the formal group $\cA(\fG)$ without mentioning $\fG$.
\end{rem}

\paragraph By~\cite[II.2.32]{zink:cartiertheorie}, \label{par:fg.from.fglaw}
$\fG\rightsquigarrow\cA(\fG)$ is an equivalence of categories.  In other words,
if $G$ is a $g$-dimensional formal group then there is a formal group law
$F$ on $R\ps\bX$ and a functorial isomorphism of abelian groups
$G(\cN)\cong\cN^{\dsum g}$, where the group law on $\cN^{\dsum g}$ is
defined by $F$ as above.  
In the sequel we will be given an explicit description of an (abstract)
formal group $G$ for which we will want to recover such a formal group
law.  The remainder of this section is devoted to the 
construction of a (non-canonical) formal Lie group $\fG$ such that
$G\cong\cA(\fG)$ for a formal group $G$.  

\paragraph To \label{par:fg.from.fglaw2}
motivate this construction, suppose that $G = \cA(\fG)$ for some formal 
Lie group $\fG = \Spf(A)$ over $R$.  We can recover $\fG$ from $G$ in the
following way.  Let $I\subset A$ be the augmentation ideal, let 
$I'= \prod_{i=1}^\infty\Sym^i(t_G^*)$, and let $A' = R\dsum I'$.
We have a natural isomorphism $t_G^* \coloneq \Hom_R(t_G,R)\cong I/I^2$,
and hence a (non-canonical) isomorphism of $R$-algebras $\alpha: I'\isom I$.  
The choice of isomorphism $\alpha$ gives rise to functorial isomorphisms
\[ \Hom_R(I', \cN) \isom \Hom_R(I,\cN) = G(\cN), \]
which recovers the structure of formal Lie group on
$\Spf(A')\cong\fG$.  
Choosing a basis for $t_G$ then gives a choice of parameters
$A'\cong R\ps\bX$, giving rise to a functorial 
isomorphism of pointed sets $G(\cN)\cong\cN^{\dsum g}$.  The formal group
law $F$ defining the group structure on $\fG$ can then be recovered by
substituting $\cN = (\bX,\bY)R\ps{\bX,\bY}$ and calculating
$F(\bX,\bY) = \bX +_{G(\cN)} \bY \in\cN^{\dsum g}$.

\begin{defn} 
  Let $G$ be a $g$-dimensional formal group, and let  
  \label{defn:choice_of_params}
  $I = \prod_{i=1}^\infty\Sym^i(t_G^*)$.
  A \emph{choice of Hopf algebra} for $G$ is
  functorial isomorphism of pointed sets 
  \[ \Hom_{\Nil_R^\wedge}(I,\cN)\isom G(\cN), \]
  and \emph{choice of parameters} for $G$ is a functorial isomorphism of
  pointed sets
  \[ \cN^{\dsum g}\isom G(\cN). \] 
\end{defn}

\paragraph
Let $I$ be a choice of Hopf algebra for a formal group $G$, 
let $A = R\dsum I$, and let $\fG=\Spf(A)$.  Then $\fG$ has a structure of
formal Lie group such that $G\cong\cA(\fG)$, as
in~\parref{par:fg.from.fglaw2}. 
A choice of parameters for $G$ is equivalent to a choice of Hopf algebra
along with a choice of basis for $t_G$, which is in turn equivalent to a
choice of parameters for $\fG$.  Conversely, let $\fG = \Spf(A)$ be a
formal Lie group, and let $I\subset A$ be the augmentation ideal.
An isomorphism 
$I\cong\prod_{i=1}^\infty\Sym^i(I/I^2)$ gives rise to a choice
of Hopf algebra for $G = \cA(\fG)$; this  along with a choice of basis
for $I/I^2$ gives rise to a choice of parameters for $G$.
Since $\cA$ is an equivalence of categories~\parref{par:fg.from.fglaw},
every formal group $G$ admits a choice of Hopf algebra and a choice of
parameters. 

\begin{eg}
  If $\fG = \Spf(R\ps\bX)$ then a
  choice of parameters for $\cA(\fG)$ is equivalent to a choice of elements
  $\bY = Y_1,\ldots,Y_g\in (\bX)R\ps\bX$ inducing an isomorphism
  $R\ps\bY\isom R\ps\bX$.  
\end{eg}

\paragraph
By Yoneda's lemma, any natural transformation of pointed set-valued functors
$\Hom_{\Nil_R^\wedge}(I,\cdot)\to G$ is of the form
$f\mapsto G(f)(\alpha)$ for a unique element
$\alpha\in G(I)$.  We will use the following criterion for a particular
$\alpha\in G(I)$ to determine a choice of Hopf algebra.  

\begin{lem} 
  With the notation in Definition~\ref{defn:choice_of_params}, choose
  $\alpha\in G(I)$, and let $\Phi:\Hom_{\Nil_R^\wedge}(I,\cdot)\to G$ be the
  natural transformation $f\mapsto G(f)(\alpha)$.  Then $\Phi$ is a choice
  of Hopf algebra for $G$ if and only if the composite map 
  \[ \theta:~ t_G \cong \Hom_{\Mod_R}(t_G^*,R)
  = \Hom_{\Nil^\wedge_R}(I,R) \overset\Phi\To G(\Lambda_R) = t_G \]
  is bijective.
  \label{lem:alphacriterion}
\end{lem}

\pf
  Since $G$ admits a choice of parameters, we may assume that 
  $G(\cN) = \Hom_{\Nil_R^\wedge}((\bX),\cN)$, where
  $(\bX) = (X_1,\ldots,X_g)\subset R\ps{X_1,\ldots,X_g}$.
  Using the identification of $(\bX)/(\bX)^2$ with $t_G^*$, we also
  identify $I$ with $(\bX)$.  
  Then $\alpha\in G(I) = \Hom_{\Nil_R^\wedge}((\bX),(\bX))$, and 
  $\theta$ is a bijection if and only if the $R$-module endomorphism
  $\alpha_1$ of $(\bX)/(\bX)^2$ induced by $\alpha$ is bijective.  
  This is equivalent to $\alpha$ being an $R$-algebra automorphism of
  $(\bX)$. 
\qed

\subparagraph
Let $\alpha_1\in G(t_G^*)$ be the image of $\alpha$ under the
quotient map $I\to t_G^*$, so $f\mapsto G(f)(\alpha_1)$ is a
functorial isomorphism of $R$-algebras 
$\Hom_R(t_G^*,\cN)\isom G(\cN)$ when $\cN^2=0$.  The map $\theta$ in
Lemma~\ref{lem:alphacriterion} only depends on $\alpha_1$, so any
$\alpha\in I$ lifting $\alpha_1$ gives rise to a choice of Hopf algebra
for $G$. 

\paragraphnn
For completeness we remark that the functor $\cA$ respects extension of
scalars, in the following sense: 

\begin{defn} 
  Let $\phi: R\to R'$ be a ring homomorphism, and let
  $G$ be a formal group over $R$.  Define
  the \emph{base change $G_{R'}$ of $G$ to $R'$} to be the functor on
  $\Nil_{R'}$ defined by $G_{R'}(\mc N) = G(\mc N)$, where we view a
  nilpotent $R'$-algebra $\mc N$ as a nilpotent $R$-algebra via $\phi$.
  \label{defn:fgbasechange}
\end{defn}

\subparagraph
Note that $G_{R'}$ is likewise a formal group over $R'$.  
If $G\cong\cA(\fG)$ for a formal Lie group $\fG\cong\Spf(R\ps\bX)$ over
$R$ with formal group law $F$, then 
$\fG_{R'} \coloneq \fG\hat\tensor_R R'$ is a formal Lie group over $R'$
with formal group law $\phi(F)$, and $G_{R'}\cong\cA(\fG_{R'})$.

\section{A choice of parameters from a display}\label{sec:params}

\paragraph[Displays] In  \label{par:displays.first}
this section, $R$ is a $p$-adically complete and separated local ring.
Recall from \parref{par:witt.notation} that $W(R)$ denotes the ring of
$p$-Witt vectors over $R$ and that $I_R = \ker(w_0)\subset W(R)$.
A \emph{display} over $R$ consists of the data
$\mc P = (P,Q,F,V\inv)$, where $P$ is a finite projective $W(R)$-module,
$Q\subset P$ is a $W(R)$-submodule, and $F:P\to P$, $V\inv:Q\to P$ are
${}^F$-linear maps, satisfying the following conditions:
\begin{enum}
\item We can write $P$ 
  as a direct sum of projective $W(R)$-modules $P = T\dsum L$ such that 
  $Q = I_R T\dsum L$.
\item The $W(R)$-linear map $W(R)\tensor_{W(R),{}^F} Q\to P$ induced by
  $V\inv$ is surjective. 
\item For $x\in P$ and $w\in W(R)$, we have
  $V\inv(({}^Vw) x) = w F x$.
\end{enum}
It follows from~(iii) that
\begin{equation} \label{eq:VF.p}
  Fy = V\inv({}^V 1\, y) = {}^{FV} 1\,V\inv y = p\,V\inv y
\end{equation}
for $y\in Q$.
A decomposition as in (i) is called a \emph{normal decomposition}%
\footnote{Be aware that Zink~\cite{zink:display} writes a normal
  decomposition as $P = L\dsum T$, yet always chooses a basis
  $e_1,\ldots,e_h$ for $P$ such that
  $e_1,\ldots,e_g$ is a basis for $T$ and $e_{g+1},\ldots,e_h$ is a basis
  for $L$.} 
of $P$.  The $W(R)$-rank of $P$ is called the \emph{height} of $\cP$, and the
$R$-rank of $P/Q$ (or equivalently the $W(R)$-rank of $T$) is its
\emph{dimension}.  There is an evident notion of morphism of displays, so
we can speak of the category of displays over $R$.  

\subparagraph
A \emph{nilpotent display} is a display that satisfies an additional
nilpotence condition; see~\cite[Definitions~1,~11,~13]{zink:display}.
Nilpotent displays correspond to connected $p$-divisible groups, so all
displays in this paper will be assumed to be nilpotent.

\paragraph By \label{par:W(R)local} 
\cite[Proposition~3]{zink:display}, $W(R)$ is complete and separated in
the $I_R$-adic topology.  Hence $I_R$ is contained in the Jacobson
radical of $W(R)$, so the maximal ideals of $W(R)$ are the same as the
maximal ideals of $W(R)/I_R = R$.  Therefore $W(R)$ is local ring, so
any projective $W(R)$-module is
free~\cite[Theorem~2.5]{matsumura:crt}. This allows us to choose bases  
for our displays.

\paragraph
If $R$ is a perfect ring of characteristic $p$ then (nilpotent) displays
over $R$ and 
(covariant) Dieudonn\'e modules over $R$ such that $V$ is topologically
nilpotent are equivalent notions~\cite[Proposition~15]{zink:display}.
In fact, if $\cP=(P,Q,F,V\inv)$ is a display over such a ring then $V\inv$
has an ${}^{F\inv}$-linear inverse $V:P\to Q$.  When $R$ is not a perfect
ring of characteristic $p$ then ${}^F: W(R)\to W(R)$ is not an
automorphism, and an operator $V$ does not exist in general.  Nonetheless
it is helpful when doing calculations to 
think of displays as Dieudonn\'e modules and $V\inv$ as the inverse
of $V$. 

\begin{eg} The \label{eg:display.gm}
  datum $\cG_m = (W(R), {}^V W(R), {}^F, {}^{V\inv})$ forms a display of
  height $1$ and dimension $1$.  See \cite[Example~16]{zink:display}, as
  well as Examples~\ref{eg:display.gm.2} and~\ref{eg:display.gm.3}.
\end{eg}

\paragraph
Let $\cP$ be a display over $R$ and let $R\to R'$ be a
homomorphism of $p$-adically separated and complete local rings.  The display  
$\cP_{R'} = (P_{R'},Q_{R'},F_{R'},V\inv_{R'})$ obtained from $\cP$ by
base change is defined as follows (see~\cite[Definition~20]{zink:display}): 
\begin{bullets}
\item $P_{R'} = W(R')\tensor_{W(R)} P$
\item $Q_{R'} = \ker[w_0\tensor\text{(proj)}: W(R')\tensor_{W(R)} P\to R'\tensor_R (P/Q)]
  \\{}\quad= \im[(W(R')\tensor_{W(R)} Q)\dsum(I_{R'}\tensor_{W(R)} P)\to P_{R'}]$
\item $F_{R'} = {}^F(\cdot)\tensor F:P_{R'}\to P_{R'}$
\item $V\inv_{R'}:Q_{R'}\to P_{R'}$ is the unique ${}^F$-linear
  homomorphism satisfying
  \begin{equation}
    \label{eq:Vinv.basechange}
    V\inv_{R'} (w\tensor y) = {}^F w\tensor V\inv y \qquad
    V\inv_{R'} ({}^V w\tensor x) = w\tensor Fx 
  \end{equation}
  for all $w\in W(R')$, $y\in Q$, and $x\in P$
  (recall~\parref{par:witt.notation} that ${}^V(W(R)) = \ker(w_0)$).
\end{bullets}
This sets up a functor from the category of displays over $R$ to the
category of displays over $R'$.  Note that if $P = T\dsum L$ is a normal
decomposition then 
\[ Q_{R'} = I_{R'} (W(R') \tensor_{W(R)} T)\dsum (W(R')\tensor_{W(R)} L). \]

\begin{defn}
  If $A$ is any ring and $\mc N\in\Nil_A$, we let $\hat W(\mc N)$ denote the
  $W(A)$-algebra of finite-length Witt vectors: that is,
  \[ \hat W(\mc N) = \{(x_0,x_1,\ldots)~:~x_n\in\mc N,~
  x_n = 0 \text{ for large enough } n\}. \]
  We will view $\hat W$ as an infinite-dimensional formal group over $A$
  under addition of Witt vectors, as in
  \cite[\S17.1.8]{hazewinkel:formal_groups}. 
  If $\cN$ is a nilpotent $A$-algebra and $x\in\cN$, we denote by $[x]$ the
  Witt vector $(x,0,0,\ldots)\in\hat W(\cN)$.
\end{defn}

\paragraph
There is a functor
$\BT:\cP\rightsquigarrow\BT_\cP$ from the category of (nilpotent) displays
over a ring $R$ as in~\parref{par:displays.first} to
the category of $p$-divisible formal groups over
$R$, which is an equivalence of categories~\cite[Theorem~1.1]{lau:displays}.  
This functor is compatible with change of base and respects the
notions of height and dimension.  The $p$-divisible formal group $\BT_\cP$
is given in terms of its corresponding formal group; it is constructed as 
follows.  Suppose for the moment that $p^nR = 0$ for some $n\geq 1$.
Let $\cP = (P,Q,F,V\inv)$ be a display over $R$, let $\cN$ be a nilpotent
$R$-algebra, and define 
\[\begin{split}
  \hat P_\cN &= \hat W(\cN)\tensor_{W(R)}P \\
  \hat Q_\cN &= \ker\big[w_0\tensor\text{(proj)}: 
  \hat W(\cN)\tensor_{W(R)} P \To \cN\tensor_R(P/Q)\big].
\end{split}\]
Let $R' = R\dsum\cN$ be the augmented $R$-algebra as defined
in Remark~\ref{rem:after.nilr}(i).  We can regard 
$\hat Q_\cN$ (resp. $\hat P_\cN$) as a submodule of the base change
$Q_{R'}$ (resp. $P_{R'}$).  As explained in~\cite[\S3.1]{zink:display},
the map  $V_{R'}\inv: Q_{R'}\to P_{R'}$ restricts to a map 
$V_\cN\inv:\hat Q_\cN\to\hat P_\cN$.
We define $\BT_\cP(\cN)$ to be the cokernel of the map
$V_\cN-\Id:\hat Q_\cN\to\hat P_\cN$, where $\Id:\hat Q_\cN\to\hat P_\cN$
is the natural inclusion.  By~\cite[Corollary~84]{zink:display}
the sequence of abelian groups
\begin{equation} \label{eq:BT.P.defn}
  0 \To \hat Q_\cN \xrightarrow{V_\cN\inv-\Id} \hat P_\cN\To
  \BT_{\cP}(\cN)\To 0
\end{equation}
is exact on the left as well.
We will often write $V\inv$ for $V_\cN\inv$ when it is unlikely to cause
confusion. 

Now suppose that  $p$ is not nilpotent in $R$.  Let $\cP$ be a display
over $R$,  let $R_n = R/p^{n+1}R$, let $\cP_n$ be the
base change of $\cP$ to $R_n$, and 
let $\BT_{\cP_n}$ be the associated $p$-divisible formal group.  Then 
$\{\BT_{\cP_n}\}_{n=1}^\infty$ is a compatible inverse system of
$p$-divisible formal groups, so by~\cite[Lemma~II.4.16]{messing:thesis}
there is a unique $p$-divisible formal group over $R$ whose base
change to $R_n$ is isomorphic to $\BT_{\cP_n}$ for all $n$.  We define
$\BT_\cP$ to be this $p$-divisible formal group.  Concretely, given a
choice of compatible systems of parameters $\bX = X_1,\ldots,X_g$ for
the formal groups $\BT_{\cP_n}$ --- i.e., a choice of formal group laws $F_n$
on $R_n\ps\bX$ giving rise to $\BT_{\cP_n}$ and such that 
$F_n = F_{n+1}\pmod{p^{n+1}}$ for all $n\geq 0$ --- the $p$-divisible formal
group $\BT_\cP$ is given by the formal group law $F = \lim F_n$ on
$R\ps\bX$.

\begin{eg} Suppose \label{eg:display.gm.2}
  that $p$ is nilpotent in $R$.
  Let $\cG_m$ be the display from Example~\ref{eg:display.gm}, and let 
  $G = \BT_{\cG_m}$.  If $\cN$ is a nilpotent $R$-algebra then
  \[ G(\cN) = \hat W(\cN)/(V\inv-\Id)\,{}^V(\hat W(\cN)). \]
  We claim that $G(\cN)\cong(1+\cN)^\times\subset(R\dsum\cN)^\times$.  
  The Artin-Hasse exponential~\cite[\S17.5]{hazewinkel:formal_groups} is
  the power series 
  \[ \hexp(X) = \exp\left(X + \frac{X^p}p + \frac{X^{p^2}}{p^2} +
    \frac{X^{p^3}}{p^3} + \cdots\right)\in 1+\Z_{(p)}\ps X. \]
  Define a map $h: \hat W(\cN)\to(1+\cN)^\times$ by
  $h(x_0,x_1,\ldots) = \prod_{i=0}^\infty \hexp(x_i)$.
  Formally we have
  \[ h(x_0,x_1,\ldots) =
  \exp\left(\sum_{i=0}^\infty\frac{w_n(x_0,x_1,\ldots)}{p^n}\right), \]
  so $h$ is a group homomorphism.  Then $h(V\inv x) = h(x)$ for
  $x\in{}^V\hat W(\cN)$ since 
  $V\inv(0,x_0,x_1,\ldots) = (x_0,x_1,\ldots)$, so $h$ descends to a map
  $G(\cN)\to(1+\cN)^\times$.  Since $\hexp$ is an 
  invertible power series, the map
  $y\mapsto h([y]) = \hexp(y): \cN\to(1+\cN)^\times$ is bijective, so by
  Proposition~\ref{prop:paramsfordisplays1} below, $h$ is an isomorphism.
  This shows that $G\cong\hat\bG_m =\MU_{p^\infty}$.
\end{eg}

\begin{prop} Suppose \label{prop:paramsfordisplays1} 
  that $p$ is nilpotent in $R$. Let 
  $\cP = (P,Q,F,V\inv)$ be a display over $R$, let $P = T\dsum L$ be a
  normal decomposition, and let $e_1,\ldots,e_g$ be a basis for $T$.
  The map
  \begin{equation}\label{eq:params2}
    (x_1,x_2,\ldots,x_g) \mapsto 
    [x_1]\tensor e_1 + [x_2]\tensor e_2 + \cdots + [x_g]\tensor e_g 
    \mod{(V_\cN\inv-\Id)\hat Q_{\mc N}} 
  \end{equation}
  a functorial bijection of pointed sets 
  $\mc N^{\dsum g}\isom\BT_\cP(\mc N)$; i.e., it is a choice of parameters for
  $\BT_\cP$. 
\end{prop}

\subparagraph
To prove Proposition~\ref{prop:paramsfordisplays1}, we will make use of Zink's
$\exp$ map, which is constructed as follows (still under the assumption
that $p$ is nilpotent in $R$).  Let $M$ be an $R$-module, thought of
as a square-zero $R$-algebra, and let $R' = R\dsum M$, so $M$ is
identified with the augmentation ideal in $R'$.  
As in~\cite[\S1.4]{zink:display} (or as an exercise) one can show that if
$w = (w_0,w_1,\ldots)\in\hat W(M)$ then ${}^F w = p(w_1,w_2,\ldots)$.
By~\eqref{eq:VF.p}
and~\eqref{eq:Vinv.basechange}, for $w\in\hat W(M)$ and $y\in Q$ we have 
\[ V_M\inv(w\tensor y) = {}^Fw\tensor V\inv y
= (w_1,w_2,\ldots) \tensor p\cdot V\inv y 
= (w_1,w_2,\ldots) \tensor Fy, \]
and for $x\in P$, 
\[ V_M\inv({}^V w\tensor x)
= V_M\inv((0,w_0,w_1,\ldots)\tensor x)
= (w_0,w_1,\ldots)\tensor Fx. \]
Hence it is natural to extend $V_M\inv$ to an endomorphism of 
$\hat P_M = \hat W(M)\tensor_{W(R)} P$ by the formula
\begin{equation} \label{eq:vinv}
  V_M\inv((w_0,w_1,\ldots)\tensor x) = (w_1,w_2,\ldots)\tensor Fx 
\end{equation}
(compare Lemma~38 and the proof of Lemma~83 in~\cite{zink:display}).  
Define a functorial homomorphism
\[ \exp: M\tensor_R (P/Q) \To \BT_\cP(M) \]
via the commutative diagram of exact sequences
\[\xymatrix @C+1cm{
  0 \ar[r] &
  {\hat Q_{M}} \ar@{=}[d] \ar[r] &
  {\hat W(M)\tensor_{W(R)} P} \ar[d]_{V_M\inv-\Id}
  \ar[r]^{w_0\tensor\text{(proj)}} & 
  {M\tensor_R(P/Q)} \ar[r] \ar[d]_\exp & 0 \\
  0 \ar[r] &
  {\hat Q_{M}} \ar[r]^(.3){V_M\inv-\Id} &
  {\hat W(M)\tensor_{W(R)} P} \ar[r] &
  {\BT_\cP(M)} \ar[r] & {0.}
}\]
Zink~\cite[Proof of Theorem~81]{zink:display} shows that $\exp$ is 
an $R$-linear isomorphism.

\subparagraph
\pf[of Proposition~\ref{prop:paramsfordisplays1}]
  Let $G = \BT_\cP$ (so $t_G = \BT_\cP(\Lambda_R)$), and let
  $I = \prod_{i=1}^\infty\Sym^i(t_G^*)$.  Let $\bar e_i\in P/Q$ be the
  residue of $e_i$, so $\bar e_1,\ldots,\bar e_g$ is an $R$-basis for $P/Q$, and hence
  $\exp(\bar e_1),\ldots,\exp(\bar e_g)$ form a basis of $t_{G}$.
  Let $\epsilon_1,\ldots,\epsilon_g\in t_{G}^*$ be the dual basis, and let
  \[ \bar\alpha_1
  = \epsilon_1\tensor\bar e_1 + \epsilon_2\tensor\bar e_2 +
  \cdots + \epsilon_g\tensor\bar e_g\in t_{G}^*\tensor_{R}(P/Q). \]
  Then for $x = \sum_{i=1}^g a_i\bar e_i\in P/Q$, if
  $\ev_x:t_{G}^*\to R$ denotes the evaluation map at $\exp(x)\in t_G$, we
  have 
  \[ (\ev_x\tensor\Id)\left(\sum_{i=1}^g \epsilon_i\tensor \bar e_i\right)
  = \sum_{i=1}^g \epsilon_i(\exp(x))\, \bar e_i = \sum_{i=1}^g a_i\,\bar e_i = x. \]
  Let $\alpha_1 = \exp(\bar\alpha_1)\in G(t_{G}^*)$.  By
  Lemma~\ref{lem:alphacriterion} and the above calculation, any lift of
  $\alpha_1$ to $G(I)$ determines a choice of Hopf algebra for $G$.

  Next we calculate a lift $\td\alpha_1$ of $\alpha_1$ to
  $\hat W(t_{G}^*)\tensor_{W(R)} P$ using the commutative square 
  \[\xymatrix{
    {\hat W(t_{G}^*)\tensor_{W(R)} P} \ar[d]_{V\inv-\Id}
    \ar[r] & 
    {t_{G}^*\tensor_{R}(P/Q)} \ar[d]_\exp^\cong \\
    {\hat W(t_{G}^*)\tensor_{W(R)} P} \ar[r] &
    {G(t_{G}^*)} 
  }\]
  Let
  \[ \td\alpha_1
  = -\big([\epsilon_1]\tensor e_1 + [\epsilon_2]\tensor e_2 +
  \cdots + [\epsilon_g]\tensor e_g\big)
  \in \hat W(t_{G}^*)\tensor_{W(R)} P. \]
  By construction, $\td\alpha_1$ lifts $-\bar\alpha_1$ under
  the quotient map 
  $\hat W(t_{G}^*)\tensor_{W(R)} P\surject
  t_{G}^*\tensor_{R}(P/Q)$, so 
  $(V\inv-\Id)(-\td\alpha_1)$ lifts $\alpha_1$.
  But $V\inv((\epsilon_i,0,0,\ldots)\tensor e_i) = 0$ by \eqref{eq:vinv},
  so $V\inv\td\alpha_1 = 0$.  Therefore, 
  $\td\alpha_1 = (V\inv-\Id)(-\td\alpha_1)$ lifts
  $\alpha_1$.

  Using Remark~\ref{rem:after.nilr}(v) to extend $G$ and $\hat W$ to
  functors on the category $\Nil_R^\wedge$, we have
  \[ \hat W(I) = \{(x_0,x_1,x_2,\ldots)\in W(I)~:~x_n\to 0\text{ as }n\to\infty\} \]
  and a natural map $\phi:\hat W(I)\tensor_{W(R)} P\to G(I)$.
  Using the inclusion $t_G^*\inject I$ to view $\hat W(t_G^*)$
  (resp. $G(t_G^*)$) as a subgroup of $\hat W(I)$ (resp. $G(I)$), if
  \[ \td\alpha = -([\epsilon_1]\tensor
  e_1+\cdots+[\epsilon_g]\tensor e_g)
  \in\td W(I)\tensor_{W(R)} P \]
  we see that $\alpha\coloneq\phi(\td\alpha)\in G(I)$ lifts
  $\alpha_1$ because $\td\alpha$ lifts $\td\alpha_1$.  Therefore,
  $\alpha$ determines a choice of Hopf algebra for $G$, so by the
  results of \S\ref{sec:gfgp}, for any nilpotent
  $R$-algebra $\cN$ the map
  \[ f\mapsto G(f)(\alpha):~ 
  \Hom_{\Nil_R^\wedge}(I,\cN)\To G(\cN) \]
  is a bijection of pointed sets.  Since $G(f)(\alpha)$ is the image of
  \[ (\td W(f)\tensor\Id)(\td\alpha)
  = (\td W(f)\tensor\Id)\left(-\sum_{i=1}^g[\epsilon_i]\tensor e_i\right)
  = -\sum_{i=1}^g [f(\epsilon_i)]\tensor e_i\in \hat W(\mc N)\tensor_{W(R)} P
  \]  
  in $G(\cN)$, we have that the map
  \[ (x_1,\ldots,x_g)\mapsto \sum_{i=1}^g [x_i]\tensor e_i
  \mod{(V\inv_\cN-\Id)\hat Q_{\cN}}:~
  \cN^{\dsum g}\isom G(\cN) \]
  is a bijection.
\qed

\paragraph
Now suppose that $p$ is not necessarily nilpotent in $R$,
and let $\cP = (P,Q,F,V\inv)$ be a display over $R$.
We would like an analogue of Proposition~\ref{prop:paramsfordisplays1} for
$\cP$, which is not generally true as stated in this context since the
associated $p$-divisible formal group is not constructed in the same way.  

\begin{rem}
  Let $\cN$ be a nilpotent $R$-algebra such that the natural map
  $\cN\to\invlimm_n\cN/p^n\cN$ is an isomorphism, and define
  \[ \td W(\cN) = \invlimm_n \hat W(\cN/p^n\cN)
  =  \{ (x_0,x_1,x_2,\ldots)\in W(\cN)~:~x_i \to 0 \}, \]
  where the convergence is taken in the $p$-adic topology.
  Let $\cP = (P,Q,F,V\inv)$ be a display over $R$, and set
  \[\begin{split}
    \td P_\cN &= \td W(\cN)\tensor_{W(R)}P \\
    \td Q_\cN &= \ker\big[w_0\tensor\text{(proj)}: 
    \td W(\cN)\tensor_{W(R)} P \To \cN\tensor_R(P/Q)\big].
  \end{split}\]
  With some work one can show that $\BT_\cP(\cN)$ is naturally isomorphic
  to $\td P_\cN/(V\inv-\Id)\td Q_\cN$.  With a great deal of work one can
  even prove an analogous formula for $\cN[p\inv]$ in place of $\cN$ when
  $\cN$ has no $p$-torsion.  As we do not need these refinements, we will
  not say more about them here.
\end{rem}

\paragraph
Let $P = T\dsum L$ be a normal decomposition, let 
$e_1,\ldots,e_g$ be a $W(R)$-basis for $T$, and let $e_{g+1},\ldots,e_h$
be a $W(R)$-basis for $L$.  Let $\cN$ be a nilpotent $R$-algebra such that 
$p^n\cN = 0$ for some $n$.  Our choice of basis allows us to identify
$\hat P_\cN$ with $\hat W(\cN)^{\dsum h}$ and 
$\hat Q_\cN$ with the subgroup 
${}^V(\hat W(\cN))^{\dsum g}\dsum\hat W(\cN)^{\dsum(h-g)}$; under these
identifications the exact sequence~\eqref{eq:BT.P.defn} becomes an exact
sqequence 
\[ 0\To {}^V\hat W(\cN)^{\dsum g}\dsum\hat W(\cN)^{\dsum(h-g)}
\xrightarrow{V_\cN\inv-\Id} \hat W(\cN)^{\dsum h}\To \BT_{\cP}(\cN)\To 0 \]
which is functorial in $\cN$.  It is not
hard to show that the above sequence uniquely extends to
a diagram of homomorphisms of formal groups over $R$ (i.e. of functors
$\Nil_R\to\Ab$) 
\begin{equation} \label{eq:seqkerlog}
  {}^V\hat W^{\dsum g}\dsum\hat W^{\dsum(h-g)}
  \xrightarrow{V\inv-\Id}\hat W^{\dsum h}\overset\pi\To\BT_{\cP} 
\end{equation}
whose composite is zero.

\begin{thm} Let \label{thm:paramsfordisplays}
  $R$ be a $p$-adically complete and separated local ring.
  Let $\cP = (P,Q,F,V\inv)$ be a display over $R$, let $P = T\dsum L$ be a
  normal decomposition, let $e_1,\ldots,e_g$ be a $W(R)$-basis for $T$,
  and let $e_{g+1},\ldots,e_h$ be a $W(R)$-basis for $L$.  Let 
  $\pi:\hat W^{\dsum h}\to\BT_{\cP}$ be the homomorphism
  \eqref{eq:seqkerlog} determined by our 
  choice of basis.  For a nilpotent $R$-algebra
  $\cN$, the map 
  \begin{equation}\label{eq:paramsfordisplays} 
    (x_1,\ldots,x_g) \mapsto \pi([x_1],\ldots,[x_g],0,\ldots,0):~
    \cN^{\dsum g} \To \BT_\cP(\cN) 
  \end{equation}
  is a choice of parameters for $\BT_\cP$ in the sense of 
  Definition~\ref{defn:choice_of_params}.
\end{thm}

\pf
  For $n\geq 0$ let $R_n = R/p^{n+1}R$ and let $\cP_n$ be the base change
  to $R_n$.  When $p^{n+1}\cN = 0$ for some $n$ (i.e., $\cN$ is an
  $R_n$-algebra), the map~\eqref{eq:paramsfordisplays} agrees
  with the bijection~\eqref{eq:params2} after 
  identifying $\BT_\cP(\cN)$ with $\BT_{\cP_n}(\cN)$.  If $G$ is any
  formal group over $R$ then a choice of parameters for $G$ is equivalent
  to a compatible system of choices of parameters for the formal groups
  $G_{R_n}$, so our basis does give rise to a choice of parameters for
  $\BT_\cP$.  Furthermore, a functorial map $\cN^{\dsum g}\to G(\cN)$
  is determined by its behavior on nilpotent algebras $\cN$ such that
  $p^{n+1}\cN=0$ for some $n$ (in fact one only needs to consider $\cN$ of the
  form $(X_1,\ldots,X_g) R_n\ps{X_1,\ldots,X_g}/(X_1,\ldots,X_g)^m$),
  so~\eqref{eq:paramsfordisplays} agrees with the
  choice of parameters for $\BT_\cP$ induced by our basis.
\qed

\paragraph
In choosing a basis for a display $\cP$ over $R$, we will implicitly make
the corresponding choice of parameters for $\BT_\cP$ given by
Theorem~\ref{thm:paramsfordisplays}.

\section{The logarithm of a display} \label{sec:logarithm}

\paragraph
In this section $R$ is a complete rank-$1$ valuation ring of mixed
characteristic $(0,p)$ with field of fractions $K$.
Fix a display $\cP = (P,Q,F,V\inv)$ over
$R$, a normal decomposition $P = T\dsum L$ (so $Q = I_R T\dsum L$), and $W(R)$-bases
$e_1,\ldots,e_g$ for $T$ and $e_{g+1},\ldots,e_h$ for $L$.  Let 
$G = \BT_{\mc P}$ be the formal
group over $R$ associated with $\mc P$.  With respect to our basis
the ${}^F$-linear maps $F:P\to P$ and $V\inv:Q\to P$ are determined by formulas
\begin{equation} \label{eq:alpha.ij} \begin{split}
  Fe_j &= \sum_{i=1}^h \alpha_{ij} e_i\quad j = 1,\ldots,g \\
  V\inv e_j &= \sum_{i=1}^h \alpha_{ij} e_i\quad j = g+1,\ldots,h
\end{split}\end{equation}
for $\alpha_{ij}\in W(R)$.  The matrix
$M = (\alpha_{ij})\in\M_h(W(R))$ is called the \emph{structure matrix} for
$\mc P$ with respect to $e_1,\ldots,e_h$; its determinant $\det(M)$ is a
unit in $W(R)$ by \cite[Lemma~9]{zink:display}.   

\subparagraph
In terms of the structure matrix $M$, the map 
$V\inv:{}^V\hat W^{\dsum g}\dsum\hat W^{\dsum(h-g)}\to\hat W^{\dsum h}$
of~\eqref{eq:seqkerlog} is given by
\[ V\inv({}^Vx_1,\ldots,{}^Vx_g,x_{g+1},\ldots,x_h)
= M(x_1,\ldots,x_g,\,{}^Fx_{g+1},\ldots,{}^Fx_h)^t \]
for $\cN\in\Nil_R$ and $x_1,\ldots,x_h\in\hat W(\cN)$; here we are using
the relation
$V\inv(({}^V x_i)\tensor e_i) = x_i\tensor F e_i\in\hat P_\cN$ 
for $1\leq i\leq g$ (cf.~\eqref{eq:Vinv.basechange}).

\begin{notn}
  If $A = (\beta_{ij})$ is a matrix with coefficients in $W(R)$, we write
  $w_n(A)$ for the matrix over $R$ whose entries are $w_n(\beta_{ij})$.
  Note that $w_n(AB) = w_n(A)w_n(B)$ since $w_n$ is a ring homomorphism.
\end{notn}

\paragraph
By \cite[Corollary~11.1.6]{hazewinkel:formal_groups},
the $g$-dimensional formal group $G_K$ over the field $K$ of characteristic
zero is isomorphic (as a formal group) to the formal additive group
$\hat\bG_{a,K}^{\dsum g}$; this 
isomorphism is unique if we require (as we may) that a choice of
parameters for $G_K$ map to the standard parameters for
$\hat\bG_{a,K}^{\dsum g}$.  We call such an 
isomorphism $G_K\isom\hat\bG_{a,K}^{\dsum g}$  a
\emph{logarithm}.  After making a choice of  
parameters, we can view the logarithm as a collection of $g$ power series
$\log_1,\ldots,\log_g$ in $g$ variables with entries in $K$.  The
valuations of the coefficients of the logarithm will be important in the
sequel; we will now calculate these coefficients in terms of the
structure coefficients $\alpha_{ij}$ from~\eqref{eq:alpha.ij}.

\paragraph
Let $\log = (\log_1,\ldots,\log_g):G_K\isom\hat\bG_{a,K}^{\dsum g}$ be the
unique logarithm mapping 
the $i$th parameter provided by Theorem~\ref{thm:paramsfordisplays} and
our choice of basis $e_1,\ldots,e_h$ for $P$ to the $i$th standard
parameter.  We will denote by 
$\bflog$ the composition of $\log$ with the map
$\pi: \hat W^{\dsum h}\to G_K$ from \eqref{eq:seqkerlog}.  We
can think of $\bflog$ as a $g\times h$ matrix of homomorphisms 
$\log_{ij}:\hat W\to\hat\bG_{a,K}$: that is,
\[ \bflog(x_1,\ldots,x_h) = \left(
  \sum_{j=1}^h \log_{1j}(x_j),~\ldots,~\sum_{j=1}^h \log_{gj}(x_j) \right)
\]
for a nilpotent $R$-algebra $\cN$ and elements $x_1,\ldots,x_h\in\hat W(\cN)$.
It is a basic theorem of Cartier
theory (``Cartier's first theorem'',
\cite[Theorem~27.7.5]{hazewinkel:formal_groups}) that  
$\hat W$ represents the ``formal curves functor'' for a formal group.  What
this means concretely in the case of $\hat\bG_{a,K}$ is that if
$f:\hat W\to\hat\bG_{a,K}$ is a homomorphism of formal groups then there
are uniquely determined $a_n\in K$ such that 
$f(x) = \sum_{n=0}^\infty a_n w_n(x)$ for all nilpotent $R$-algebras $\cN$
and all $x\in\hat W(\cN)$
(this is a finite sum because for all $x\in\hat W(\cN)$ we have $w_n(x) =
0$ when $n\gg 0$).  Hence we may write
$\log_{ij}(x) = \sum_{n=0}^\infty a_{n,ij} w_n(x).$
Letting $\ba_n$ be the $g\times h$ matrix $(a_{n,ij})$, we have
\begin{equation}
  \label{eq:bflog}
  \bflog(x_1,\ldots,x_h) 
  = \sum_{n=0}^\infty \ba_n 
  \begin{bmatrix}
    w_n(x_1) \\ \vdots \\ w_n(x_h)
  \end{bmatrix}.
\end{equation}
In terms of the formal group parameters 
$(y_1,\ldots,y_g)\mapsto\pi([y_1],\ldots,[y_g],0,\ldots,0)$
on $G$ (and hence on $G_K$) provided by
Theorem~\ref{thm:paramsfordisplays}, we have  
\[\begin{split}
  \log(y_1,\ldots,y_g) &= \bflog([y_1],\ldots,[y_g],0,\ldots,0) \\
  &= \sum_{j=1}^g(\log_{1j}([y_j]),~\ldots,~\log_{gj}([y_j])) \\
  &= \sum_{n=0}^\infty\sum_{j=1}^g
  (a_{n,1j}y_j^{p^n},~\ldots,~a_{n,gj}y_j^{p^n})
  \qquad\text{ since } w_n([y_j]) = y_j^{p^n} \\
  &= \sum_{n=0}^\infty \ba_n(y_1^{p^n},\ldots,y_g^{p^n},0,\ldots,0)^t.
\end{split}\]
In particular, since $\log(y) = y + O(y^2)$ we have that $a_{0,ij}=\delta_{ij}$ for
$1\leq i,j\leq g$.

\paragraph
By \eqref{eq:seqkerlog}, the image of
\[ V\inv-\Id:~ {}^V\hat W^{\dsum g}\dsum\hat W^{\dsum(h-g)}
\To \hat W^{\dsum h} \]
is contained in the kernel of $\bflog$.  
Hence for every nilpotent $K$-algebra $\cN$ and all 
$x,y\in\hat W(\mc N)$, we have 
\begin{equation} \label{eq:log.relations} \begin{split}
  \bflog(({}^Vx)\,\be_i) &= \bflog(V\inv(({}^Vx)\,\be_i))
  = \bflog(x\,F\be_i) = \bflog(x\, M \be_i) \\
  \bflog(y\,\be_j) &= \bflog(V\inv(y\,\be_j)) = \bflog(({}^Fy)\,V\inv \be_j)
  = \bflog(({}^Fy)\, M\be_j)
\end{split}\end{equation}
for $i = 1,\ldots,g$ and $j = g+1,\ldots,h$, where $\be_i$ is the
the $i$th standard basis vector in $W(R)^{\dsum h}$.
Expanding $\bflog$ using \eqref{eq:bflog}, the left sides
of~\eqref{eq:log.relations} are
\[\begin{split}
  \bflog(({}^V x)\,\be_i) &= 
  \sum_{n=0}^\infty w_n({}^Vx)\,\ba_n \be_i
  = \sum_{n=1}^\infty p\,w_{n-1}(x)\,\ba_n \be_i \\
  \bflog(y\be_j) &= \sum_{n=0}^\infty w_n(y)\,\ba_n \be_j
\end{split}\]
and the right sides of~\eqref{eq:log.relations} are
\[\begin{split}
  \bflog(x\, M \be_i) 
  &= \sum_{n=0}^\infty w_n(x)\,\ba_n w_n(M)\be_i \\
  \bflog(({}^Fy)\, M\be_j)
  &= \sum_{n=0}^\infty w_{n+1}(y)\,\ba_n w_n(M)\be_j.
\end{split}\]
Therefore,
\begin{equation} \label{eq:log.relations.2} \begin{split}
  \sum_{n=1}^\infty p\, w_{n-1}(x)\,\ba_n \be_i
  &= \sum_{n=0}^\infty w_n(x)\, \ba_n w_n(M)\be_i \\
  \sum_{n=0}^\infty w_n(y)\, \ba_n \be_j
  &= \sum_{n=0}^\infty w_{n+1}(y)\, \ba_n w_n(M)\be_j.
\end{split}\end{equation}
As these equalities hold for all $\mc N$, they are in fact
equalities of tuples of power series in the Witt coordinates of $x$ and $y$.
Examining the equation
\[ w_n(x_0,x_1,\ldots) = x_0^{p^n} + px_1^{p^{n-1}} + \cdots
+ p^{n-1}x_{n-1}^p + p^nx_n, \]
we see that the sets of monomials appearing in the polynomials
$w_n(x_0,x_1,\ldots)$ and $w_{n'}(x_0,x_1,\ldots)$ are disjoint when
$n\neq n'$.  Thus we may separate the sums in~\eqref{eq:log.relations.2}
to obtain:
\[\begin{split}
  p\, w_n(x)\,\ba_{n+1}\be_i &= w_n(x)\,\ba_n w_n(M)\be_i \\
  w_{n+1}(y)\,\ba_{n+1}\be_j &=  w_{n+1}(y)\,\ba_n w_n(M)\be_j \\
  \ba_0\be_j &= 0
\end{split}\]
for all $n\geq 0$, $1\leq i\leq g$, and $g+1\leq j\leq h$.  Cancelling the
$w_n(x)$ and $w_{n+1}(y)$, and recalling that $a_{0,ij}=\delta_{ij}$ for
$1\leq i,j\leq g$, we obtain the fundamental recursive equation: 
\begin{equation} \label{eq:fundeq}
\begin{split}
  \ba_0 &= [ I_g ~ 0 ] \\
  \ba_{n+1} &= \ba_n w_n(M)\bp\inv,
\end{split}
\end{equation}
where $I_g$ is the identity matrix and $\bp$ is the $h\times h$
diagonal matrix whose diagonal entries 
are $(p,p,\ldots,p,1,\ldots,1)$ ($g$ entries are $p$).
Note that \eqref{eq:fundeq} uniquely determines $\bflog$ in terms of the
structure matrix $M$ of $\mc P$ with respect to $e_1,\ldots,e_h$.

\begin{eg} Let \label{eg:display.gm.3}
  $\cG_m$ be the display from Examples~\ref{eg:display.gm}
  and~\ref{eg:display.gm.2}.  Its structure matrix is $M=(1)$ with respect
  to the canonical basis $e_1 = 1$ of $P = W(R)$, so by~\eqref{eq:fundeq},
  its logarithm is 
  \[ \log(x) = x + \frac{x^p}p + \frac{x^{p^2}}{p^2} + \frac{x^{p^3}}{p^3}
  + \cdots. \]
  We recognize this as the $p$-typical logarithm for the formal group
  $\hat\bG_m$, in the sense of~\cite[\S16.3]{hazewinkel:formal_groups}.
  This means that $\hexp = \exp\circ\log:\BT_{\cG_m}\to\hat\bG_m$ is
  an isomorphism of formal groups defined over $R$, as in
  Example~\ref{eg:display.gm.2}. 

\end{eg}

\section{The rigid generic fiber of a $p$-divisible formal group}
\label{sec:genericfiber}

\paragraph
In this section we give a geometric interpretation of a $p$-divisible
formal group $G$, and in particular of the kernel of its logarithm.  We
also derive the well-known fact that $\ker(\log) = G[p^\infty]$, in an
appropriate geometric sense.

\paragraph
Let $K$ be a field of characteristic zero that is complete with respect to
a non-Archimedean valuation $\ord: K\to\R\cup\{\infty\}$
satisfying $\ord(p) = 1$.   Let $|x| = p^{-\ord(x)}$ be the
associated absolute value.  Let $R$ be the ring of integers in $K$ and let
$k$ be its residue field.  It is convenient at this point to set our notation
involving rigid geometry; our primary reference for rigid $K$-analytic
spaces is~\cite{bgr:nonarch}. 

\begin{notn}
  Let $|K^\times| = \{|x|~:~x\in K^\times\}$ be the \emph{value group} of
  $K$, and let
  \[ \sqrt{|K^\times|} = |\bar K^\times|
  = \{x\in\R_{>0}~:~ x^n\in|K^\times|\text{ for some } n \}, \]
  where $\bar K$ is an algebraic closure of $K$.
  For $\rho\in\sqrt{|K^\times|}$ we define 
  \[ T_{K,g,\rho} = T_{g,\rho} =\left\{\sum a_\nu X^\nu\in K\ps{X_1,\ldots,X_g}
    ~:~|a_\nu|\rho^{|\nu|}\to 0\text{ as }|\nu|\to\infty\right\} \]
  (we omit $K$ from the notation when the ground field is clear from the context),
  which is equipped with the norm 
  $|\sum a_\nu X^\nu|_\rho = \max\{|a_\nu|\rho^{|\nu|}\}$, where
  $|\nu| = \nu_1+\cdots+\nu_g$.
  The associated affinoid space is the closed $g$-ball of radius $\rho$,
  and is denoted $\B^g_K(\rho) = \Sp(T_{K, g,\rho})$.  Likewise we set
  $\D^g_K(\rho) = \bigcup_{\mu < \rho} \B^g_K(\mu)$, the open ball of
  radius $\rho$.  For brevity we write $T_g = T_{K,g} = T_{K,g,1}$, 
  $\B^g_K = \B^g_K(1)$, and $\D^g_K = \D^g_K(1)$.

  If $\sX$ is a rigid space and $x\in\sX$ is a point then we denote
  the residue field at $x$ by $\kappa(x)$; this is a finite extension of $K$.
\end{notn}

\paragraph Let \label{par:generic.fiber}
$G$ be a $p$-divisible formal group over $R$ of height $h$ and
dimension $g > 0$, and let $\fG \cong \Spf(R\ps\bX)$ (where
$\bX=X_1,\ldots,X_g$) be the associated formal Lie group over $R$ with
$G\cong\fG[p^\infty] = \{\fG[p^n]\}_{n\in\N}$.  In particular,
$G[p^n]\cong\fG[p^n] = \Spec(A_n)$, where $A_n = R\ps\bX/[p^n](\bX)$ 
is a free $R$-module of rank $p^{nh}$.  Let $F$ be the formal group law on
$R\ps\bX$ determining the group structure on $\fG$.
Let $\sG=\D^g_K$, and endow $\sG$ with the structure of 
$K$-analytic group via the convergent power series $F$.  We call $\sG$ the
\emph{rigid generic fiber} of $G$ (with respect to a choice of parameters
for $G$).  The $p^n$-torsion $\sG[p^n]$ is the closed analytic subspace 
of $\sG$ defined by the equations $[p^n](\bX)$.
For any $\rho\in\sqrt{|K^\times|}$, $0<\rho<1$ the natural inclusion
$R\ps\bX\to T_{g,\rho}$ induces a homomorphism 
$K\tensor_R A_n \to T_{g,\rho}/[p^n](\bX)$, whence we obtain a
natural morphism of $K$-analytic groups 
$\sG[p^n]\cap\B^g_K(\rho)\to G[p^n]\tensor_R K$. Passing to the direct
limit, we have a natural morphism $\sG[p^n]\to G[p^n]\tensor_R K$.

\begin{lem} The \label{lem:generic.fiber}
  natural map $\sG[p^n]\to G[p^n]\tensor_R K$ of finite (\'etale)
  $K$-analytic groups is an isomorphism for all $n\geq 1$.
\end{lem}

\pf
  For $\xi\in G[p^n](\bar K)$ let $\xi_i = X_i(\xi)\in\bar K$.  Since
  $A_n$ is local, we have $|\xi_i| < 1$.  Let
  $\rho\in\sqrt{|K^\times|}$, $0<\rho<1$ be such that 
  $|\xi_i|\leq\rho$ for all $\xi\in G[p^n](\bar K)$.  Let
  $x_i$ be the image of $X_i$ under the quotient map
  $R\ps\bX\to A_n$, and define a homomorphism
  $T_{g,\rho}\to K\tensor_R A_n$ by 
  $\sum a_\nu X^\nu\mapsto\sum a_\nu x^\nu$.  This is well-defined because
  $K\tensor_R A_n\cong\prod K_i$ is a product of finite field extensions
  $K_i$ of $K$, and by the above the series $\sum a_\nu x^\nu$ converges
  in each $K_i$.  Thus we obtain a map
  $T_{g,\rho}/[p^n](\bX)\to K\tensor_R A_n$, which is easily seen to be
  inverse to the natural map defining the morphism
  $\sG[p^n]\cap\B^g_K(\rho)\to G[p^n]\tensor_R K$.
\qed

\begin{rem}
  Suppose that the valuation on $K$ is discrete.  Then
  Berthelot has defined a ``rigid generic fiber'' functor
  $\mf X\rightsquigarrow\mf X_\rig$ from the
  category of locally noetherian adic formal schemes $\mf X$ over $\Spf(R)$
  whose reduction is a scheme locally of finite type over $\Spec k$, to the
  category of rigid analytic spaces over $K$; this is explained
  in~\cite[\S7]{dejong:crystalline}.  
  In this case $\sG$ is identified with $\fG_\rig$, and
  Lemma~\ref{lem:generic.fiber} follows from the fact that 
  Berthelot's functor is compatible with fiber products.

  Using Berthelot's functor, it is not necessary to choose parameters for
  $\fG$ in order to define $\sG$.  We prefer to use an ad-hoc construction
  in this case instead of imposing noetherian hypotheses on $R$,
  especially since we will need to choose parameters for $\fG$ 
  in \S\ref{sec:canonsubgp}.

\end{rem}

\paragraph
By Lemma~\ref{lem:generic.fiber}, $\sG[p^n]$ is a finite \'etale
$K$-analytic group of order $p^{nh}$.  The following proposition 
relates the logarithm from \S\ref{sec:logarithm} with
$\sG[p^\infty] \coloneq \bigcup_{n=1}^\infty\sG[p^n]$
when $G$ comes from a display.

\begin{prop} Let \label{prop:kerlog}
  $\cP = (P,Q,F,V\inv)$ be a display over $R$, let
  $P = T\dsum L$ be a normal decomposition, and choose bases for $T$ and
  $L$. Let $\log:G_K\isom\hat\bG_{a,K}^{\dsum g}$ be the logarithm of 
  $G = \BT_\cP$ with respect to the induced choice of parameters for $G$, as in
  \S\ref{sec:logarithm}.   Let $F$ be the formal group
  law for $G$ and let $\sG$ be its rigid generic fiber, as above.

  The $g$-tuple of power series $\log = (\log_1,\ldots,\log_g)$ converges
  on the analytic open unit ball $\D_K^g$ and defines a homomorphism
  $\log: \ms G\to\sG_a^{\dsum g}$ of $K$-analytic groups, where $\sG_a$
  is the additive group on the rigid affine line.  Let
  $\ker(\log)$ be its kernel, and let $\rho\in\sqrt{|K^\times|}$, 
  $0<\rho < 1$.  For $n\gg 0$ (depending on $\rho$) we have an equality 
  \[ \ker(\log)\cap\B^g_K(\rho) = \ms G[p^n]\cap\B^g_K(\rho) \]
  of (finite) $K$-analytic subspaces of $\B^g_K(\rho)$.
\end{prop}

\pf
  Using \eqref{eq:fundeq}, one sees that 
  $\ord a_{n,ij} \geq -n$, from which it follows immediately that the
  radius of convergence of $\log$ is at least $1$, i.e., that 
  $\log_i\in\Gamma(\D^g_K,\O_{\D^g_K})$.  Therefore $\log_i$ defines a map
  $\D^g_K\to\ms G_a$, so we obtain a $K$-analytic morphism
  $\log:\ms G\to\ms G_a^g$.  Since $\log(F(X,Y)) = X+Y$,
  $\log$ is a homomorphism of $K$-analytic groups.

  To prove the second assertion, first we show that
  $\ker(\log)(\bar K) = \bigcup_{n\geq 1}\sG[p^n](\bar K)$.
  As $\ms G^g_a$ has no additive torsion, any point
  $x\in\D_K^g(\bar K)$ such that $[p^n](x) = 0$ must satisfy $\log(x) = 0$.  
  Conversely, suppose that 
  $\log(x) = 0$ for $x = (x_1,\ldots,x_g)\in\D_K^g(\bar K)$.  By the
  non-Archimedean inverse  
  function theorem, the power series inverse $\exp$ to $\log$ has nonzero
  radius of convergence, so $\log$ is injective on some ball
  $\B^g_K(\rho_0)$ around $0$ with $\rho_0 > 0$.   Since 
  $[p^n](\bX) = p^n\bX + O(\bX^2)$ as power series over $R$, we have
  $[p^n] x\to 0$ as $n\to\infty$, so
  for some $n$ we have $[p^n] x\in\B^g_K(\rho_0)$.
  Then $0 = [p^n]\log(x) = \log([p^n]x)$ implies $[p^n]x = 0$.
  Moreover the value of $n$ can be bounded in terms of
  $|x| = \max\{|x_1|,\ldots,|x_g|\}$, so
  $\sG[p^\infty]\cap\B^g_K(\rho) \subset\sG[p^n]$ for
  large enough $n$.

  Since $\log$ has a local inverse near the identity, the $K$-analytic
  group $\ker(\log)$ is \'etale.  Thus for fixed $\rho$ and sufficiently
  large $n$ depending on $\rho$, 
  $\ker(\log)\cap\B^g_K(\rho)$ and $\ms G[p^n]\cap\B^g_K(\rho)$
  are equal, being \'etale subspaces of $\B^g_K(\rho)$
  with the same geometric points.
\qed

\begin{rem}
  Proposition~\ref{prop:kerlog} is not surprising in view
  of the fact that 
  $\log = \lim_{n\to\infty}{[p^n]}/{p^n}$
  in the sense of \cite[p.64]{hazewinkel:formal_groups}. 
\end{rem}

\section{Statement of the main theorem and preliminary reductions}
\label{sec:canonsubgp}  

\paragraph
In this section we define the canonical subgroup of a truncated
$p$-divisible group $G$, we define a measure of
non-ordinarity of $G$, and we state our main theorem
(Theorem~\ref{thm:strong.mainthm}) relating these.  We 
also make some preliminary reductions regarding the proof.

\begin{notn}
For the rest of this paper
$R$ is a complete rank-$1$ valuation ring of mixed characteristic
$(0,p)$ with fraction field $K$, residue field $k$, valuation $\ord$
normalized so that $\ord(p) = 1$, and absolute value 
$|\cdot| = p^{-\ord(\cdot)}$.  We do \emph{not} assume that $R$ is
noetherian or that its residue field is perfect.
\end{notn}

\paragraph
Let $G$ be a connected $p$-divisible group over $R$ of dimension $g$ and
let $\sG\cong\D^g_K$ be its rigid generic fiber in the sense
of~\parref{par:generic.fiber}. 
As indicated in the introduction, we would like to define the level-$n$
canonical subgroup of $G$ to be a $K$-subgroup of 
$G[p^n]\tensor_R K = \sG[p^n]$ with geometric structure
$(\Z/p^n\Z)^g$ whose geometric points are smaller than all
other points of $\sG[p^n](\bar K)$ with respect to the size function 
$|(x_1,\ldots,x_g)| = \max\{|x_i|\}_{i=1}^g$.  This size function is
in fact determined by the $R$-group structure on $G[p^n]$, as follows.

\paragraph
Let $G = \Spec(A)$ be a finite flat connected $R$-group scheme, let $I$ be
the augmentation ideal of $A$, and let $\xi\in G(\bar K)$.  
Since $A$ is local, we have $|f(\xi)| < 1$ for all 
$f\in I$.  For $r_1,\ldots,r_n\in A$ and $f_1,\ldots,f_n\in I$ we have
$|\sum_{i=1}^n r_i f_i(\xi)| \leq \max\{|f_i(\xi)|\}_{i=1}^n$, so if the
residues of $f_1,\ldots,f_n$ generate $I/I^2$ as an $R$-module then 
\[ |\xi| \coloneq \sup_{f\in I} |f(\xi)| = \max\{|f_i(\xi)|\}_{i=1}^n < 1. \]
Suppose that $\xi,\xi'\in G(\bar K)$ and $|\xi|,|\xi'|\leq\rho$ for some
$0 < \rho \leq 1$.  Let $f\in I$, and let 
$1\tensor f+f\tensor 1+ \sum_{i=1}^n f_i\tensor f_i'\in A\tensor_R A$ be
the image of $f$ under the comultiplication map, where $f_i,f_i'\in I$.
Then
\[ |f(\xi+\xi')| 
= \left|f(\xi) + f(\xi') + \sum_{i=1}^n f_i(\xi)f_i'(\xi')\right|
\leq \rho, \]
so $\{\xi\in G(\bar K)~:~|\xi|\leq\rho\}$ is a subgroup of $G(\bar K)$. 
Since $|\xi| = |\sigma(\xi)|$ for all $\sigma\in\Gal(\bar K/K)$, this
subgroup descends to a $K$-subgroup of $G\tensor_R K$, denoted
$G_{\leq\rho}$.  To summarize:

\begin{defn}
  Let $G = \Spec(A)$ be a finite flat connected $R$-group scheme, let $I$
  be the augmentation ideal of $A$, and let $\xi\in G(\bar K)$.  We define
  the \emph{size of $\xi$} to be
  \[ |\xi| = \sup_{f\in I} |f(\xi)| < 1, \]
  and for $0<\rho\leq 1$ we let $G_{\leq\rho}$ be the $K$-subgroup of
  $G\tensor_R K$ whose geometric points are
  \[ G_{\leq\rho}(\bar K) = \{ \xi\in G(\bar K)~:~ |\xi| \leq \rho \}. \]
\end{defn}

\begin{remsub}
  Let $G$ be a finite flat connected $R$-group scheme.
  \begin{enum}
  \item Let $H\subset G$ be a finite flat closed $R$-subgroup scheme. For
    $\xi\in H(\bar K)$ with image $\xi'$ in $G(\bar K)$,
    we have $|\xi| = |\xi'|$.  
  \item Let $K'$ be a complete valued field extension of $K$ with
    valuation ring $R'$ and let $G' = G\tensor_R R'$.  Let 
    $\xi\in G(\bar K)$, and let $\xi'$ be the image of $\xi$ in 
    $G(\bar K') = G'(\bar K')$ under some embedding $\bar K\inject\bar K'$.
    Then $|\xi| = |\xi'|$.  In particular, if $0 < \rho \leq 1$ then 
    $G_{\leq\rho}\tensor_K K' = G'_{\leq\rho}$.
  \item If $\rho \leq \mu$ then $G_{\leq\rho}$ is a $K$-subgroup of
    $G_{\leq\mu}$.     
  \end{enum}
\end{remsub}

\subparagraph \label{par:size.pdivgp}
  Let $G$ be a $g$-dimensional connected $p$-divisible group over $R$,
  let $\fG\cong\Spf(R\ps\bX)$ be the associated formal Lie group with
  $G = \fG[p^\infty]$, and
  let $\sG\cong\D^g_K$ be its rigid generic 
  fiber.  By Lemma~\ref{lem:generic.fiber}, for all $n\geq 1$ we have a
  natural isomorphism $G[p^n]\tensor_R K\isom\sG[p^n]$.  By definition
  $G[p^n] = \Spec(A)$ where $A = R\ps\bX/[p^n](\bX)$, and the augmentation
  ideal $I$ of $A$ is generated by the residues of
  $X_1,\ldots,X_g$.  Hence if $\xi\in G[p^n](\bar K)$ and if 
  $x = (x_1,\ldots,x_g)$ is its image in $\D^g_K$ then 
  $|\xi| = \max\{|x_i|\}_{i=1}^g$.  It follows that
  $G[p^n]_{\leq\rho} = \sG[p^n]\cap\B^g_K(\rho)$.  

\begin{defn} Let \label{defn:canon.subgp}
  $N\geq 1$, and let $G$ be a $\BT_N$ over $R$ of dimension
  $g$.  If there exists $\rho \leq 1$ such that 
  $G^\circ_{\leq\rho}(\bar K)\cong(\Z/p^N\Z)^g$ then we call
  $G^\circ_{\leq\rho}$ the \emph{level-$N$ canonical subgroup of $G$}, and
  we say that this canonical subgroup admits the \emph{radius} $\rho$.

  If $G$ is a $p$-divisible group over $R$ of dimension $g$ then
  the \emph{level-$N$ canonical subgroup of $G$} is defined to be the level-$N$
  canonical subgroup of $G[p^N]$, if it exists.
\end{defn}

\subparagraph
Note that Definition~\ref{defn:canon.subgp} only depends on
$G^\circ$ and is insensitive to valued field extension of $K$.

\begin{remsub}
  Let $G$ be a connected $\BT_N$ over $R$ of dimension $g$ and height $h$.
  \begin{enum}
  \item Note that
    $G_{\leq 1}(\bar K) = G(\bar K)\cong(\Z/p^n\Z)^h$, and 
    $G_{\leq\rho}(\bar K) = \{0\}$ for small enough $\rho > 0$.  In
    particular, if $g = h$ then $G_{\leq 1}$ is the level-$N$ canonical
    subgroup of $G$.

  \item If $G_{\leq\rho}$ is the level-$N$ canonical subgroup of $G$ then 
    $G[p^n]_{\leq\rho}$ is the level-$n$ canonical subgroup of 
    $G[p^n]$ for all $1\leq n\leq N$.
  \end{enum}
\end{remsub}

\paragraph
The main theorem of this paper is a nontrivial condition on a truncated
$p$-divisible group $G$ of level $n$ that is sufficient for the canonical
subgroup of level $n$ to exist.  This involves the Hasse invariant of $G$,
defined as follows.

\paragraph Let \label{par:dV}
$G = \Spec(A)$ be a connected $\BT_n$ over $R$ for some
$n\geq 1$.  Let $R_0 = R/pR$, let $A_0 = A\tensor_R R_0$, let $I_0$ be the
augmentation ideal of $A_0$, let 
$G_0 = G\tensor_R R_0 = \Spec(A_0)$, and let
$F:G_0\to G_0^{(p)}$ and $V: G_0^{(p)}\to G_0$ be the relative Frobenius
and Verschiebung homomorphisms, respectively.
By~\cite[Propositions~II.2.1.2 and~II.3.3.11]{messing:thesis}, 
$I_0/I_0^2$ is a finite free $R_0$-module, so $\Lie(G_0)$
(resp. $\Lie(G_0^{(p)})$) is a finite free
$R_0$-module of the same rank.  After choosing bases for 
$\Lie(G_0)$ and $\Lie(G_0^{(p)})$, we can regard the map
$\dV: \Lie(G_0^{(p)})\to\Lie(G_0)$ as a square matrix with entries in
$R_0$; the determinant $\det(\dV)$ is then defined up to multiplication by
a unit in $R_0$.
For $r\in R_0$ let $\td r\in R$ be any lift of $r$, and note that
$\min\{\ord(\td r), 1\}$ only depends on the ideal $rR_0$.

\begin{defn}
  If $G$ is a truncated $p$-divisible group over $R$, 
  the \emph{Hasse invariant} of $G$ is defined to be
  \[ H(G) = H(G^\circ) = \min\{\ord(\det(\dV)^\sim), ~1\}. \]
  If $G$ is a $p$-divisible group over $R$, the \emph{Hasse invariant} of
  $G$ is defined to be $H(G[p])$.
\end{defn}

\begin{remsub} \label{rem:after.hasse}
\begin{enum}
\item
  If $G$ is a truncated $p$-divisible group over $R$ then 
  $H(G) = H(G[p])$  because $\Lie(G_0) = \Lie(G_0[p])$;
  cf. Remark~\ref{rem:captures.Lie}. 

\item
    Let $G$ be a connected $\BT_N$ over $R$ of dimension $g>0$ and height
    $h$ for some  $N\geq 1$. 
    We have $H(G) = 0$ if and only if the Verschiebung $V$ is an isomorphism on
    the tangent space of $G_0$, which is to say that $V:G_0^{(p)}\to G_0$
    is an isomorphism, or equivalently the Frobenius $F$ is an
    isomorphism on the Cartier dual group $G_0^*$.  Equivalently the dual
    $G^*$ is \'etale, which is to say $\dim(G^*)=0$.  As 
    $\dim(G^*) + \dim(G) = h$, this says exactly that  
    $h = g$, i.e., that $G$ is multiplicative.  In this case we say that $G$
    is \emph{ordinary}; the theory of the canonical subgroup is
    uninteresting when $G$ is ordinary, as $\sG[p^N]_{\leq 1}$ is trivially
    the level-$N$ canonical subgroup $G$ for all $n\geq 1$.  
    In general we consider $H(G)$  as a measure of non-ordinarity of $G$.

\item
    As a converse to~(ii), if $G$ is a $p$-divisible group over $R$
    such that the canonical subgroup of level $N$ exists for all $N\geq 1$, then
    $H(G)=0$; cf. Remark~\ref{rem:all.cs.ordinary}.

\item
    Let $\cP = (P,Q,F,V\inv)$ be a display over $R$, and let 
    $G = \BT_\cP$ be the associated $p$-divisible group.  We can calculate
    the Hasse invariant of $G$ in terms of the display $\cP$ as follows.
    Choose a normal  decomposition $P = T\dsum L$, choose
    $W(R)$-bases for $T$ and $L$, and let $M = (\alpha_{ij})$ be the
    structure matrix for $\cP$ with respect to these bases, as
    in~\eqref{eq:alpha.ij}.  Let $A = (\alpha_{ij})_{i,j=1}^g$ be the
    upper-left $g\times g$ submatrix of $M$, the matrix
    for the composite map $T\to P\overset{F}\To P\to T$.  The tangent
    space $\Lie(G_0)$ of $G_0=G\tensor_R R_0$ is identified with
    $(P/Q)\tensor_R R_0\cong T_0\coloneq T\tensor_R R_0$, and likewise
    $\Lie(G_0^{(p)})\cong T_0\tensor_{R_0,\Frob} R_0$.  Hence
    $\dV:\Lie(G_0^{(p)})\to\Lie(G_0)$ 
    is given by the matrix $w_0(A)\mod p$ with respect to our choice of
    basis, so
    \[ H(G) = \min\{\ord(\det(w_0(A))),~ 1\}. \]
    See~\cite[Example~23]{zink:display}.
\end{enum}
\end{remsub}

\paragraphnn
Our goal is to prove the following

\begin{thm} Let \label{thm:strong.mainthm}
  $G$ be a truncated $p$-divisible group of level $N\geq 1$ over $R$.
  If $H(G) < (p-1)/p^N$ then the level-$N$ canonical subgroup of $G$ exists,
  and it admits the radius $\rho = p^{-r}$ where
  \[ r = \frac 1{p^{N-1}(p-1)} - \frac {H(G)}{p-1}. \]
\end{thm}

\begin{remsub} \label{rem:after.mainthm}
  \begin{enum}
  \item It is important to note that the bound in
    Theorem~\ref{thm:strong.mainthm} is independent of the height and
    dimension of $G$.  

  \item As mentioned in~\parref{par:katz.bound}, when $g=1$ and $h\leq 2$ 
    Katz \emph{et al} have shown that the canonical subgroup of
    level $N$ exists if and only if $H(G) < 1/p^{N-2}(p+1)$ (with the
    same radius).  This bound is larger than the
    bound in Theorem~\ref{thm:strong.mainthm} by a factor of
    $p^2/(p^2-1)$; we do not know if it holds in higher
    dimensions.  See Remark~\ref{rem:why.our.bound}.

  \item Since the Hasse invariant and canonical subgroup of $G$ are
    intrinsic to $G^\circ$, Theorem~\ref{thm:strong.mainthm} is in fact a
    statement about \emph{connected} truncated Barsotti-Tate groups.
    Furthermore, if $K'$ is a complete valued field extension of $K$ with
    valuation ring $R'$ and if $G' = G\tensor_R R'$, then $H(G) = H(G')$ and 
    $G^\circ_{\leq\rho}$ is the level-$N$ canonical subgroup of $G$ if and only
    if $(G')^\circ_{\leq\rho}$ is the level-$N$ canonical subgroup of
    $G'$, so Theorem~\ref{thm:strong.mainthm} is insensitive to 
    valued field extensions.

  \item Theorem~\ref{thm:strong.mainthm} can be extended to give a
    criterion for the existence of the canonical subgroup of an abelian
    variety over $K$ of arbitrary reduction type.  In addition, 
    a relative version of Theorem~\ref{thm:strong.mainthm}
    holds for algebraic families of abelian varieties over $K$.  This is a
    subject of forthcoming work with Brian  Conrad. 
  \end{enum}

\end{remsub}

\paragraph We \label{par:reduction1}
would like to use the theory of logarithms of displays developed in
\S\ref{sec:logarithm} to prove Theorem~\ref{thm:strong.mainthm}.
However, since we are not assuming the residue field of $R$ to be perfect, we
cannot apply~\cite[Theorem~4.4(e)]{illusie:pdivgp.deformation} to identify
$G$ with the $p^N$-torsion of a $p$-divisible group over $R$.  In
\S\ref{sec:not.noetherian} we will 
prove (Proposition~\ref{prop:comes.from.display}) that any connected
level-$N$ truncated $p$-divisible group $G$ is 
isomorphic to $\BT_\cP[p^N]$ for some display $\cP$ over $R$ when $K$ is
algebraically closed; passing to the completion of the algebraic closure
of $K$ and using Remark~\ref{rem:after.mainthm}(iii), we see that it
suffices to prove the slightly weaker

\begin{thm} Let \label{thm:canonicalsubgp}
  $\cP$ be a display over $R$, and let $G = \BT_\cP$.
  If $H(G) < (p-1)/p^N$ then the level-$N$ canonical subgroup of $G$ exists,
  and it admits the radius $\rho = p^{-r}$ where
  \[ r = \frac 1{p^{N-1}(p-1)} - \frac {H(G)}{p-1}. \]
\end{thm}

\paragraph
The strategy for proving Theorem~\ref{thm:canonicalsubgp} is as follows.  
Let $\cP = (P,Q,F,V\inv)$ be a display over $R$, let $P = T\dsum L$ be a
normal decomposition, let $e_1,\ldots,e_g$ and $e_{g+1},\ldots,e_h$ be
$W(R)$-bases for $T$ and $L$ respectively, and let 
$M = (\alpha_{ij})$ be the structure matrix for $\cP$ with respect to this
basis, as in~\eqref{eq:alpha.ij}.  Let $G = \BT_\cP$ be the associated
$p$-divisible formal group over $R$ (of dimension $g$ and height $h$),
and let $\sG$ be its rigid generic fiber.
Recall that the choice of basis $e_1,\ldots,e_g$ of $T$ provides a 
homomorphism of $K$-analytic groups 
$\log: \sG\To\sG_{a,K}^{\dsum g}$,
which is given by $g$ power series $(\log_1,\ldots,\log_g)$ in $g$
variables calculated in terms of $\cP$ using \eqref{eq:fundeq}.

The kernel of $\log$ being the $p$-power torsion of
$\sG$, we would like to understand the valuations of the coordinates of
the \emph{points} of $\ker(\log)$.  The theory of tropical geometry is
ideally set up to solve such a problem: in what is essentially a
higher-dimensional analogue of a Newton polygon, the ultrametric
inequality will allow us to pinpoint the only possible points 
$x = (x_1,\ldots,x_g)\in (\R\cup\{\infty\})^g$ such that 
$x_i = \ord(\xi_i)$ for some simultaneous root
$\xi = (\xi_1,\ldots,\xi_g)\in\D^g_K(\bar K)$ of $\log_1,\ldots,\log_g$.
More precisely, we will be able to say where these points \emph{cannot}
lie, and then make a deformation to count the number of points of
$\ker(\log)$ contained in $\B^g_K(\rho)$ for suitable $\rho$.  This will
allow us to prove that $\sG[p^\infty]\cap\B^g_K(\rho)$ is the level-$N$
canonical subgroup of $G$.

\paragraphnn
We start by making a preliminary reduction.

\begin{lem} Let \label{lem:upper.triangular}
  $T$ be a finite free $W(R)$-module of rank $g$, and let $F:T\to T$
  be an ${}^F$-linear map.  After possibly making a finite extension of
  $K$, there is a $W(R)$-basis $e_1,\ldots,e_g$ for
  $T$ such that, if $A$ is the matrix for $F$ with respect to this basis,
  then $w_0(A)$ is upper triangular modulo $pR$.
\end{lem}

\pf
  Let $T_0 = T\tensor_{W(R),w_0} R_0$ where
  $R_0 = R/pR$.  Since ${}^F(\cdot)$ preserves the ideal 
  $W(pR) + I_R\subset W(R)$, $F$ preserves the submodule 
  $(W(pR) + I_R) T$ and hence $F$ induces an additive map
  $F_0:T_0\to T_0$ satisfying $F(\lambda v) = \lambda^p F( v)$ for
  $\lambda\in R_0$ and $v\in T_0$.  By \parref{par:W(R)local} and
  Nakayama's lemma, it suffices to find a basis for $T_0$ with respect
  to which the matrix for $F_0$ is upper-triangular.  This statement is
  intrinsic to the $p$-linear endomorphism $F_0$ of the finite free 
  $R_0$-module $T_0$, so it suffices to prove the following:
  \begin{deflist}
  \item[(*)~] Let $M_0$ be a finite free $R_0$-module of rank $g$ and let
    $F_0:M_0\to M_0$ be a $p$-linear endomorphism.  After possibly
    extending scalars to $R'/pR'$, where $R'$ is the ring of integers in 
    a finite field extension of $K$, there is a basis
    for $M_0$ with respect to which $F_0$ is upper-triangular.
  \end{deflist}
  By the standard inductive argument, we
  need only exhibit a ``basis eigenvector'' for $F_0$, i.e., an element
  $v\in M_0$ such that $F_0 v = \lambda v$ for some $\lambda\in R_0$, and
  such that $\{v\}$ extends to a basis for $M_0$.   
  
  Let $M$ be a finite free $R$-module such that $M/pM\cong M_0$,
  and denote the map $M\to M_0$ by $v\mapsto\bar v$.
  Let $x_1,\ldots,x_g\in M$ be any basis, so
  $\bar x_1,\ldots,\bar x_g$ is a basis of $M_0$, and
  $F_0(\bar x_i) = \sum_{i=1}^g \bar\mu_{ij} \bar x_j$ for some 
  $\bar\mu_{ij}\in R_0$.  Let $\mu_{ij}\in R$ be a lift of
  $\bar\mu_{ij}$ (which may be chosen such that $\det(\mu_{ij})\neq0$),
  and define a map $F:M\to M$ by  
  \begin{equation}\label{eq:F_R}
    F\left(\sum_{i=1}^g a_i x_i\right)
    = \sum_{i=1}^g a_i^p\sum_{j=1}^g \mu_{ij} x_j. 
  \end{equation}
  Then $F_0(\bar v) = F(v)\mod p$ and
  $F(\lambda v) = \lambda^p F(v)$ for $v\in M$ and $\lambda\in R$. 
  Suppose that $F(w) = \lambda w$ for $w\in M$ nonzero and
  $\lambda\in R$.  Write $w = \sum_{i=1}^g a_i x_i$, and let
  $\pi\in K$ be an element with $\ord(\pi) = -\max\{\ord(a_i)\}_{i=1}^g$.
  Let $v = \pi w$.  Note that 
  \[ \lambda\pi\inv\, v = \lambda w = F(w) = F(\pi\inv v) = \pi^{-p} F(v)
  \implies F(v) = \pi^{p-1}\lambda \,v, \]
  with $\pi^{p-1}\lambda\in R$ because $\pi^{p-1}\lambda v\in M$. 
  Since $M/(Rv)$ is a finite torsion-free $R$-module, it is flat
  (use~\cite[Theorem~7.8(3)]{matsumura:crt}, noting that every
  finitely generated ideal of $R$ is principal) and hence free
  by~\cite[Theorem~7.10]{matsumura:crt}, so $\{v\}$ extends to a
  basis of $M$ and hence $\bar v$ is a basis eigenvector for $F_0$.
  Thus we are reduced to finding $w\in M$ such that
  $F(w) = \lambda w$.  If $F(w) = 0$ for some nonzero $w\in M$ then
  we are done, so we may assume that there is no such $w$, even after
  a finite field extension of $K$.

  Let $M_K = M\tensor_R K$, and let $\bP = \bP(M_K)\cong\bP^g_K$ be the
  associated projective space over $\Spec(K)$.  Since $F(w)\neq 0$ for all
  nonzero $w\in M$ even after finite extension of $K$, the homogeneous
  polynomials \eqref{eq:F_R} define a morphism $F_K:\bP\to\bP$.  By 
  \cite[Example~16.2.2]{fulton:itheory}, any self-map of $\bP$ has a fixed
  point, so after extending scalars we may assume that $F_K$ has a
  rational fixed point.  Of
  course the rational points of $\bP$ correspond to the lines in $M_K$, so
  there exists a nonzero element $w\in M$ such that $F_K(w) = \lambda w$
  for some  $\lambda\in K^\times$.  Choosing $\pi\in R$ such that 
  $\pi^{p-1}\lambda\in R$, we
  have $F(\pi w) = \pi^{p-1}\lambda\cdot\pi w$, as required.
\qed

\paragraph Recall \label{par:assume.uppertriang}
that we have chosen a normal decomposition $P = T\dsum L$, and that
$T$ is a rank-$g$ free module over $W(R)$.  Let 
$A$ be the matrix for the composite
$T\inject P\overset{F}\To P\to T$, so 
$H(G) = \min\{\ord(\det(w_0(A))),1\}$ by Remark~\ref{rem:after.hasse}(iv).
After making a finite extension of $K$
and choosing a basis for $T$ as in Lemma~\ref{lem:upper.triangular}, we will 
\emph{assume from now on that $w_0(A)$ is upper-triangular modulo $pR$}.
In this case,  
\[ H(G) = \ord(\det(w_0(A))) 
= \sum_{i=1}^g \ord(w_0(\alpha_{ii})) \]
assuming $H(G) < 1$.

\section{A survey of some concepts from tropical geometry}
\label{sec:tropical}

\paragraph
In this section we review the dual concepts of the tropical hypersurface
and the Newton polytope of a power series.  This theory can be seen as a
direct generalization of the theory of Newton polygons; see
Example~\ref{eg:newtonpoly}.  For a more complete theory
see~\cite{jdr:trop_ps} and the references contained therein.

\begin{notn}
  Let $S$ be a set and $T\subset S$ a subset, and let $f:S\to\R$ be any
  function.  Define
  \[ \minset(f,T) \coloneq \big\{ t\in T~:~ f(t) = \inf_{t'\in T} f(t')
  \big\}. \]
  Note that this set could be empty.
\end{notn}

\paragraph 
First we introduce some notions from convex geometry. 
We will make a great deal of definitions, stating facts without proof.
A good reference for this material is~\cite{barvinok:convexity}.
Let $\angles{\cdot,\cdot}$ denote the standard inner product on $\R^g$.
The convex hull of a set of points $S\subset\R^g$ will be denoted
$\conv(S)$.
A \emph{polyhedron} $P$ is a nonempty intersection of finitely many
half-spaces in $\R^g$; i.e., a nonempty subset of the form
\[ P = \bigcap_{i=1}^r \{ x\in\R^g~:~\angles{x,u_i} \geq a_i \} \]
for some elements $u_1,\ldots,u_r\in\R^g$ and $a_1,\ldots,a_r\in\R$.  A
\emph{polytope} is a bounded polyhedron.
For $w\in\R^g$ we define
\[ \face_w(P) = \minset(\angles{\cdot,w},P); \]
this is the locus in $P$ where a linear form attains its infimum on $P$.
A \emph{face} of a polyhedron $P$ is a nonempty
subset of $P$ of the form $\face_w(P)$; this is again a polyhedron.  A
\emph{vertex} of $P$ is a one-point face; we let $\vertices(P)$ be the set
of vertices of $P$.  A polyhedron has finitely
many faces (hence finitely many vertices), and a polytope is the convex
hull of its vertices. 
The \emph{affine span} of a polyhedron $P$ is the smallest affine subspace
containing $P$.  The \emph{dimension} of $P$ is the dimension of
its affine span, and the \emph{relative interior} $\relint(P)$ of $P$ is
the interior of $P$ as a subset of its affine span.  
The relative interior of a face of $P$ is the set of
points not contained in a strictly smaller face.

A \emph{polyhedral complex} is a nonempty
finite collection $\Pi$ of polyhedra satisfying
\begin{enum}\label{page:polycom}
\item If $P\in\Pi$ then every face of $P$ is in $\Pi$.
\item If $P,Q\in\Pi$ and $P\cap Q\neq\emptyset$ then $P\cap Q$ is a face
  of $P$ and a face of $Q$. 
\end{enum}
An element of $\Pi$ is called a \emph{cell}; the
\emph{support} $|\Pi|$ of $\Pi$ is the union of its cells.  A
\emph{polytopal complex} is a polyhedral complex whose cells are
polytopes (i.e., bounded).

\paragraph The \label{par:choice.of.f}
following example illustrates the idea behind using
tropical geometry to determine the valuations of the zeros of a
power series.  To motivate our choice of power series,
let $\cP$ be the display of dimension $2$ and height $2$ over the
valuation ring $R$ of $K$ given by the structure matrix 
$M = \smallmat 1101$, and let $\log = (\log_1,\log_2)$ be the logarithm
of $\BT_\cP$ as in \S\ref{sec:logarithm}.  Using~\eqref{eq:fundeq}
we see that 
\[ \log_1(X,Y)\equiv X + p\inv X^p + p\inv Y^p + p^{-2} X^{p^2} + 2p^{-2}
Y^{p^2} \pmod{X^{p^3},Y^{p^3}}. \]

\begin{eg} Let \label{eg:pretrop}
  $f = X + p\inv X^p + p\inv Y^p\in K[X,Y]$, and choose
  $(\xi,\eta)\in(\bar K^\times)^2$.  If $f(\xi,\eta) = 0$ then at least
  two of the three quantities $|\xi|$, $|p\inv\xi^p|$, and $|p\inv\eta^p|$
  are equal and are at least as large as the third.  In terms of valuations, if
  $(u,v) = (\ord(\xi),\ord(\eta))$, this says that at least two of the
  quantities $u,pu-1,pv-1$ are equal to $\min\{u, pu - 1, pv - 1\}$.
  The locus where $pu-1 = pv-1$ is the line $L_1 = \{u = v\}$, the locus
  where $u = pu-1$ is the line $L_2 = \{u = 1/(p-1)\}$, and the locus
  where $u = pv-1$ is the line $L_3 = \{v = (u-1)/p\}$; we have
  $u = pu-1 = pv-1$ at the point $(u,v) = (1/(p-1),1/(p-1))$.  This is
  illustrated in the right side of Figure~\ref{fig:tropeg1}; we now
  explain the significance of the half-lines in that picture.

  When $(u,v)\in L_1$ and $u = v > 1/(p-1)$ then $u = v > pv-1$, so
  $(\xi,\eta)$ cannot be a zero of $f$; hence
  if $(\xi,\eta)$ is a zero of $f$ and $(u,v)\in L_1$ then $(u,v)$ is
  contained in the ray $R_1 = (1/(p-1),1/(p-1)) + \R_{\geq 0} (-1,-1)$.
  Similar reasoning 
  shows that if $(\xi,\eta)$ is a zero of $f$ then $(u,v)$ is contained in
  one of the rays $R_1$, $R_2 = (1/(p-1),1/(p-1)) + \R_{\geq 0}(0,1)$, or 
  $R_3 = (1/(p-1),1/(p-1)) + \R_{\geq 0} (p,1)$.  In other words, there
  are ``piecewise linear'' necessary conditions on $(u,v)$ for
  $(\xi,\eta)$ to be a zero of $f$; we will use tropical geometry as a
  language to express these conditions.  (In fact, we will see in a moment
  that the ``tropicalization'' of $f$, to be defined shortly, is
  essentially equal to $R_1\cup R_2\cup R_3$.)  
\end{eg}

\paragraph \label{par:trop.f}  
Choose $\rho > 0$, let $r = -\log_p(\rho)\in\R$, and let
$\T^g_K(\rho) = (\B^1_K(\rho)\setminus\{0\})^g$.  
For nonzero $f = \sum_\nu a_\nu X^\nu \in T_{g,\rho}$ let 
\[ H(f) = \{(\nu, \ord(a_\nu))~:~a_\nu\neq 0\}\subset\Z_{\geq 0}^g\times\R; \]
this is called the \emph{height graph} of $f$.
For $w\in\R_{\geq r}^g$ let 
$\Inn_w(f) = \minset(\angles{(w,1),\cdot}, H(f))$, where we are denoting
the inner product on $\R^{g+1}$ by $\angles{\cdot,\cdot}$ as well;
by~\cite[\S8]{jdr:trop_ps}, $\Inn_w(f)$ is a nonempty finite set for all
$w\in\R^g_{\geq r}$.  In Example~\ref{eg:pretrop2} we work this out for
$f$ as in Example~\ref{eg:pretrop}.
For $w\in\R_{\geq r}^g$, define the \emph{initial form}
$\inn_w(f)$ of $f$ to be 
\[ \inn_w(f) = \sum_{(\nu,\ord(a_\nu))\in\Inn_w(f)} a_\nu X^\nu, \]
so $H(\inn_w(f)) = \Inn_w(f)$.
In other words, $\inn_w(f)$ is the (nonzero) sum of those monomials
$a_\nu X^\nu$ such that
\[ \ord(a_\nu) + \angles{w,\nu} 
= \min_{\mu\in\Z_{\geq 0}^g}\{\ord(a_\mu)+\angles{w,\mu}\}. \]
In particular, if $w = (\ord(\xi_1),\ldots,\ord(\xi_g))$ for 
$\xi = (\xi_1,\ldots,\xi_g)\in\T^g_K(\rho)(\bar K)$ then $\inn_w(f)$ is
the sum of those monomials $a_\nu X^\nu$ such that $a_\nu \xi^\nu$ has
minimal valuation among all monomials of $f(\xi)$.  
Define the \emph{tropicalization} of $f$ to be
\[ \Trop(f) = \{w\in\R^g_{\geq r} ~:~ \inn_w(f)\text{ is not a monomial}\}. \]

\begin{eg} Continuing \label{eg:pretrop2}
  with Example~\ref{eg:pretrop}, let $\rho = 1$ (so $r=0$) and $g = 2$, so
  $f\in T_{g,\rho}$ and 
  \[ H(f) = \{ (1,0,0),~ (p,0,-1),~ (0,p,-1) \}. \]
  For $w$ in the part of the relative interior of $R_1$ that lies in
  $\R^2_{\geq 0}$ we have  
  $\Inn_w(f) = \{ (p,0,-1), (0,p,-1) \}$ and
  $\inn_w(f) = p\inv X^p + p\inv Y^p$; equivalently, if
  $w = (u,v) = (\ord(\xi),\ord(\eta))$ then 
  $\ord(p\inv\xi^p)=\ord(p\inv\eta^p)<\ord(\xi)$.  Likewise, for
  $w\in\relint(R_2)\subset\R^2_{\geq 0}$ we have 
  $\Inn_w(f) = \{ (1,0,0), (p,0,-1) \}$ and $\inn_w(f) = X + p\inv X^p$,
  for $w\in\relint(R_3)\subset\R^2_{\geq 0}$ we have 
  $\Inn_w(f) = \{ (1,0,0), (0,p,-1) \}$ and
  $\inn_w(f) = X + p\inv Y^p$, and if $w = (1/(p-1),1/(p-1))$
  then $\Inn_w(f) = H(f)$ and $\inn_w(f) = f$.  As we saw in
  Example~\ref{eg:pretrop}, if $w = (\ord(\xi),\ord(\eta))$ for
  $(\xi,\eta)\in\T^2_K(0)(\bar K)$ then 
  $\inn_w(f)$ is the sum of those monomials of $f$ for which
  some ``internal cancellation'' has to occur in order for $(\xi,\eta)$ to
  be a zero of $f$.

  On the other hand, 
  let $U_{ij}$ denote the interior of the polytope bounded by the rays
  $R_i$ and $R_j$ for $i,j=1,2,3,~i\neq j$.  If 
  $w\in U_{12}\cap\R^2_{\geq 0}$ then $\inn_w(f) = X$, if $w\in U_{23}$ then 
  $\inn_w(f) = p\inv X^p$, and if $w\in U_{13}\cap\R^2_{\geq 0}$ then
  $\inn_w(f) = p\inv Y^p$.  This proves that
  $\Trop(f) = (R_1\cup R_2\cup R_3)\cap\R^2_{\geq0}$.  As we worked out in
  Example~\ref{eg:pretrop}, if $(\xi,\eta)\in\T^2_K(1)(\bar K)$ is a zero
  of $f$ then $(\ord(\xi),\ord(\eta))\in\Trop(f)$.  This is true in
  general, as we will see in~\parref{par:zeros.in.trop}.
\end{eg}

\paragraph \label{par:zeros.in.trop}
We use the notation in~\parref{par:trop.f}.
Let $\ord:\T^g_K(\rho)(\bar K)\to\R^g_{\geq r}$ be the map defined by
\[ \ord(\xi_1,\ldots,\xi_g) = (\ord(\xi_1),\ldots,\ord(\xi_g)). \]
If $\xi = (\xi_1,\ldots,\xi_g)\in\T^g_K(\rho)(\bar K)$ is a zero of 
$f\in T_{g,\rho}$ then we claim that $\ord(\xi)\in\Trop(f)$.  Letting 
$w = \ord(\xi)$, we have 
$\ord(a_\nu\xi^\nu) = \ord(a_\nu) + \angles{w,\nu}$, so if 
$\inn_w(f)$ is a monomial $a_\nu X^\nu$ then $a_\nu\xi^\nu$ has strictly
smaller valuation than $a_\mu\xi^\mu$ for all $\mu\neq\nu$, and thus
$|f(\xi)| = |a_\nu\xi^\nu|\neq 0$ by the ultrametric inequality.  Hence
$\ord$ restricts to a map 
$\ord: V(f)(\bar K)\to\Trop(f)$, where $V(f)$ is the closed subspace of
$\T^g_K(\rho)$ cut out by $f$.  It is a fundamental fact that the closure of
$\ord(V(f)(\bar K))$ is exactly $\Trop(f)$, and if $w\in\Trop(f)$ has
coordinates in $\ord(\bar K^\times)$ then there is a zero $\xi$ of $f$
such that $\ord(\xi) = w$; see~\cite[\S8]{jdr:trop_ps} for a proof in this
context.  This is why we are interested in $\Trop(f)$. 

\paragraph
For $f\in T_{g,\rho}$ nonzero,
the set $\Trop(f)$ has the following natural polyhedral complex structure.
For $w\in\R^g_{\geq r}$ define 
\[\begin{split}
  P_w &= \{ w'\in\R^g_{\geq r} ~:~ \Inn_{w'}(f)\supset\Inn_{w}(f) \}  \\
  &= \bigcap_{\substack{(\nu,\ord(a_\nu))\in\Inn_w(f)\\\mu\in\Z_{\geq 0}^n}}
\{ w'\in\R^g_{\geq r}~:~ \angles{w',\nu-\mu} \leq \ord(a_\mu)-\ord(a_\nu) \}; 
\end{split}\]
note that $w\in P_w$ and if $w'\in P_w$ then $P_{w'}\subset P_w$.
It is proved in~\cite[\S8]{jdr:trop_ps} that 
$\bigcup_{v\in\R^g_{\geq r}} \Inn_v(f)$ is a finite set.  It
follows that $P_w$ is a polyhedron, and that the collection 
$\cC = \{P_w~:~\inn_w(f)\text{ is not a monomial}\}$ of polyhedra in
$\R^g_{\geq r}$ is finite.
In fact $\cC$ is a polyhedral complex (at least in $\R^g_{>r}$), and its
support is clearly 
\[ |\cC| = \{w\in\R^g_{\geq r}~:~ \inn_w(f)\text{ is not a monomial} \}
= \Trop(f). \]
From now on we will use $\Trop(f)$ to denote both $\cC$ and its support.
See Example~\ref{eg:trop1} for a description of the polyhedra $P_w$ for
$f$ as in Example~\ref{eg:pretrop}.

\paragraph
Recall that the height graph $H(f)$ for nonzero $f\in T_{g,\rho}$ lies in
$\Z^g_{\geq 0}\times\R$, with the initial $g$ coordinates encoding the
monomials appearing in $f$, and that the $\Inn_w(f)$ are subsets of
$H(f)$. Let $\pi:\R^{g+1}\to\R^g$ be the projection onto the
first $g$ coordinates.  The \emph{Newton complex} $\New(f)$ of
$f$ is the collection of polytopes in $\R^g$ consisting of the convex hulls
\[ C_w = \conv(\pi(\Inn_w(f))) \]
of the finite sets $\pi(\Inn_w(f))$ for $w\in\R^g_{\geq r}$;
it is shown in~\cite[\S8]{jdr:trop_ps} that the Newton complex is
``almost'' a polytopal complex.  
What we will use is the fact that the complex $\New(f)$ controls the
combinatorics of $\Trop(f)$, in the following sense.

\subparagraph
By definition $P_w = P_{w'}$ if and only if $\Inn_{w'}(f)\subset\Inn_w(f)$
and $\Inn_w(f)\subset\Inn_{w'}(f)$; i.e., if and only if
$\Inn_w(f) = \Inn_{w'}(f)$.  We leave the equality
\begin{equation} \label{eq:Pw.from.Cw}
\Inn_w(f) = H(f)\cap\conv\{x\in H(f)~:~\pi(x)\in\vertices(C_w)\} 
\end{equation}
as an exercise to the reader; in particular,
$C_w=C_{w'}$ if and only if $\Inn_w(f)=\Inn_{w'}(f)$.  Therefore the
association $P_w \leftrightarrow C_w$ for $w\in\Trop(f)$ is a bijective
correspondence between the cells of $\Trop(f)$ and the
positive-dimensional cells of $\New(f)$ (of course $C_w$ is a point if and
only if $\inn_w(f)$ is a monomial).  We call $C_w$ the 
\emph{dual cell} to $P_w$, and vice-versa; this terminology is 
justified by the following paragraphs.

\subparagraph We \label{par:explain.trop}
can  interpret~\eqref{eq:Pw.from.Cw} as follows.  Let $C = C_w$ be
a cell of $\New(f)$, and let $\nu_1,\ldots,\nu_n$ be the vertices of $C$.  
Then $\Inn_w(f) = H(f)\cap\conv\{(\nu_i,\ord(a_{\nu_i}))~:~i=1,\ldots,n\}$, and
hence the dual cell $P = P_w$ is exactly
\[\begin{split}
  P &= \big\{
  w'\in\Trop(f)~:~\pi(\Inn_{w'}(f))\supset\{\nu_1,\ldots,\nu_n\} \big\} \\ 
  &= \big\{ w'\in\Trop(f) ~:~ \angles{w',\nu_i} + \ord(a_{\nu_i})  
  = \min_{\nu\in\Z^g_{\geq 0}} \{\angles{w',\nu} + \ord(a_\nu) \} 
  \text{ for all } i=1,\ldots,n\big\}. 
\end{split}\]
This allows us to recover $P$ from $C$ and $f$ without reference to $w$.
\emph{Note in particular that} 
\begin{equation} \label{eq:Pw.from.Cw.2}
  \angles{w',\nu_1} + \ord(a_{\nu_1}) = \cdots = 
  \angles{w',\nu_n} + \ord(a_{\nu_n}) 
\end{equation}
for all $w'$ in the cell $P$ of $\Trop(f)$ dual to the cell $C$ of
$\New(f)$ with vertices $\{\nu_1,\ldots,\nu_n\}$.

\subparagraph For \label{par:orthogonal}
$w\in\Trop(f)$ the cells $C_w$ and $P_w$ are orthogonal to each other
(in the sense that the vector subspaces of $\R^g$ associated to their
affine spans are orthogonal); this can be seen as follows.
Let $(\nu_1,\ord(a_{\nu_1}))$ and $(\nu_2,\ord(a_{\nu_2}))$ be in
$\Inn_w(f)$ and let $w_1,w_2\in P_w$.  Then
$(\nu_i,\ord(a_{\nu_i}))\in\Inn_{w_j}(f)$ for $i,j=1,2$, so since
\[ \Inn_{w_j}(f) =
\minset\big((\nu,\ord(a_\nu))\mapsto\ord(a_\nu)+\angles{w_j,\nu},~H(f)\big), \]
it follows that
\[\begin{split}
  \ord(a_{\nu_1}) + \angles{w_1,\nu_1} &= \ord(a_{\nu_2})+\angles{w_1,\nu_2} \\
  \ord(a_{\nu_1}) + \angles{w_2,\nu_1} &= \ord(a_{\nu_2})+\angles{w_2,\nu_2},
\end{split}\]
and hence $\angles{w_1-w_2,\nu_1-\nu_2} = 0$.  We leave as an exercise to
the reader to derive that $\angles{w_1-w_2,v_1-v_2}=0$ for all 
$v_1,v_2\in C_w$ using the fact that any
element $v\in C_w$ is of the form $v = \sum_{i=1}^r \alpha_i \mu_i$ where
$\alpha_i\geq 0$, $\sum_{i=1}^r\alpha_i = 1$, and
$(\mu_i,\ord(a_{\mu_i}))\in\Inn_w(f)$ for $i=1,\ldots,r$.

\subparagraph
The ``duality'' between $\Trop(f)$ and $\New(f)$ satisfies other nice
properties.  For example, $\dim(C_w) + \dim(P_w) = g$ (at least when $P_w$
is not contained in the boundary of $\R^g_{\geq r}$ in $\R^g$), and $C_w$
is a face of $C_{w'}$ if and only if $P_{w'}$ is a face of $P_w$ for
$w,w'\in\Trop(f)$.  However, $\Trop(f)$ and $\New(f)$ are not 
dual in any intrinsic way --- it is more accurate to say that they are
both expressions of the combinatorial properties of $f$ (really of
$H(f)$) that live in dual vector spaces.

\begin{eg} Continuing \label{eg:trop1}
  Example~\ref{eg:pretrop2}, we have
  $\Inn_w(f) = \{(p,0,-1),(0,p,-1)\}$ for
  $w\in\relint(R_1)\cap\R^2_{\geq0}$, and hence 
  $P_w = \{w'\in\Trop(f)~:~\{(p,0,-1),(0,p,-1)\}\subset\Inn_{w'}(f)\} =
  R_1\cap\R^2_{\geq 0}$.
  Likewise $P_w = R_2$ for $w\in\relint(R_2)$, 
  $P_w = R_3$ for $w\in\relint(R_3)$, and 
  $P_{(1/(p-1),1/(p-1))} = \{(1/(p-1),1/(p-1))\}$.
  
  For $w\in\relint(R_1)\cap\R^2_{\geq 0}$ we have
  \[ C_w = \conv(\pi(\Inn_w(f))) = \conv(\pi(\{(p,0,-1),(0,p,-1)\}))
  = \conv\{(p,0), (0,p)\}; \]
  this is the line segment $R_1'\coloneq\bar{(p,0)(0,p)}$.  Likewise
  $C_w = R_2'\coloneq\bar{(1,0)(p,0)}$ for $w\in\relint(R_2)$,
  $C_w = R_3'\coloneq\bar{(1,0)(0,p)}$ for $w\in\relint(R_3)$, and
  $C_{(1/(p-1),1/(p-1))}$ is the triangle $\tau$ with vertices
  $(1,0),(0,p),(0,p)$.  For $i=1,2,3$ the cells $R_i$ and $R_i'$ are dual;
  note that they are orthogonal and have complementary dimension.  The
  triangle $\tau$ is dual to the vertex
  $\{(1/(p-1),1/(p-1))\}\in\Trop(f)$.  The vertices of $\New(f)$
  correspond to the connected components of 
  $\R^2_{\geq 0}\setminus\Trop(f)$ (these are the $U_{ij}$ of
  Example~\ref{eg:pretrop2}); they do not have dual cells in
  $\Trop(f)$.  See
  Figure~\ref{fig:tropeg1}. 

  \genericfig[ht]{tropeg1}{The Newton complex and
    tropicalization of the polynomial $f$ of Examples~\ref{eg:pretrop},
    \ref{eg:pretrop2}, and~\ref{eg:trop1}.  The
    cell $R_i$ of $\Trop(f)$ is dual to the line segment $R_i'\in\New(f)$
    for $i=1,2,3$.} 

\end{eg}

\begin{eg} Extending \label{eg:trop2}
  the previous example, suppose that $p\neq 2$, and let
  \[ f = X + p\inv X^p + p\inv Y^p + p^{-2} X^{p^2} + 2p^{-2} Y^{p^2}, \]
  (cf.~\parref{par:choice.of.f}), so
  \[ H(f) = \{(1,0,0),~ (p,0,-1),~ (0,p,-1),~ (p^2,0,-2),~ (0,p^2,-2) \}.
  \]
  We claim that the line segment $\bar{(p,0)(p^2,0)}$ is a cell of
  $\New(f)$, i.e., that there exists $w\in\R^2_{\geq 0}$ such that
  $\Inn_w(f) = \{(p,0,-1),(p^2,0,-2)\}$.  Any such $w = (u,v)$ must satisfy
  $pu - 1 = p^2u - 2$ (cf.~\eqref{eq:Pw.from.Cw.2}), i.e., $u = 1/p(p-1)$,
  in addition to $pu - 1 < \min\{u,pv-1,p^2v-2\}$; this is true for all
  $v\gg 0$, which proves the claim.  We can
  calculate all of the cells of $\New(f)$ in the same way: first we choose
  a subset $S\subset H(f)$, then we solve a system of linear equations and
  inequalities to determine if there exists $w\in\R^g_{\geq 0}$ such that
  $S = \Inn_w(f)$; if so, then $\conv(\pi(S))$ is a cell of $\New(f)$.
  The complex $\New(f)$ is drawn in Figure~\ref{fig:tropeg2}.  

  Once we have drawn $\New(f)$, it is much easier to draw $\Trop(f)$.  
  For instance, the cell $P$ dual to $\bar{(0,p^2)(p^2,0)}$ is contained
  in the line 
  \[ \{(u,v)~:~p^2u-2=p^2v-2\} = \{(u,v)~:~u=v\} \quad
  \text{ (again cf.~\eqref{eq:Pw.from.Cw.2})}. \]
  Since $\bar{(0,p^2)(p^2,0)}$ is a face of the trapezoid with vertices
  $\{(p,0),(p^2,0),(0,p),(0,p^2)\}$, we know that $P$ has the unique
  vertex $(u,v)$ satisfying $pu-1 = pv-1 = p^2u-2 = p^2v-2$.  Since
  $p^2u-2\ll 0$ when $u\ll 0$, we see that
  $P$ is the ray $(1/p(p-1),1/p(p-1)) + \R_{\geq 0}(-1,-1)$ intersected
  with $\R^2_{\geq 0}$.  One can calculate all of the cells of $\Trop(f)$
  in this way: any $P\in\Trop(f)$ is dual to a unique positive-dimensional
  cell $P'\in\New(f)$, and $P$ can be recovered from $P'$ as
  in~\parref{par:explain.trop}.  The tropicalization of $f$ is also
  illustrated in Figure~\ref{fig:tropeg2}.

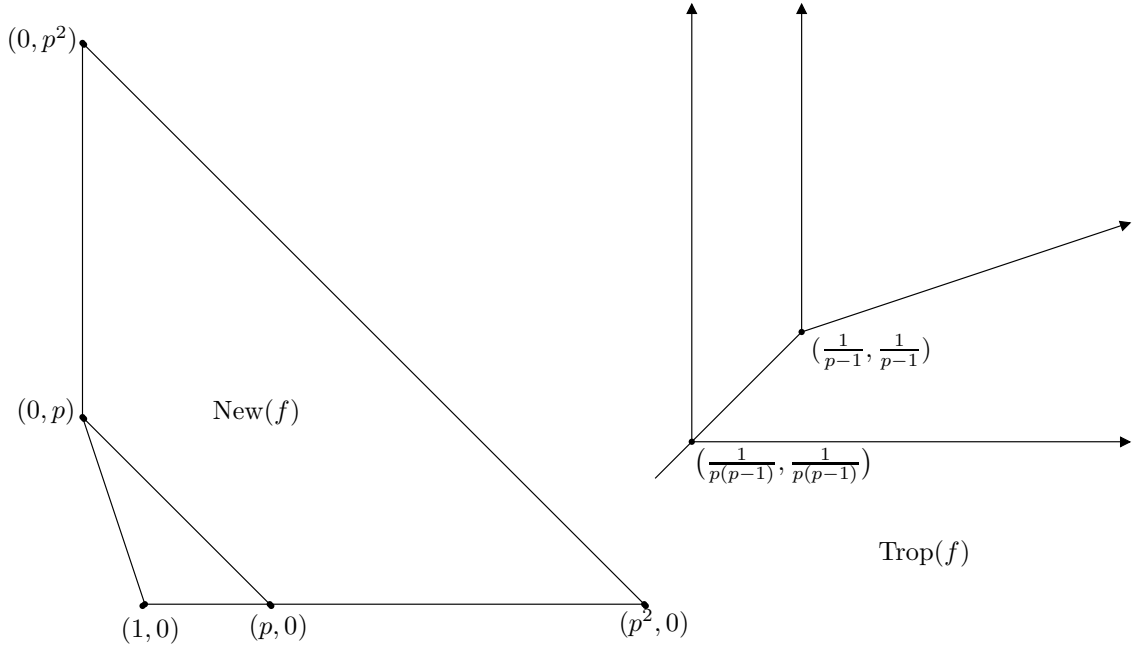
  \genericfig[ht]{tropeg2}{The Newton complex and
    tropicalization of the polynomial $f$ of Example~\ref{eg:trop2}.} 

\end{eg}

\begin{eg} When \label{eg:newtonpoly}
  $g = 1$, tropical geometry basically reduces to the
  classical theory of the Newton polygon, as follows.  
  For a nonzero $f = \sum_{\nu\geq 0} a_\nu X^\nu\in T_{1,\rho}$, the
  Newton polygon of $f$ is by definition the set
  $\cN(f) = \bigcup_{w\in\R_{\geq r}} \conv(\Inn_w(f))$.  
  Hence the support of Newton complex is the projection of
  $\cN(f)$ onto $\R$; the polytopal complex structure simply remembers the
  horizontal lengths of the line segments $\conv(\Inn_w(f))$.  As
  illustrated in Figure~\ref{fig:newtonpoly}, an element $w\in\R_{\geq r}$ is in
  $\Trop(f)$ if and only if $-w$ is the slope of a line segment 
  (of positive length) in $\cN(f)$; this recovers the fact that if $\xi$
  is a zero of $f$ then $-\ord(\xi)$ is the slope of a line segment of
  $\cN(f)$.  In fact, the tropicalization of $f$ along with its Newton
  complex is enough to recover the valuations and multiplicities
  of all of the zeros of $f$ using the  theory of the Newton
  polygon. 
  
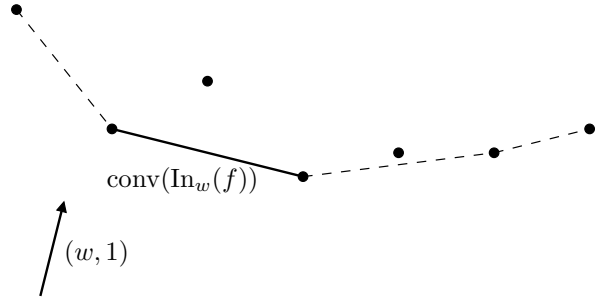
  \genericfig[ht]{newtonpoly}{The Newton polygon of a polynomial $f$.  The
    vertices in the figure are the points of the height graph $H(f)$, and
    the union of the (dotted and solid) lines is the Newton
    polygon of $f$.  For the choice of $w$ indicated in the figure,
    $\Inn_w(f)$ contains two 
    points, so $w\in\Trop(f)$; the slope of the line segment
    $\conv(\Inn_w(f))$ joining these two points is $-w$.} 

  There is an analogue of this fact in higher dimensions.  Namely, if 
  $f_1,\ldots,f_g\in T_{g,\rho}$ are nonzero and
  $w\in\bigcap_{i=1}^g\Trop(f_i)$ is an isolated point contained in
  $\R^g_{>r}$ then there is a formula based on the ``sizes'' of the cells
  of the Newton complexes of the $f_i$ that  
  calculates the (finite, nonzero) number of points 
  $\xi\in\T^g_K(\rho)(\bar K)$ (counted with multiplicity) such that
  $\ord(\xi) = w$ and $f_i(\xi) = 0$ for all $i$.
  This follows from an analogous result for Laurent polynomials proved
  by Eric Katz, Sam Payne, and Brian Osserman; see~\cite[\S11]{jdr:trop_ps} for
  a statement and proof.
\end{eg}

\paragraph
The tropicalization and Newton complex of a nonzero power series
$f$ only depend on its height graph, and are therefore insensitive to
extension of the ground field.  More precisely,
let $K'$ be a complete valued extension of $K$, let $f\in T_{K,g,\rho}$, and let
$f'$ be the image of $f$ in $T_{K',g,\rho}$.
Then $\Trop(f) = \Trop(f')$ and $\New(f) = \New(f')$.  

\paragraph
Let $f$ be a nonzero power series converging on the open unit $g$-ball $\D^g_K$,
and for $0 < \rho < 1$ let $f_\rho$ be the image of $f$ in $T_{g,\rho}$.  
Let $0 < \rho' < \rho$, let $r = -\log_p(\rho) > 0$, and 
let $r' = -\log_p(\rho') > r$.  It is clear that 
$\Trop(f_\rho)\cap\R^g_{\geq r'} = \Trop(f_{\rho'})$, so if we define
\[ \Trop(f) = \bigcup_{\rho > 0} \Trop(f_\rho) \]
(as a subset of $\R^g_{>0}$) then 
$\Trop(f)\cap\R^g_{\geq r} = \Trop(f_\rho)$ for all $\rho$.  This set is
the union $\bigcup_{w\in\R^g_{>0}} P_w$, where
$P_w = \{w'\in\R^g_{>0}~:~\inn_w(f)\subset\inn_{w'}(f)\}$,
the intersection of a polyhedron with $\R^g_{>0}$.
Hence $\Trop(f)$ is a ``polyhedral complex'' with infinitely many cells.
We also let  
$\New(f) = \bigcup_{\rho > 0} \New(f_\rho)$; again this would ``almost'' be a
polyhedral complex if it had finitely many cells.  In any case we have a
bijective correspondence $P_w \leftrightarrow C_w$, so we will still say that
$P_w$ and $C_w$ are \emph{dual cells}.  We
will call $\Trop(f)$ the tropicalization of $f$ and $\New(f)$ the infinite
Newton complex of $f$. 
In what follows all power series $f$ will converge on $\D^g_K$, so this
ambiguity of notation will hopefully not cause confusion.

Alternatively, in the calculations in the following sections
involving $\Trop(f)$ and $\New(f)$, we could always work with
$\Trop(f_\rho)$ and $\New(f_\rho)$ for a fixed (small) value of $\rho$.

\paragraph Suppose \label{par:trop.intersect}
that $f_1,\ldots,f_g$ are nonzero power series in $g$ variables
converging on $\D^g_K$.  Clearly if 
$\xi = (\xi_1,\ldots,\xi_g)\in(\D^1_K(\bar K)\setminus\{0\})^g$ is a common root
of the $f_i$ then $\ord(\xi)$ must be contained in the intersection
$\bigcap_{i=1}^g\Trop(f_i)$.  
Let $\sX_i\subset(\D^1_K\setminus\{0\})^g$
be the hypersurface cut out by $f_i$ and let $\sX = \bigcap_{i=1}^g \sX_i$.
We will denote by $\ord$ the map 
$\sX(\bar K)\to\bigcap_{i=1}^g\Trop(f_i)$.

\section{The tropicalizations of $\log_1,\ldots,\log_g$}
\label{sec:thepicture}

\paragraph
Let $\cP$ be a display over $R$ and let $\log = (\log_1,\ldots,\log_g)$ be its
logarithm with respect to some choice of basis, as in
\S\ref{sec:logarithm}.  
In view of Proposition~\ref{prop:kerlog} and~\parref{par:trop.intersect}, 
we would like to draw the tropical hypersurfaces
$\Trop(\log_i)$ and locate their intersection points explicitly.  
However, this is a difficult problem in geometry and combinatorics 
when $g > 2$, so instead we will learn enough about the
structure of $\Trop(\log_i)$ to allow us to perturb the coefficients of
$\log_i$ without changing the number of roots in a suitable
$\B^g_K(\rho)$.  In this spirit, we will study the tropical hypersurfaces
of general $g$-tuples of convergent power series  
\begin{equation}\label{eq:fi}
  f_i(X_1,\ldots,X_g) = \sum_{j=1}^g\sum_{n=0}^\infty a_{n,ij}\, X_j^{p^n} 
  \in K\ps{X_1,\ldots,X_g}\qquad i=1,\ldots,g
\end{equation}
under specific hypotheses on the valuations
$\ord(a_{n,ij})\in\R\cup\{\infty\}$ of the coefficients
$a_{n,ij}$ for $n\geq 0$, $0\leq i,j\leq g$ which will be satisfied by the
perturbations to be considered later on the coefficients of the $\log_i$.
To this end, we fix: 

\begin{genthm}{Data}{true}{}{paragraph}\rmfamily\upshape
  \label{data:data}
  \begin{enum}
  \item  numbers
    $u_{n,ij}\in \R\cup\{\infty\}$ for $n = 0,1,\ldots$ and 
    $i,j = 1,2,\ldots,g$ ($u_{n,ij}$ will be $\ord(a_{n,ij})$);
  \item an integer $N\geq 1$ (to be the level of the canonical subgroup we wish to
    find); 
  \item a nonnegative real number $H < (p-1)/p^N$ (to be the Hasse invariant);
  \item nonnegative real numbers $U_1,\ldots,U_g$ with 
    $\sum_{i=1}^g U_i = H$ ($U_i$ will be $\ord(w_0(\alpha_{ii}))$).
  \end{enum}
\end{genthm}

\subparagraph
Let $N_i$ be the largest integer such that 
\[ \frac{p^{N_i}-1}{p-1}U_i < 1 \]
(if $U_i = 0$ we set $N_i = \infty$). Since
\[ \frac{p^N-1}{p-1}U_i \leq \frac{p^N-1}{p-1} H 
< \frac{p^N-1}{p^N} < 1, \]
we see that $N_i\geq N$ for all $i$.  We require: 

\begin{genthm}{Hypotheses}{true}{}{equation}\rmfamily\upshape 
  \label{hyp:hypotheses}
  \begin{enum}
  \item $u_{0,ii} = 0$ for $i=1,\ldots,g$ and $u_{0,ij} = \infty$ for
    $i\neq j$.
  \item 
    $u_{n,ii} = \frac{p^n-1}{p-1}U_i-n$ for $1\leq n \leq N_i$.
  \item 
    $u_{n,ii}\geq 1 - n$ whenever $n > N_i$.
  \item $u_{n,ij}\geq -n$ for all $1\leq i,j\leq g$ and all $n\geq 0$, 
    and furthermore $u_{n,ij}\geq 1-n$ when $i > j$.
  \end{enum}
\end{genthm}

\subparagraph
Let $f_i$ be any power series of the form~\eqref{eq:fi} such that
the data $(u_{n,ij},N,H,U_i)$ satisfy
Hypotheses~\ref{hyp:hypotheses} with $u_{n,ij} = \ord(a_{n,ij})$.
Note that Hypotheses~(ii) and~(ii) imply Hypothesis~(iv) for $i=j$, and
Hypothesis~(iv) guarantees that the $f_i$ 
converge on the open unit $g$-ball.  As another immediate consequence
of Hypotheses~\ref{hyp:hypotheses}, we have:

\begin{lemsub} Let \label{lem:fromhypotheses}
  $\mc U = (u_{n,ij},N,H,U_i)$ be a set of data 
  satisfying Hypotheses~\ref{hyp:hypotheses}.
  \begin{enum}
  \item If $1\leq N' < N$ then $(u_{n,ij},N',H,U_i)$ is also a set of data
    satisfying Hypotheses~\ref{hyp:hypotheses}.
  \item If $u_{n,ij}' \geq u_{n,ij}$ for all $n\geq 0$ and all $i\neq j$ and
    $u_{n,ii}'=u_{n,ii}$ for all $n\geq 0$ and all $i$ then
    $(u_{n,ij}',N,H,U_i)$ is also a set of data 
    satisfying Hypotheses~\ref{hyp:hypotheses}.
  \item Suppose that $g > 1$.  Choose $i\in\{1,\ldots,g\}$, and let 
    $H' = \sum_{j\neq i} U_j$.  Then 
    \[ (\{u_{n,i'j}\}_{i',j\neq i},N,H',\{U_{i'}\}_{i'\neq i}) \]
    is a set of data satisfying Hypotheses~\ref{hyp:hypotheses} for $g-1$.
  \end{enum}

\end{lemsub}

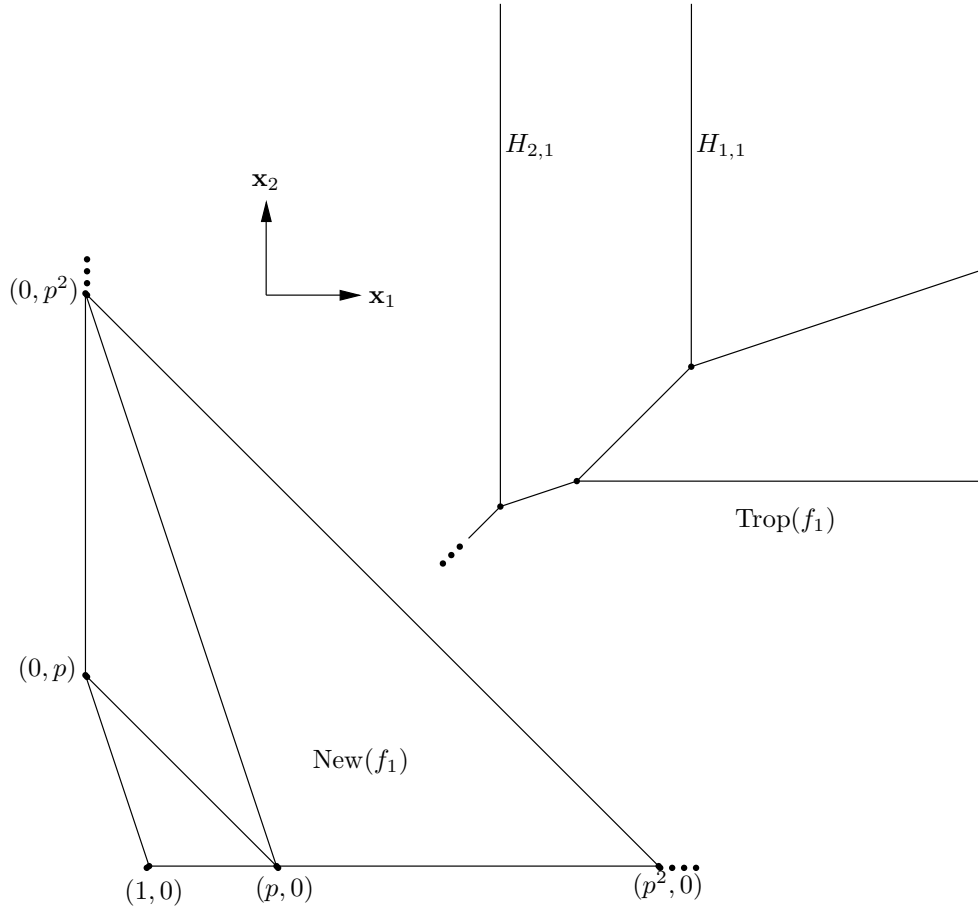
\genericfig[ht]{newtonpoly2}{An example of what finite approximations to
  $\New(f_1)$ and $\Trop(f_1)$ can look like for $f_1$ as
  in~\eqref{eq:fi} when $g = 2$ (the ellipses 
  indicate that the full infinite Newton complex,
  resp.\ tropicalization, continue in the given direction).  Since
  there are no mixed 
  monomial terms in the power series $f_1$, all of the vertices of
  $\New(f_1)$ are contained on the coordinate axes.  
  The infinite Newton complex of the power series $\td f_1$
  of~\eqref{eq:f.i.tilde} is identified with the part of $\New(f_1)$
  contained in the line $\bx_2=0$, as in Proposition~\ref{prop:dim1}.
  The ray $H_{n,1}$ is the dual cell to the line segment
  $\bar{(p^{n-1},0)(p^n,0)}$ in $\New(f_1)$.} 

\paragraph
The first step in understanding $\Trop(f_i)$ is to analyze the 
infinite Newton complex $\New(f_i)$.  
It is clear that all of the vertices of 
$\New(f_i)$ are contained on the coordinate axes, since there are no
mixed monomials in \eqref{eq:fi}.  
Let $\bx_i$ be $i$th standard basis vector in $\R^g$ and let
$L_i = \R\bx_i$ be the $\bx_i$-axis.  Let
\[ \Delta_i = \{ \sigma\cap L_i~:~ \sigma\in\New(f_i)\text{ and }
\sigma\cap L_i\neq\emptyset \}. \]
Then $\Delta_i$ is a collection of line segments and
vertices in $L_i$.  As an exercise one can prove 
that $\Delta_i$ is naturally identified with the infinite Newton
complex of the power series 
\begin{equation} \label{eq:f.i.tilde}
\td f_i(X) = f_i(0,\ldots,X,\ldots,0) = \sum_{n=0}^\infty a_{n,ii} X^{p^n},  
\end{equation}
and that any cell in $\Delta_i$ is a cell in $\New(f_i)$.
Note that by Lemma~\ref{lem:fromhypotheses}(iii), $\td f_i$ satisfies
Hypotheses~\ref{hyp:hypotheses} for $g = 1$.
The following Proposition is illustrated in Figure~\ref{fig:newtonpoly2}.   

\begin{prop} Choose \label{prop:dim1}
  $1\leq i\leq g$ and $n$ such that $0\leq n < N_i$.  
  The line segment joining $p^n\bx_i$ and $p^{n+1}\bx_i$ is a cell
  in $\Delta_i$.
\end{prop}  

\pf
  By the above remarks, we may replace $f_i$ with 
  $f_i(0,\ldots,X_i,\ldots,0)$; hence we may and do assume that
  $g = 1$ and $f = f_1$.  This is now a question about the Newton polygon
  of $f$ (see Example~\ref{eg:newtonpoly}).  Let $u_n = u_{n,11}$ and 
  $U = U_1$.  We need to show that 
  the line passing through the points $(p^n,u_n)$ and $(p^{n+1},u_{n+1})$
  lies strictly below every point $(p^m, u_m)$ with $m\neq n,n+1$, i.e.,
  that the linear form $\omega: (x,y)\mapsto (u_{n}-u_{n+1})x+p^n(p-1)y$
  evaluated at points $(p^m,u_m)$ for $m\geq 0$ achieves its minimum 
  at $(p^n,u_{n})$ and $(p^{n+1},u_{n+1})$.  If $1\leq n < N_i$ then
  \[ u_{n}-u_{n+1} = \frac{p^n-1}{p-1}U - n
  - \frac{p^{n+1}-1}{p-1}U + n+1 = 1 - p^nU, \]
  and if $n = 0$ then
  \[ u_0 - u_1 = 0 - (U - 1) = 1 - p^nU \]
  as well.
  Choose any point $(m,u_{m})$ with $m\neq n,n+1$.  We want the quantity
  \begin{equation}\label{eq:diff}
    \omega(p^m,u_m) - \omega(p^n,u_{n}) 
    = (1-p^nU)(p^m-p^n) + p^n(p-1)(u_{m} - u_n) 
  \end{equation}
  to be positive.  

\subparagraph
  We treat the case $n = 0$ separately.  In this case $m > 1$, and
  \eqref{eq:diff} becomes 
  \[ \omega(p^m,u_m) - \omega(1,0) 
  = (1-U)(p^m-1) + (p-1)u_m. \]
  Suppose that
  $u_m = \frac{p^m-1}{p-1}U - m$.  Then \eqref{eq:diff} is equal to
  \[ (1-U)(p^m-1) + (p^m-1)U - m(p-1) = p^m-1 - m(p-1), \]
  which is positive because $p^m-1  > m(p-1)$ for $m\geq 2$, since
  $\frac{p^m-1}{p-1}=1+p+\cdots+p^{m-1}$.  If
  $u_m\neq\frac{p^m-1}{p-1}U - m$ then $m > N_1$, and by 
  Hypothesis~\ref{hyp:hypotheses}(iii) we have $u_m\geq 1-m$.  Hence
  \eqref{eq:diff} becomes 
  \[ (1-U)(p^m-1) + (p-1)u_m \geq (1-U)(p^m-1) - (p-1)(m-1). \]
  Now $U \leq H < (p-1)/p^N \leq (p-1)/p$, so $1-U\geq 1/p$.  Hence
  \[ (1-U)(p^m-1) + (p-1)u_m 
  \geq \frac 1p(p^m-1) - (p-1)(m-1). \]
  Again one can check that $p^m-1 > p(p-1)(m-1)$ for $m\geq 2$.

\subparagraph
  Now assume $n\geq 1$, and suppose that 
  $u_{m} = \frac{p^m-1}{p-1}U - m.$
  Then \eqref{eq:diff} is equal to
  \begin{equation} \label{eq:diff2} \begin{split}
    (1-p^nU)(p^m-p^n) &+ 
    p^n(p-1)\left(\frac{p^m-1}{p-1}U - \frac{p^n-1}{p-1}U - m+n\right) \\
    &= p^m - p^n + p^n(p-1)(n-m).
  \end{split}\end{equation}
  When $m \geq n+2$ this quantity is positive because
  \[ \frac{p^{m-n}-1}{p-1} = p^{m-n-1} + p^{m-n-2} + \cdots + p + 1 >
  m-n. \]
  When $m < n$, \eqref{eq:diff2} is also positive since
  \[ (p-1)(n-m) > 1 - p^{m-n}. \]

\subparagraph
  If $u_{m} \neq \frac{p^m-1}{p-1}U - m$ then either $m=0$ (so $u_m=0$) or
  else $m>N_1$ (so $m\geq n+2$), in which case $u_{m}\geq 1-m$.
  First suppose that $m>N_1$, so \eqref{eq:diff} is greater than 
  or   equal to
  \[\begin{split}
    (1-p^nU)(p^m-p^n) &- p^n(p^n-1)U + p^n(p-1)(n-m+1) \\
    &= (1-p^nU)p^m - (1-U)p^n + p^n(p-1)(n-m+1) \\
    &> \left(1-\frac{p-1}p\right)p^m - p^n + p^n(p-1)(n-m+1) \\
    &= p^{m-1} - p^n + p^n(p-1)(n-m+1),
  \end{split}\]
  where the second-to-last step comes from the fact that
  \[ U\leq H < \frac{p-1}{p^N}\leq\frac{p-1}{p^{n+1}}. \]
  But $p^{m-1} - p^n + p^n(p-1)(n-m+1)$ is nonnegative because
  \[ \frac{p^{m-n-1}-1}{p-1} \geq m-n-1 \]
  for all $m \geq n+2$.
  If instead $m=0$ then \eqref{eq:diff} is 
  \[\begin{split}
    (1-p^nU)(1-p^n) - p^n(p-1)u_n
    &= (1-p^nU)(1-p^n) - p^n(p-1)\left(\frac{p^n-1}{p-1}U - n\right) \\
    &= 1 + p^n(n(p-1)-1) 
    \geq 1
  \end{split}\]
  because $n(p-1)\geq 1$ when $n\geq 1$.
\qed

\begin{eg}
  Suppose that $g=1$, and let $U = U_1 = H$ and $u_n = u_{n,11}$.  Let
  $N = 1$, and assume that $u_0 = 0,~u_1 = U-1$, and $u_2 = -1$.  If
  the line segment joining $\bx_1$ and $p\bx_1$ is a cell in
  $\Delta=\Delta_1$ then the point $(p,u_1)$ is below the line
  segment joining $(1,u_0)$ and $(p^2,u_2)$; an elementary calculation
  shows that this happens precisely when $U = H < p/(p+1)$.  It is not a
  coincidence that $p/(p+1)$ is Katz's bound on 
  the Hasse invariant of an elliptic curve that suffices for its level-$1$
  canonical subgroup to exist; cf.~\parref{par:katz.bound},
  Remark~\ref{rem:after.mainthm}(ii), and Remark~\ref{rem:g.1.done}.  
\end{eg}

\paragraph
Proposition~\ref{prop:dim1} says that the line segment 
joining $p^{n-1}\bx_i$ and $p^n\bx_i$ is a cell in $\New(f_i)$ when 
$1\leq n \leq N$; let $H_{n,i}$ be its dual cell in $\Trop(f_i)$. 

\begin{lem} Let \label{lem:locate.Hni}
  $1\leq i\leq g$ and let $1\leq n\leq N$ (so $n\leq N_i$).
  \begin{enum}
  \item The cell $H_{n,i}$ is contained in the affine hyperplane defined
    by the equation
    \[ x_i = \frac 1{p^{n-1}(p-1)} - \frac{U_i}{p-1}
    \geq \frac 1{p^{N-1}(p-1)} - \frac H{p-1} > \frac 1{p^N(p-1)}. \]

  \item If $H$ is the dual cell in $\Trop(f_i)$ to any other cell of
    $\Delta_i$ then it is contained in the affine hyperplane defined
    by the equation $x_i = r$ for some $r \leq \frac 1{p^N(p-1)}$.
  \end{enum}

\end{lem}

\pf
  For $m',m\in\Z_{\geq 0}$ with $m'<m$ let
  $I_{m',m}$ be the line segment joining 
  $p^{m'}\bx_i$ and $p^{m}\bx_i$, and when
  $u_{m',ii},u_{m,ii}\neq\infty$ let
  $r_{m',m} = (u_{m',ii}-u_{m,ii})/p^{m'}(p^{m-m'}-1)$, so $-r_{m',m}$ is
  the slope of the line segment joining $(p^{m'},u_{m',ii})$ and
  $(p^m,u_{m,ii})$ in $\R^2$.  (For the first part of the Lemma, we will be
  interested in the case $m' = m-1$ with $1\leq m\leq N$.)
  As explained in~\parref{par:explain.trop},
  if $I_{m',m}$ is a cell in $\Delta_i$ and   
  $x = (x_1,\ldots,x_g)\in\R^g_{>0}$ is contained in its dual cell $H$
  then $p^{m'} x_i + u_{m',ii} = p^m x_i + u_{m,ii}$; i.e., 
  $x_i = r_{m',m}$.  Hence by Example~\ref{eg:newtonpoly}, this is a
  question about the Newton polygon of the power series
  $f_i(0,\ldots,X_i,\ldots,0)$, so as in the proof of 
  Proposition~\ref{prop:dim1} we can assume $g=1$.  Let 
  $\Delta = \Delta_1$, let $u_n = u_{n,11}$ and $U = U_1$.
  Part~(i) now follows immediately from the calculation
  \[ r_{n-1,n} = \frac{u_{n-1}-u_{n}}{p^{n-1}(p-1)}
  = \frac{1-p^{n-1}U}{p^{n-1}(p-1)} 
  = \frac 1{p^{n-1}(p-1)} - \frac{U}{p-1}. \]

  Let $I_{m,m'}$ be a cell in $\Delta$ for some
  $m,m'\in\Z$ such that $N \leq m' < m$, so its dual cell is the point
  $\{r_{m',m}\}$.  If $I_{m'',m'}$ is also a cell in $\Delta$ for some
  $N\leq m'' < m'$ then  
  $-r_{m'',m'} < -r_{m',m}$ (i.e., $r_{m',m} < r_{m'',m'}$) because the
  slopes of the line segments 
  in a Newton polygon are monotonically increasing; thus we may
  assume that $m' = N$.  If $N < N_i$ then $m=N+1$ and 
  $r_{N,N+1} = 1/p^N(p-1) - U/(p-1)$ as above, 
  so we may assume that $N = N_i$ (unless $N_i=\infty$, in which case we
  are done), and hence $u_m\geq 1-m$ by
  Hypothesis~\ref{hyp:hypotheses}(iii). 
  We therefore have (by Hypotheses~\ref{hyp:hypotheses}(ii,iii))
  \[\begin{split} r_{N,m} &= \frac{u_N-u_m}{p^N(p^{m-N}-1)} \\
  &\leq \frac{\frac{p^N-1}{p-1}U - N + m - 1}{p^N(p^{m-N}-1)} \\
  &< \frac{m-N}{p^N(p^{m-N}-1)}
  \qquad\left(\text{using } U\leq H < \frac{p-1}{p^N} \right).
  \end{split}\]
  But this quantity is at most $1/p^N(p-1)$ since for $k\geq 1$ we have
  \[ \frac k{p^N(p^k-1)}\leq\frac 1{p^N(p-1)}, \]
  as one easily checks.
\qed

\begin{rem} When \label{rem:g.1.done}
  $\cP$ be a display of dimension $g = 1$ over $R$, 
  Proposition~\ref{prop:dim1} and Lemma~\ref{lem:locate.Hni} effectively
  calculate the Newton polygon of the logarithm $\log = \log_1$ with
  respect to some choice of basis, as follows.
  We will prove (Proposition~\ref{prop:satisfies_hypotheses}) that if
  $\log = \sum a_n X^n$ then the data $(\ord(a_n),N,H,H)$ 
  satisfy Hypotheses~\ref{hyp:hypotheses}, so
  Proposition~\ref{prop:dim1} shows that for $1\leq n\leq N$ the line
  segment $\sigma_n$ joining $(p^{n-1},\ord(a_{n-1}))$ and
  $(p^n,\ord(a_n))$ is in $\cN(\log)$.  By Lemma~\ref{lem:locate.Hni} the
  slope of $\sigma_n$ is $-(1/p^{n-1}(p-1)-H/(p-1)) < -1/p^N(p-1)$, and the
  slope of any other line segment in $\cN(\log)$ is at least $-1/p^N(p-1)$.
  Hence $\log$ has $p^n-p^{n-1}$ roots with valuation 
  $1/p^{n-1}(p-1) - H/(p-1)$, and this accounts for all nonzero roots of
  $\log$ with valuation at least $r \coloneq 1/p^{N-1}(p-1) - H/(p-1)$.
  Therefore $\log$ has a total of 
  \[ (p^n-p^{n-1}) + (p^{n-1}-p^{n-2}) + \cdots + (p^2-p) + (p-1) + 1 = p^n \]
  roots in the ball $\B^g_K(\rho)$ where $\rho=p^{-r}$.  
  We will also show (Lemma~\ref{lem:pdoesntkill}) that 
  $(\sG[p^\infty]\cap\B^g_K(\rho))(\bar K)\cong\Z/p^n\Z$ 
  (so $\sG[p^\infty]\cap\B^g_K(\rho) = G[p^N]_{\leq\rho}$), which proves
  Theorem~\ref{thm:canonicalsubgp} for $g=1$. 

\end{rem}

\begin{prop} Let \label{prop:diagonals}
  $r = 1/p^N(p-1)$, let $1\leq i\leq g$, and let 
  $x = (x_1,\ldots,x_g)\in\Trop(f_i)\cap\R^g_{>r}$.  Suppose that 
  $x\notin H_{n,i}$ for $n = 1,\ldots,N$ (so $g\geq 2$; cf.
  Remark~\ref{rem:g.1.done}).  Then there exists $j\neq i$ such that 
  \[ x_i - x_j \geq
  \frac {\epsilon_{i>j}}{p^N} - \frac{U_i}{p-1}, \]
  where $\epsilon_{i>j}=1$ if $i>j$ and is $0$ otherwise.
\end{prop}

Proposition~\ref{prop:diagonals} is illustrated in
Figure~\ref{fig:thepicture2}. 

\pf
  To say that $x\in\Trop(f_i)$ means that $\inn_x(f_i)$ is not a monomial.
  Hence there are $m,m'$ and $j,j'$ with $(m,j)\neq(m',j')$ such
  that $p^m x_j + u_{m,ij} = p^{m'}x_{j'} + u_{m',ij'}$, and
  $p^m x_j + u_{m,ij} \leq p^{m''} x_{j''} + u_{m'',ij''}$ for all 
  $m'',j''$.  Since $x_i > r$, either $j\neq i$ or $j'\neq i$, because
  otherwise $x$ would be contained in the dual cell to a cell of
  $\Delta_i$ and hence would lie on some $H_{n,i}$ by
  Lemma~\ref{lem:locate.Hni}(ii).  Therefore, without loss of generality
  we can assume that $j\neq i$.  Then we have
  \[ p^m x_j + u_{m,ij} \leq p^n x_i + u_{n,ii} \]
  for all $n\geq 0$.  We now consider two separate cases.
  \begin{enum}
  \item Suppose that $m\leq N$ (so $m\leq N_i$).  Then 
    $p^m x_j + u_{m,ij} \leq p^m x_i + u_{m,ii}$ implies
    \[\begin{split}
      x_i - x_j &\geq \frac{u_{m,ij}-u_{m,ii}}{p^m} \\
      &\geq \frac{\epsilon_{i>j} -m+m- \frac{p^m-1}{p-1} U_i}{p^m}\qquad
      \text{(by Hypotheses~\ref{hyp:hypotheses}(ii,iii))}\\
      &= \frac{\epsilon_{i>j}}{p^m} -
      \frac{p^m-1}{p^m}\cdot\frac{U_i}{p-1} 
      \geq \frac{\epsilon_{i>j}}{p^m} - \frac{U_i}{p-1}
      \geq \frac{\epsilon_{i>j}}{p^N} - \frac{U_i}{p-1}.
    \end{split}\]

  \item Now suppose that $m > N$.  Then
    $p^m x_j + u_{m,ij} \leq p^N x_i + u_{N,ii}$ implies
    \[\begin{split}
      x_i - x_j 
      &\geq (p^{m-N}-1)x_j + \frac 1{p^N}(u_{m,ij}-u_{N,ii}) \\
      &\geq \frac{p^{m-N}-1}{p^N(p-1)} 
      + \frac{\epsilon_{i>j} - m + N - \frac{p^N-1}{p-1}U_i}{p^N} \\
      &\geq \frac{p^{m-N}-1-(p-1)(m-N)}{p^N(p-1)}
      + \frac{\epsilon_{i>j}}{p^N} - \frac{U_i}{p-1},
    \end{split}\]
    where we have used $x_j > r = 1/p^N(p-1)$.  Hence we need only show that
    for $k = m-N \geq 1$ we have $p^k-1 \geq (p-1)k$, which is easily
    checked.  

  \end{enum}
\qed

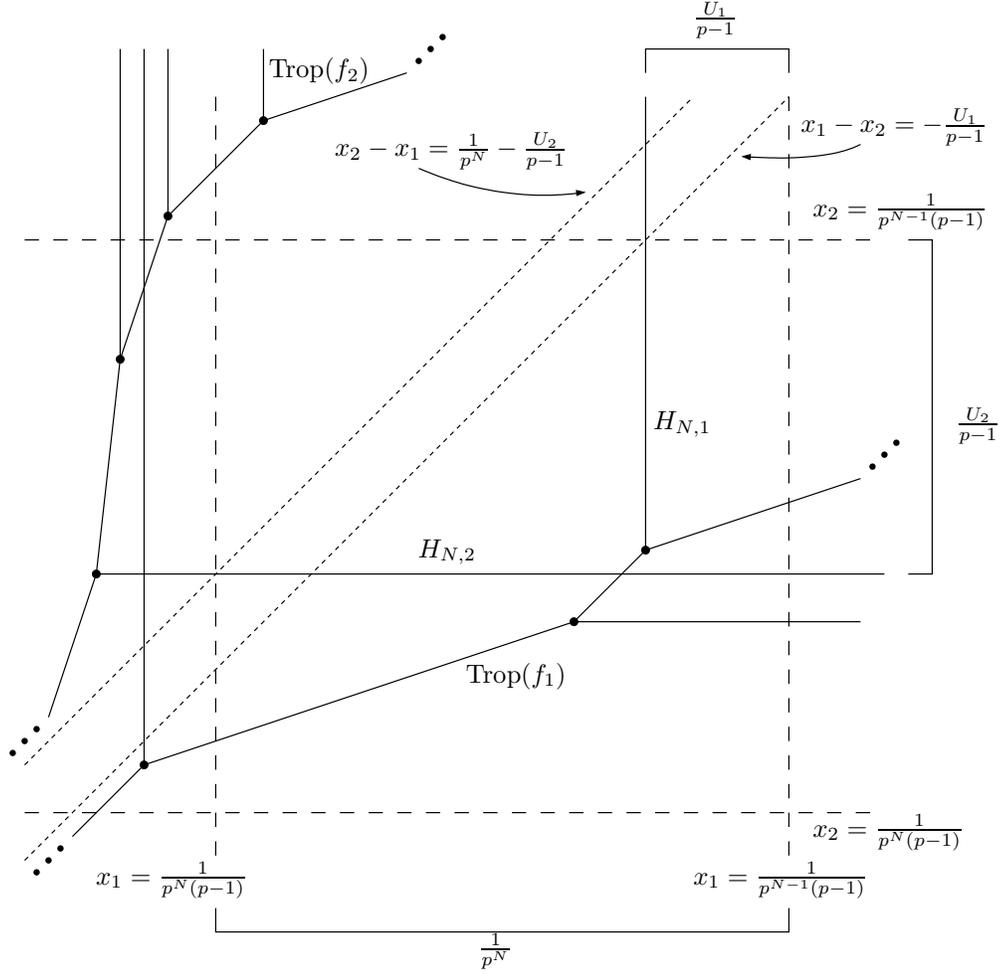
\genericfig[ht]{thepicture2}{This figure illustrates
  Proposition~\ref{prop:diagonals}.  Let $r = 1/p^N(p-1)$.  The vertical
  rays in $\Trop(f_1)\cap\R^2_{>r}$ are the cells
  $H_{n,1}$, and the horizontal rays in $\Trop(f_2)\cap\R^2_{>r}$ are the
  cells $H_{n,2}$, cf.~\parref{par:orthogonal} and Lemma~\ref{lem:locate.Hni}.
  In the region $x_1,x_2 > r$, every point on $\Trop(f_1)$ not 
  contained in $H_{n,1}$ for $1\leq n\leq N$ lies below the dotted line
  $x_1 - x_2 = -U_1/(p-1)$, and any point 
  on $\Trop(f_2)$ not contained in $H_{n,2}$ for $1\leq n\leq N$ lies
  above the dotted line $x_2-x_1 = 1/p^N - U_2/(p-1)$.  
  As the regions $\{(x_1,x_2~:~x_1 - x_2 \geq -U_1/(p-1)\}$ and
  $\{(x_1,x_2)~:~x_2-x_1 \geq 1/p^N - U_2/(p-1)\}$ do not intersect,
  this implies that any
  point of $\Trop(f_1)\cap\Trop(f_2)\cap\R^2_{>r}$ is
  located on some ray $H_{n,1}$ or $H_{n,2}$, as in the proof of
  Theorem~\ref{thm:properness}.  (Note that many such intersection points
  will lie above or to the right of the Figure.)
}

\paragraphnn
The main result concerning the structure of $\Trop(f_i)$ is as follows.

\begin{thm} Let \label{thm:properness}
  \[ r = \frac 1{p^N(p-1)}\quad\text{ and }\quad
  r' = \frac 1{p^{N-1}(p-1)} - \frac{H}{p-1}, \]
  so $r' > r$ since $H < (p-1)/p^N$.
  Then 
  \[ \bigcap_{i=1}^g\Trop(f_i)\cap\R^g_{>r} 
  = \bigcap_{i=1}^g\Trop(f_i)\cap\R^g_{\geq r'} \quad\text{ and }\quad
  \bigcap_{i=1}^g\Trop(f_i)\cap\R^g_{>1/(p-1)} = \emptyset. \]
\end{thm}

\pf
  Let $x = (x_1,\ldots,x_g)\in\bigcap_{i=1}^g\Trop(f_i)\cap\R^g_{>r}$.
  Choose $1\leq i_1\leq g$; we want to show that $x_{i_1}\geq r'$.
  If $x\in H_{m,i_1}$ for some $1\leq m\leq N$ then $x_i\geq r'$ by
  Lemma~\ref{lem:locate.Hni}, so we may restrict our attention to the case
  $x\notin\bigcup_{m=1}^N H_{m,i_1}$.  By
  Proposition~\ref{prop:diagonals}, there exists $1\leq i_2\leq g$ with
  $i_2\neq i_1$ such that  
  \[ x_{i_1} - x_{i_2} \geq 
  \frac{\epsilon_{i_1>i_2}}{p^N} - \frac{U_{i_1}}{p-1}. \]
  Continuing in this fashion, there are pairwise distinct integers 
  $i_1,i_2,\ldots,i_n\in\{1,\ldots,g\}$ for some $2\leq n\leq g$ such
  that for each $j<n$ we have 
  \begin{equation}  \tag{*}
    x\notin\bigcup_{m=1}^N H_{m,i_j} \quad\text{ and }\quad
    x_{i_j} - x_{i_{j+1}} \geq 
    \frac{\epsilon_{i_j>i_{j+1}}}{p^N} - \frac{U_{i_j}}{p-1}. 
  \end{equation}
  We may assume that there is no $i_{n+1}\notin\{i_1,\ldots,i_n\}$
  such that $i_1,\ldots,i_{n+1}$
  satisfies~(*) for all $j < n+1$, which is to say that 
  either (i) there exists an $i_{n+1} = i_\ell$ for some $\ell < n$
  satisfying~(*) for $j=n$, 
  or (ii) $x\in\bigcup_{m=1}^N H_{m,i_n}$.  We claim that case~(i) cannot
  happen.  If~(i) occurs then
  \[ 0 = \sum_{j=\ell}^n (x_{i_j} - x_{i_{j+1}})
  \geq \frac{\sum_{j=\ell}^n\epsilon_{i_j>i_{j+1}}}{p^N} - 
  \frac{\sum_{j=\ell}^n U_{i_j}}{p-1}. \]
  Since $i_\ell = i_{n+1}$, at least one $\epsilon_{i_j>i_{j+1}}$ is $1$,
  and since $\sum_{j=1}^g U_j = H$ we have 
  $\sum_{j=\ell}^n U_{i_j}\leq H$.  Hence
  \[ 0 \geq \frac 1{p^N} - \frac H{p-1} > 0, \]
  a contradiction.  Therefore $x\in H_{m,i_n}$ for some
  $1\leq m \leq N$, so
  \[ x_{i_n} = \frac 1{p^{m-1}(p-1)} - \frac{U_{i_n}}{p-1}
  \geq \frac 1{p^{N-1}(p-1)} - \frac{U_{i_n}}{p-1}. \]
  Hence
  \[\begin{split}
    x_{i_1} &= (x_{i_1} - x_{i_n}) + x_{i_n} 
    = \sum_{j=1}^{n-1} (x_{i_j} - x_{i_{j+1}}) + x_{i_n} \\
    &\geq \frac{\sum_{j=1}^{n-1}\epsilon_{i_j>i_{j+1}}}{p^N} 
    - \frac{\sum_{j=1}^{n-1} U_{i_j}}{p-1} 
    + \frac 1{p^{N-1}(p-1)} - \frac{U_{i_n}}{p-1} \\
    &\geq \frac 1{p^{N-1}(p-1)} 
    - \frac{\sum_{j=1}^n U_{i_j}}{p-1}
    \geq \frac 1{p^{N-1}(p-1)} 
    - \frac H{p-1} = r'.
  \end{split}\]

  For the final assertion, since 
  \[ \frac 1{p-1} \geq \frac 1{p^{m-1}(p-1)} - \frac {U_i}{p-1} \]
  for $1\leq m\leq N$, any point 
  $x\in\bigcap_{i=1}^g\Trop(f_i)\cap\R^g_{>1/(p-1)}$ is not contained on
  any hyperplane $H_{m,i}$ by Lemma~\ref{lem:locate.Hni}(i).  But we
  showed above that this is impossible. 
\qed

\section{Counting the common roots of $(\log_1, \ldots, \log_g)$}
\label{sec:counting}

\paragraph
Let $\cP = (P,Q,F,V\inv)$ be a display over $R$ and let $P = T\dsum L$ be a
normal decomposition, let $e_1,\ldots,e_g$ and $e_{g+1},\ldots,e_h$ be
$W(R)$-bases for $T$ and $L$, respectively, and let 
$M = (\alpha_{ij})$ be the structure matrix for $\cP$ with respect to this
basis, as in~\eqref{eq:alpha.ij}.  Let $G = \BT_\cP$ be the associated
$p$-divisible formal group over $R$ (of dimension $g$ and height $h$),
and let $\sG$ be its rigid generic fiber.
Let $\log = (\log_1,\ldots,\log_g)$ be the logarithm of $G$, calculated
in terms of $\cP$ using \eqref{eq:fundeq}. 

\begin{notn} We \label{notn:r.n.rho.n}
  will use the following notation in this section and the next. For 
  $n\geq 0$, define 
  \[ r_n = \frac 1{p^n(p-1)} \qquad\text{ and }\qquad
  r_n' = \frac 1{p^{n-1}(p-1)} - \frac H{p-1}, \]
  where $H = H(G)$.  We also set $\rho_n = p^{-r_n}$ and 
  $\rho_n' = p^{-r_n'}$.  
  Note that if $H < (p-1)/p^n$ then $r_n < r_n'$ and $\rho_n > \rho_n'$.

  If $\sX\subset\D^g_K$ is any analytic subspace and 
  $\rho\in\sqrt{|K^\times|},~0<\rho\leq 1$, we define
  \begin{equation*} 
    \sX_{\leq\rho} \coloneq \sX\cap\B^g_K(\rho) 
    \quad\text{ and }\quad
    \sX_{<\rho} \coloneq \sX\cap\D^g_K(\rho).
  \end{equation*}
\end{notn}

\paragraph
We will use the results of \S\ref{sec:thepicture} and a continuity
of roots argument to count the number of points 
in the set $\ms G[p^\infty]_{\leq\rho_n'}(\bar K)$ when $H < (p-1)/p^n$;
i.e., to count the number of common zeros of $(\log_1,\ldots,\log_g)$ whose
coordinates have absolute value at most $\rho_n'$.  (It will turn out that
$\sG[p^\infty]_{\leq\rho_n'} = G[p^n]_{\leq\rho_n'}$ when
$H(G)<(p-1)/p^n$.)  Specifically, we will prove:

\begin{prop} Fix \label{prop:countroots}
  an integer $N\geq 1$, and assume that $H(G) < (p-1)/p^N$.
  For $0 < n \leq N$, the group $\ms G[p^\infty]_{\leq\rho_n'}(\bar K)$ has
  exactly $p^{ng}$ points. 
\end{prop}

\paragraph Recall \label{par:fix.data}
that we have defined $g\times h$ matrices $\ba_n = (a_{n,ij})$ with
entries in $K$ for $n\geq 0$, calculated 
in terms of $M$ using the equations~\eqref{eq:fundeq}, such that 
\[ \log_i(X_1,\ldots,X_g) = \sum_{j=1}^g\sum_{n=0}^\infty a_{n,ij}\, X_j^{p^n} 
\in K\ps{X_1,\ldots,X_g}\qquad i=1,\ldots,g. \]
Our first task is to show that the coefficients $a_{n,ij}$ of
$\log_1,\ldots,\log_g$  satisfy 
Hypotheses~\ref{hyp:hypotheses}.  To this end, for the rest of this
section we \emph{fix an integer} $N\geq 1$, and we \emph{assume that} 
$H = H(G) < (p-1)/p^N$.  We also let
$U_i = \ord(w_0(\alpha_{ii}))$, so $\sum_{i=1}^g U_i = H$
(cf.~\parref{par:assume.uppertriang}).  Let  
$u_{n,ij} = \ord(a_{n,ij})\in\R\cup\{\infty\}$ for
$i=1,2,\ldots,g$ and $j = 1,2,\ldots,h$. 
Note that $(\{u_{n,ij}\}_{i,j=1,\ldots,g},N,H,U_i)$ forms a system of data as
in~\parref{data:data}. 

\begin{prop} The \label{prop:satisfies_hypotheses}
  data $(\{u_{n,ij}\}_{i,j=1,\ldots,g},N,H,U_i)$ satisfy 
  Hypotheses~\ref{hyp:hypotheses}.
\end{prop}

\pf
  Hypothesis~(i) is clear by \eqref{eq:fundeq}.  We will prove
  Hypotheses~(ii)--(iv) by induction on $n$, the base case $n=0$ being
  Hypothesis~(i).  We will also need the following inductive hypothesis:
  \begin{equation}\tag{*} \qquad u_{n,ij} \geq -n+1 \text{ for } j > g, 
  \end{equation}
  which is satisfied for $n = 0$.

\subparagraph
  Suppose that Hypotheses~(ii)--(iv) and~(*) are satisfied for
  some $n \geq 0$.  Express the $h\times h$ structure matrix $M$ in
  block-matrix form: 
  \[ M = 
  \begin{bmatrix}
    A & B \\ C & D
  \end{bmatrix}, \]
  where $A$ is a $g\times g$ matrix, $D$ is an $(h-g)\times(h-g)$ matrix, 
  $C$ is an $(h-g)\times g$ matrix, and $B$ is a $g\times (h-g)$ matrix.  We
  also write $\ba_n$ in the block form
  \[ \ba_n = 
  \begin{bmatrix} 
    \ba_n'  & \bb b_n
  \end{bmatrix}, \]
  where $\ba_n'$ is a $g\times g$ matrix
  and $\bb b_n$ is a $g\times (h-g)$ matrix. 
  For brevity, if $E$ is a matrix with coefficients in $R$, we will write
  $\ord(E)\geq\delta$ to mean that the valuation of every coefficient of
  $E$ is at least $\delta$.  Hence
  Hypotheses~(ii)--(iv) and~(*) imply
  \[ \ord(\ba_n') \geq -n,\quad
  \ord(\bb b_n) \geq -n+1. \]

\subparagraph
  First we prove~(*) and Hypothesis~(iv) for $n+1$.
  Writing out \eqref{eq:fundeq}, we have
  \[ 
  \begin{bmatrix} 
    \ba_{n+1}'  & \bb b_{n+1}
  \end{bmatrix} =
  \begin{bmatrix}
    p\inv\big(\ba_n'w_n(A) + \bb b_nw_n(C)\big) &
    \ba_n'w_n(B) + \bb b_n w_n(D)
  \end{bmatrix}.
  \]
  Since $\ord(\ba_n')\geq -n$, $\ord(\bb b_n)\geq 1-n$, and
  $\ord(w_n(M))\geq 0$, we see that 
  $\ord(\bb b_{n+1}) \geq -n = 1-(n+1)$, which proves (*) for $n+1$.
  The same argument shows that
  $\ord(\ba_{n+1}') \geq -(n+1)$, i.e.,
  $u_{n+1,ij}\geq -(n+1)$ for all $i,j$.  Since $\ord(\bb b_n)\geq 1-n$,
  we have
  \[ a_{n+1,ij} = \frac 1p\left(\sum_{\ell=1}^g 
    a_{n,i\ell}w_n(\alpha_{\ell j})\right) + 
  (\text{terms with valuation at least $-n$}). \]
  Hypothesis~(iv) for $n$ says that $\ord(a_{n,i\ell})\geq 1-n$ when 
  $i > \ell$, and since $w_0(A)$ is upper-triangular modulo $p$, 
  $w_n(\alpha_{\ell j})$ is
  divisible by $p$ for $\ell > j$, so $\ord(w_n(\alpha_{\ell j}))\geq 1$.
  If $i > j$ then either $i > \ell$ or $\ell > j$, so 
  $u_{n+1,ij}\geq 1-(n+1)$ in this case.  This proves Hypothesis~(iv) for $n+1$.
  
\subparagraph
  It remains to prove Hypotheses~(ii) and~(iii) for $n+1$.  We have
  \[ a_{n+1,ii} = \frac 1p\left(
    \sum_{\ell=1}^ga_{n,i\ell}w_n(\alpha_{\ell i}) 
  \right) 
  + (\text{terms with valuation at least $-n$}), \]
  with $\ord(a_{n,i\ell}w_n(\alpha_{\ell i}))\geq 1-n$ when 
  $\ell\neq i$ as above, since $\ord(a_{n,i\ell}) \geq 1-n$ when
  $i>\ell$ by Hypothesis~(iv) and $\ord(w_n(\alpha_{\ell i}))\geq 1$ when
  $\ell > i$.  Assume that $n+1\leq N_i$ (i.e. we are in the situation of
  Hypothesis~(ii)), so
  \[ \frac{p^{n+1} - 1}{p-1} U_i < 1 \quad\text{ and by induction, }\quad
  \ord(a_{n,ii}) = \frac{p^n-1}{p-1} U_i - n. \]
  Write $\alpha_{ij} = (\alpha_{ij,0},\alpha_{ij,1},\dots)$ in its Witt
  coordinates, so $U_i = \ord(\alpha_{ii,0})$.  Then
  \[ p^n\ord(\alpha_{ii,0}) = p^n U_i < \frac{p^{n+1}-1}{p-1}U_i < 1 \]
  so
  \[ \ord(w_n(\alpha_{ii})) 
  = \ord\left(\sum_{j=0}^n p^j \alpha_{ii,j}^{p^{n-j}}\right)
  = \ord\big(\alpha_{ii,0}^{p^n}\big) = p^n U_i. \]
  Hence
  \[ \ord(a_{n,ii}w_n(\alpha_{ii})) 
  = \frac{p^n-1}{p-1}U_i + p^nU_i - n
  = \frac{p^{n+1}-1}{p-1}U_i - n < 1-n. \]
  We conclude that
  \[ u_{n+1,ii} = \ord(a_{n+1,ii}) = \frac{p^{n+1}-1}{p-1}U_i - (n+1), \]
  which proves Hypothesis~(ii) for $n+1$.

  Now assume $n+1 > N_i$ (i.e. we are in the situation of
  Hypothesis~(iii)), so
  \[ \frac{p^{n+1} - 1}{p-1} U_i \geq 1. \] 
  The formula $w_n(\alpha_{ii}) = \sum_{j=0}^n p^j\alpha_{ii,j}^{p^{n-j}}$
  gives $\ord(w_n(\alpha_{ii}))\geq\min\{p^nU_i,1\}$ since
  $U_i=\ord(\alpha_{ii,0})$, so if
  $n = N_i$ then by Hypothesis~(ii) for $n$, 
  \[ \ord(a_{n,ii}w_n(\alpha_{ii,0})) 
  \geq \frac{p^n-1}{p-1}U_i - n + \min\{p^nU_i,1\}
  = \min\left\{\frac{p^{n+1}-1}{p-1}U_i,~1\right\}-n\geq 1-n.
  \]
  If $n > N_i$ then by Hypothesis~(iii) for $n$,
  \[ \ord(a_{n,ii}w_n(\alpha_{ii,0})) \geq
  1-n + \ord(w_n(\alpha_{ii,0})) \geq 1-n.
  \]
  Therefore, $u_{n+1,ii}\geq -n = 1-(n+1)$, which proves Hypothesis~(iii)
  for $n+1$.

\qed

\paragraph
As mentioned above, we will perturb the coefficients of the power series
$\log_i$ to simplify the combinatorics.  The following kind of power
series will serve as the perturbations (recall~\parref{par:fix.data} that
we have fixed $N, H,$ and $\{U_i\}_{i=1,\ldots,g}$).

\begin{defn} Let \label{def:hypotheses}
  \[ f_i(X_1,\ldots,X_g) 
  = \sum_{j=1}^g\sum_{n=0}^\infty b_{n,ij}\, X_j^{p^n}
  \in K'\ps{X_1,\ldots,X_g} \]
  be power series with coefficients in a finite extension field $K'$ of
  $K$ for $i=1,\ldots,g$.  We say that $(f_1,\ldots,f_g)$ 
  \emph{satisfies Hypotheses~\ref{hyp:hypotheses}} if the data
  $(\ord(b_{n,ij}),N,H,U_i)$ satisfy
  Hypotheses~\ref{hyp:hypotheses}.

\end{defn}

\paragraphnn
The following Lemma is basically a translation of
Theorem~\ref{thm:properness}.  We use the notation of~\parref{notn:r.n.rho.n}.

\begin{lem} Let \label{lem:proper}
  $K'$ be a finite extension of $K$, 
  let $f_1,\ldots,f_g\in K'\ps{X_1,\ldots,X_g}$ be power series satisfying
  Hypotheses~\ref{hyp:hypotheses}, and let $\sX\subset\D^g_{K'}$ be the analytic
  subspace defined by $(f_1,\ldots,f_g)$.  Then for $1\leq n\leq N$ we have
  $\sX\cap\D^g_{K'}(\rho_n) = \sX\cap\B^g_{K'}(\rho_n')$.  
  In addition, $\sX\cap\D^g_{K'}(\rho_0) = \{0\}$.

\end{lem}

\pf
  By decreasing $N$ we may assume that $n = N$;
  cf. Lemma~\ref{lem:fromhypotheses}(i).  The case $g=1$ is handled by
  Remark~\ref{rem:g.1.done}, so assume by induction on $g$ that $g>1$.
  Let $\xi = (\xi_1,\ldots,\xi_g)\in \sX\cap\D^g_{K'}(\rho_N)$. 
  If all $\xi_i$ are nonzero then
  \[ \ord(\xi) = (\ord(\xi_1),\ldots,\ord(\xi_g))
  \in\bigcap_{i=1}^g\Trop(f_i)\cap\R^g_{> r_N}, \] 
  so by Theorem~\ref{thm:properness}, $\xi\in\B^g_{K'}(\rho_N')$.
  Now suppose that $\xi_i = 0$ for some $i$.  Let
  \[ \td f_j(X_1,\ldots,\hat X_i,\ldots,X_g)
  = f_j(X_1,\ldots,X_{i-1},0,X_{i+1},\ldots,X_g) \]
  for $j\neq i$.  Then by Lemma~\ref{lem:fromhypotheses}(iii), 
  the power series $\{\td f_j\}_{j\neq i}$ satisfy
  Hypotheses~\ref{hyp:hypotheses} for 
  $g-1$.  But $(\xi_1,\ldots,\hat\xi_i,\ldots,\xi_g)$ is a root of each
  $\td f_j$, $j\neq i$, 
  so by induction on $g$ we have $\xi\in\B^g_{K'}(\rho_N')\cap\{\xi_i=0\}$. 

  The same argument shows that $\sX\cap\D^g_{K'}(\rho_0) = \{0\}$.
\qed

\subparagraph
In particular, taking $(f_1,\ldots,f_g) = (\log_1,\ldots,\log_g)$ and
using Propositions~\ref{prop:kerlog} and~\ref{prop:satisfies_hypotheses},
we see that for all $1\leq n\leq N$,
\begin{equation} \label{eq:cor.proper}
  \sG[p^\infty]_{\leq\rho_n'} = \sG[p^\infty]_{<\rho_n} \quad\text{ and }\quad
  \sG[p^\infty]_{<\rho_0} = \{0\}.
\end{equation}

\paragraphnn
Lemma~\ref{lem:proper} and the following proposition allow us to
make drastic modifications to our power series without affecting the
number of roots in $\B^g_K(\rho_n')$.

\begin{prop}[Continuity of Roots] Let \label{prop:controots1}
  $\rho,\mu\in\sqrt{|K^\times|}$ be such that $\rho < \mu$.
  Let $\sX\subset\B^g_K(\mu)\times\B^1_K$ be a closed analytic
  subspace contained in $\B^g_K(\rho)\times\B^1_K$, and let
  $\phi:\sX\to\B^1_K$ be projection onto the second factor.  Suppose that
  $\sX$ is defined by an ideal of the form
  $(f_1,\ldots,f_g)\subset T_{g,\mu}\angles{t}$.
  Then $\phi$ is a finite, flat map.  In particular, each fiber of $\sX$
  over $\B^1_K$ is finite and has the same length.
\end{prop}

\pf
  It is harmless to make a finite extension of $K$, so we can reduce to
  the case $\mu=1$.
  Let $A = T_g\angles{t}/\fa$, where $\fa$ is the ideal
  generated by $f_1,\ldots,f_g$, so $\sX = \Sp(A)$.
  First we show that $A$ is a finite $K\angles t$-module.
  By the rigid-analytic direct image
  theorem~\cite[Theorem~9.6.3/1]{bgr:nonarch}, we need only show 
  that $\phi$ is proper.  Using Kiehl's definition of properness for rigid
  spaces, it suffices to show 
  that there are affinoid generators $h_1,\ldots,h_n$ of $A$ over
  $K\angles{t}$ such that $|h_i|_{\sup} < 1$.  Since
  $\sX\subset\B^g_K(\rho)\times\B^1_K$ with $\rho<1$, the
  standard generators of $T_g$ satisfy this property.

  Since $\fa$ has $g$ generators, by Krull's principal ideal theorem
  every irreducible component of $\Spec(A)$ 
  has dimension at least one, so since the fibers of $\phi$ have dimension
  zero, $\sX$ is pure dimension $1$ by \cite[Theorem~15.1]{matsumura:crt}.
  Since $T_{g,\mu}$ is a regular ring, it is Cohen-Macaulay, so by the
  unmixedness theorem \cite[Theorem~17.6]{matsumura:crt}, $\Spec(A)$ has no
  embedded points.  Thus every associated point of $\Spec(A)$ lies over the
  generic point of $\Spec(K\angles t)$.  Since $K\angles t$ is Dedekind,
  $A$ is a flat $K\angles t$-module. 
\qed

\begin{defn}
  Let $K'$ be a finite extension of $K$,
  let $(f_1,\ldots,f_g)$ be a $g$-tuple of power series in $g$ variables
  converging on $\D^g_{K'}$, and choose $\rho\in\sqrt{|K^\times|}$, 
  $0 < \rho < 1$.  By the  
  \emph{number of roots} of $(f_1,\ldots,f_g)$ in $\B^g_{K'}(\rho)$ we mean
  the $K'$-dimension of $T_{K',g,\rho}/(f_1,\ldots,f_g)$ (which may be
  infinite).
\end{defn}

When $(f_1,\ldots,f_g)$ cut out an \'etale subspace $\sX\subset\B^g_{K'}(\rho)$,
the number of roots of $(f_1,\ldots,f_g)$ in $\B^g_{K'}(\rho)$ is simply
the number of geometric points of $\sX$.

\begin{cor} Let \label{cor:deform}
  $f_1,\ldots,f_g$ be analytic functions on $\D^g_K\times\B^1_K$ such that
  for each $t_0\in\B^1_K$, the specializations
  $f_{1,t_0},\ldots,f_{g,t_0}$ satisfy
  Hypotheses~\ref{hyp:hypotheses} over $\kappa(t_0)$ (in the sense
  of~\parref{def:hypotheses}).  Then for any 
  $1\leq n\leq N$, $(f_{1,t_0},\ldots,f_{g,t_0})$ has a finite number of roots
  in $\B^g_K(\rho_n')$, and this number is independent of $t_0\in\B^1_K$.
\end{cor}

\pf
  Decreasing $N$ if necessary, we may assume that $n = N$.
  Let $\sX\subset\D^g_K\times\B^1_K$ be the analytic subspace cut
  out by $f_1,\ldots,f_g$, and for $t_0\in\B^1_K$ let $\sX_{t_0}$ be the fiber
  over $t_0$.  By Lemma~\ref{lem:proper},
  $\sX_{t_0}\subset\B^g_{\kappa(t_0)}(\rho_N')$ for all $t_0$, so 
  $\sX\subset\B^g_K(\rho_N')\times\B^1_K$.  The result now follows from
  Proposition~\ref{prop:controots1} with $\rho = \rho_N'$ and
  any $\mu\in(\rho_N',1)\cap\sqrt{|K^\times|}$. 
\qed

\paragraph
\pf[of Proposition~\ref{prop:countroots}]
  We define
  \[ f_i(X_1,\ldots,X_g; t)
  = \sum_{n=0}^\infty a_{n,ii} X_i^{p^n} + 
  t\sum_{j\neq i}\sum_{n=0}^\infty a_{n,ij}\, X_j^{p^n}. \]
  By Lemma~\ref{lem:fromhypotheses}(ii), for every 
  $t_0\in\B^1_K$ the specializations $f_{1,t_0},\ldots,f_{g,t_0}$ satisfy
  Hypotheses~\ref{hyp:hypotheses}.  Thus by Corollary~\ref{cor:deform}, for
  every $1\leq n\leq N$, 
  $(f_{1,0},\ldots,f_{g,0})$ has the same number of roots in 
  $\B^g_K(\rho_n')$ as $(f_{1,1},\ldots,f_{g,1})=(\log_1,\ldots,\log_g)$.  Let
  \[ h_i = f_{i,0} = \sum_{n=0}^\infty a_{n,ii} X_i^{p^n}; \]
  we want to show that $(h_1,\ldots,h_g)$ has $p^{ng}$ roots in
  $\B^g_K(\rho_n')$.  Thinking of $h_i$ as a power series in one variable,
  let $\sX_i\subset\B^1_K(\rho_n')$ be the subspace cut out by $h_i$, so 
  the ideal $(h_1,\ldots,h_g)$ cuts out the product
  $\sX_1\times\cdots\times\sX_g\subset\B^g_K(\rho_n')$.  Thus it suffices
  to show that the length of $\sX_i$ is $p^n$ for $i=1,\ldots,g$, so we
  are reduced to the case $g = 1$.  This is handled by
  Remark~\ref{rem:g.1.done}.
\qed

\section{The group structure on $\ms G[p^\infty]_{\leq\rho_N'}(\bar K)$}
\label{sec:groupstruct}

\paragraph
We keep the notation of \S\ref{sec:counting}.
The only remaining step in the proof of Theorem~\ref{thm:canonicalsubgp}
is to show that the group $\ms G[p^\infty]_{\leq\rho_N'}(\bar K)$ of size
$p^{Ng}$ is in fact equal to $\ms G[p^N]_{\leq\rho_N'}(\bar K)$, and is
isomorphic to $(\Z/p^N\Z)^g$. 

\paragraph
Let $\bX=(X_1,\ldots,X_g)$ and let $\fG\cong\Spf(R\ps\bX)$ be the
canonical formal Lie group with $G\cong\fG[p^\infty]$, as
in~\parref{par:G.p.infty}.  Let 
$[p](\bX)\in(\bX)R\ps{\bX}^g$ be the $g$-tuple of power
series defining the multiplication-by-$p$ map on $\fG$.  
Let $R_0=R/pR$, let $\fG_0 = \fG\tensor_R R_0$ and $G_0=G\tensor_R R_0$, and 
let $V$ be the relative Verschiebung map $\fG_{0}^{(p)}\to\fG_{0}$ over 
$R_0$.  We have chosen a basis for the tangent space $\Lie(G_0)$ of
$G_0$, which induces a basis for
$\Lie(G_0^{(p)})\cong\Lie(G_0)\tensor_{R_0,\Frob} R_0$; with respect to
these bases, we can view
$\dV:\Lie(G_0^{(p)})\to\Lie(G_0)$ as a $g\times g$ matrix with
coefficients in $R_0$. 
Choose a lift of $\td\dV\in\M_g(R)$ of this matrix.
The following is a slight generalization
of a result in \cite[\S3.6]{katz:padicmfs}.

\begin{lem} We have \label{lem:pX} 
  \[ [p](\bX) = p\bX + \td\dV(X_1^p,\ldots,X_g^p) + pf(\bX) 
  + O(\bX^{p^2}), \]
  where $f$ is a $g$-tuple of power series with no terms of total degree
  less than $p$. 
\end{lem}

The proof is left to the reader; one uses the fact that for $\zeta\in R$ a
primitive $(p-1)$th root of unity, we have
$\log(\zeta\bX) = \zeta\log(\bX)$, so 
\[\begin{split}
  [p](\zeta\bX) &= \exp(p\cdot\log(\zeta\bX))  \\
  &= \exp(\zeta p\cdot\log(\bX)) = \exp(\zeta\log([p]\bX)) \\
  &= \exp(\log(\zeta[p](\bX))) = \zeta[p](\bX).
\end{split}\]
(In fact Lemma~\ref{lem:pX} is true for any formal Lie group $\fG$ equipped
with a choice of parameters such that the
total degree of any monomial in $\log(\bX)$ is a power of $p$.)

\paragraph
By Remark~\ref{rem:after.hasse}(iv), the matrix for $\dV$ 
(with respect to our choice of basis) is given 
$w_0(A)\mod p$ where $A$ is the upper-left $g\times g$
submatrix of the structure matrix $M$ of our display $\cP$.
Recall~\parref{par:assume.uppertriang} that we are assuming that $w_0(A)$
is upper-triangular mod $p$, so $\td\dV\pmod p$ is an
upper-triangular matrix.

\begin{notn}
  For $\xi = (\xi_1,\ldots,\xi_g)\in\B^g_K(\bar K)\setminus\{0\}$ let 
  $\size(\xi) = \min\{\ord(\xi_i)\}_{i=1}^g$ and for $n\in\Z_{\geq 0}$ let
  $\xi^{(n)} = (\xi_1^n,\ldots,\xi_g^n)$.
\end{notn}

In other words, if $\xi\in\B^g_K(\bar K)\setminus\{0\}$ and $r=\size(\xi)$ then 
$\rho = p^{-r}$ is the radius of the smallest ball $\B^g_K(\rho)$
containing $\xi$. 

\begin{rem} Lemma~\ref{lem:pX} \label{rem:all.cs.ordinary}
  implies that if $\xi\in\D^g_K(\bar K)\setminus\{0\}$ then
  \[ \size([p]\xi) \geq \min\{p\size(\xi),~1+\size(\xi)\}, \]
  or equivalently,
  \[ \size(\xi)\leq\max\{p\inv\size([p]\xi),~\size([p]\xi)-1\}. \]
  Consequently, if $\xi_1,\xi_2,\ldots\in\D^g_K(\bar K)$ satisfy
  $\xi_1\neq 0$ and $[p]\xi_{n+1} = \xi_{n}$ for all $n\geq 1$, then
  $\size(\xi_n)\to 0$ as $n\to\infty$.
 
  Suppose that the level-$n$ canonical subgroup $G_n$ of $G$ exists for
  all $n$, and that $G_n$ admits the radius $\mu_n$.
  Then $[p]:G_{n+1}(\bar K)\to G_n(\bar K)$ is a surjection, so 
  there exist $\xi_n\in G_n(\bar K) = G[p^n]_{\leq\mu_n}(\bar K)$ for 
  $n\geq 1$ with $\size(\xi_n)\to 0$.  Hence $\mu_n\to 1$, so
  $\sG[p^\infty] = \bigcup_{n\geq 1} G_n\tensor_R K$.
  This implies that $G$ is ordinary, since 
  $G[p](\bar K) \subset\bigcup_{n\geq 1} G_n[p](\bar K) = G_1(\bar K)$.
\end{rem}

\begin{lem} Choose \label{lem:pdoesntkill}
  $\xi\in\B^g_K(\bar K)\setminus\{0\}$ and suppose that for some integer
  $n\geq 1$ we have  
  \[ \frac 1{p^{n+1}(p-1)} < \size(\xi) \leq \frac 1{p(p-1)}. \]
  If $H = H(G) < (p-1)/p^{n+1}$ then $[p](\xi)\neq 0$.
\end{lem}

\pf
  First we will show by induction on $g$ that if $B = (b_{ij})\in\M_g(R)$
  is any upper-triangular matrix with nonzero diagonal entries and
  $\nu\in\B^g_K(\bar K)\setminus\{0\}$ is any point then 
  \begin{equation} \tag{*} \size(B\nu) 
  = \min\left\{\ord\left(\sum_{j\geq i}b_{ij}\nu_j\right)\right\}_{i=1}^g
  \leq h + \size(\nu), \end{equation}
  where $h = \sum_{i=1}^g\ord(b_{ii})$ (note that $B\nu\neq 0$ since
  $\det(B)\neq 0$).  The assertion is clear when
  $g=1$, so assume $g>1$.
  \begin{enum}
  \item If $\size(\nu) < \ord(\nu_1)$ then 
    $\size(\nu) = \size(\nu_2,\ldots,\nu_g)$, so we are done by induction
    on $g$.

  \item If $\size(\nu) = \ord(\nu_1)$ and 
    $\ord(\sum_{i=1}^g b_{1i}\nu_i)\leq\ord(b_{11}\nu_1)$ then we are
    done. 

  \item Otherwise, for some $i>1$ we have 
    $\ord(\nu_i)\leq \ord(\nu_i) + \ord(b_{1i})\leq \size(\nu) + \ord(b_{11})$.
    Hence if $\nu' = (\nu_2,\ldots,\nu_g)$, 
    $h' = \sum_{i=2}^g\ord(b_{ii})$, and
    $B' = (b_{ij})_{i,j\geq 2}$, then by induction we have 
    \[\size(B\nu) \leq \size(B'\nu')
      \leq h' + \size(\nu') 
      \leq h' + \ord(\nu_i) 
      \leq h' + \size(\nu) + \ord(b_{11}) 
      = h + \size(\nu). 
    \]
  \end{enum}
  This proves the assertion.
  
  Write $\dV = (a_{ij})$, and choose lifts $\td a_{ij}$ of $a_{ij}$ to $R$
  such that the matrix $\td\dV \coloneq (\td a_{ij})$ is also upper-triangular.
  Since $H = \sum_{i=1}^g\ord(\td a_{ii}) < 1$, the diagonal entries of
  $\td\dV$ are nonzero so by~(*)
  we have $\size(\td\dV(\xi^{(p)}))\leq H+p\size(\xi)$.
  We also have
  \[ H + (p-1)\size(\xi) < \frac{p-1}{p^{n+1}} + \frac 1p
  \leq \frac{p-1}{p^2} + \frac 1p
  = \frac 2p - \frac 1{p^2} \leq 1 - \frac 1{p^2} < 1, \]
  so $\size(\td\dV(\xi^{(p)})) \leq H + p\size(\xi) < \size(\xi)+1$.
  By Lemma~\ref{lem:pX}, 
  \begin{equation} \label{eq:pxi.form}
    [p](\bX) = \td\dV(X_1^p,\ldots,X_g^p) + pf(\bX) + O(\bX^{p^2})
  \end{equation}
  for some $f\in (\bX)R\ps\bX$.
  Write $[p](\bX) = ([p]_1(\bX),\ldots,[p]_g(\bX))$, and let
  $cX^\mu$ be a monomial occurring in $[p]_i(\bX)$ for some
  $1\leq i\leq g$ such that $c\xi^\mu\neq 0$ (i.e., $\xi_i\neq 0$ when
  $\mu_i\neq 0$).
  If $c\equiv 0\pmod p$ then
  $\size(\td\dV(\xi^{(p)})) < \size(\xi)+1 \leq \ord(c\xi^\mu)$.
  If instead $|\mu| \geq p^2$ then 
  \[ \ord(c\xi^\mu) \geq p^2\size(\xi) > \frac{p-1}{p^{n+1}} + p\size(\xi)
  > H + p\size(\xi) \geq \size(\td\dV(\xi^{(p)})), \]
  where the second inequality holds because 
  $p(p-1)\size(\xi) > (p-1)/p^{n+1}$ since $\size(\xi) > 1/p^{n+2}$.
  Hence in~\eqref{eq:pxi.form},
  we have $\size(f(\xi)) > \size(\td\dV(\xi^{(p)}))$, and
  the $O(\bX^{p^2})$-term also evaluates at $\bX=\xi$ with size greater
  than $\size(\td\dV(\xi^{(p)}))$.  Thus 
  $\size([p]\xi)=\size(\td\dV(\xi^{(p)}))$, and in particular
  $[p]\xi\neq 0$.
\qed

\begin{prop} Assume \label{prop:group.structure}
  that $H = H(G) < (p-1)/p^N$.  Then
  $\sG[p^\infty]_{\leq\rho_N'}(\bar K)\cong(\Z/p^N\Z)^g$.
\end{prop}

\pf
  First we show that for $1 \leq n\leq N$, 
  $\ms G[p^\infty]_{\leq\rho_n'} = \ms G[p^n]_{\leq\rho_n'}$.
  We will prove by induction on $n$ that
  $\ms G[p^\infty]_{\leq\rho_n'}(\bar K) = \sG[p^\infty]_{<\rho_n}(\bar K)$ 
  is killed by $p^n$ (cf.~\eqref{eq:cor.proper}). 
  Let $\xi\in\ms G[p^\infty]_{<\rho_n}(\bar K)$.  By Lemma~\ref{lem:pX},
  \[ \size([p]\xi)\geq\min\{p\size(\xi),~ 1+\size(\xi)\}
  > \min\{p\,r_{n}, 1+r_{n}\} 
  \geq \frac 1{p^{n-1}(p-1)}. \]
  If $n = 1$ then $\size([p]\xi) > 1/(p-1)$, so $[p]\xi = 0$
  by~\eqref{eq:cor.proper}.  Otherwise  
  $[p]\xi\in\ms G[p^\infty]_{<\rho_{n-1}}(\bar K)$, so 
  $[p^n]\xi = [p^{n-1}][p]\xi = 0$ by induction.%
  \footnote{See \cite[Lemma~2.2.6]{conrad:csg} for another proof of this
    fact.} 

  In order to prove that the $p^N$-torsion abelian group
  $\sG[p^\infty]_{\leq\rho_N'}(\bar K)$ of 
  order $p^{Ng}$ is isomorphic to $(\Z/p^N\Z)^g$, it suffices to show that 
  $\sG[p^\infty]_{\leq\rho_N'}(\bar K)$ has $p^g$ points of $p$-torsion.
  It is clear from the above that
  $\sG[p^\infty]_{\leq\rho_1'}\subset\sG[p^\infty]_{\leq\rho_N'}$ has
  order $p^g$ and is killed by $p$, so we must show that
  for $\xi\in\ms G[p^\infty]_{\leq\rho_N'}(\bar K)$, 
  $[p]\xi = 0$ implies $\xi\in\ms G[p^\infty]_{\leq\rho_1'}(\bar K)$.
  By~\eqref{eq:cor.proper},
  $\ms G[p^\infty]_{\leq\rho_1'} = \ms G[p^\infty]_{<\rho_1}$,
  so if $\xi\in\ms G[p^\infty]_{\leq\rho_N'}(\bar K)$ is not contained in
  $\ms G[p^\infty]_{\leq\rho_1'}(\bar K)$ then $\size(\xi)\leq 1/p(p-1)$.
  By Lemma~\ref{lem:pdoesntkill}, this implies that $[p]\xi\neq 0$.
\qed

The proof of Theorem~\ref{thm:canonicalsubgp} is now complete.

\section{Elimination of noetherian hypotheses}
\label{sec:not.noetherian}

\paragraph
As mentioned in~\parref{par:reduction1}, the goal of this section is
to prove:

\begin{prop} Let \label{prop:comes.from.display}
  $G$ be a connected level-$N$ truncated $p$-divisible group over $R$.
  If the fraction field $K$ of $R$ is algebraically closed then
  $G$ extends to a $p$-divisible group over $R$.
\end{prop}

In particular, there is a display $\cP$ over $R$ such that
$G\cong\BT_\cP[p^N]$. 

We will use a standard noetherian approximation argument.

\begin{prop} Suppose \label{prop:extend.resfield}
  that the fraction field $K$ of $R$ is algebraically closed.
  Let $(A,\fm)$ be a local noetherian ring with  residue
  field $\kappa$ and let $\phi:A\to R$ be a local homomorphism.  Let
  $\kappa'$ be a subfield of the residue field $k$ of $R$, and assume
  $\kappa\subset\kappa'$.  There exists a flat local noetherian
  $A$-algebra $B$ and commutative diagram of local homomorphisms
  \[\xymatrix{
    {B} \ar[r]^{\psi} & {R} \\
    {A} \ar[u] \ar[ur]_{\phi}
  }\]
  such that $\fm B$ is the maximal ideal of $B$, and the map
  $B/\fm B\to k$ has image $\kappa'$.
\end{prop}

\pf
  The existence of a flat local noetherian $A$-algebra $B$ such that
  $B/\fm B$ is isomorphic to a given field extension of $\kappa$ is a
  standard fact proved in~[EGA0${}_{\text{III}}$, Proposition~10.3.1]. 
  \nocite{egaIII_1}
  Here we give an indication of how to modify that proof to include a
  construction of the map $\psi:B\to R$. 

\subparagraph Suppose that $\kappa' = \kappa(t)$ where $t\in\kappa'$ is a
    transcendental element.  Let $A' = A[T]$, let $\fp = \fm A'$, and let
    $B = A'_\fp$.  Let $\psi:A'\to R$ be the
    $A$-algebra homomorphism sending $T\mapsto t$.  If $f(T)\notin\fp$
    then the residue of $\psi(f(T))$ in $k$ is nonzero because $t$ is
    transcendental over $\kappa$.  Thus $\psi$ extends to an
    $A$-map $B\to R$, and we are done in this case.

\subparagraph Suppose that $\kappa' = \kappa(\alpha)$ where $\alpha\in\kappa'$ is an
    algebraic element with (monic) minimal polynomial 
    $f\in\kappa[T]$.  Let $F\in A[T]$ be a monic polynomial reducing to
    $f\in k[T]$, let $A' = A[T]$, and let $B = A'/(F)$.  Since $K$
    is algebraically closed, there exists
    $\td\alpha\in R$ lifting $\alpha$ such that $F(\td\alpha)=0$.
    The map $\psi:A'\to R$ sending $T\mapsto\td\alpha$ factors through an
    $A$-map $\psi:B\to R$, which settles this case.

\subparagraph In general there exists an ordinal $\gamma$ and, for all
    ordinals $\lambda\leq\gamma$, a subfield $\kappa_\lambda$ of $\kappa'$
    containing $\kappa$ such that (a)~for each $\lambda<\gamma$,
    $\kappa_{\lambda+1}$ is an extension of $\kappa_\lambda$ generated by
    a single element, (b)~for every ordinal $\mu$ without a predecessor,
    we have $\kappa_\mu = \bigcup_{\lambda<\mu}\kappa_\lambda$, and
    (c)~$\kappa = \kappa_0$ and $\kappa'=\kappa_\gamma$.  By
    transfinite recursion we will construct local noetherian rings
    $B_\lambda$ for $\lambda\leq\gamma$ and local homomorphisms
    $\sigma_{\mu\lambda}: B_\lambda\to B_\mu$ for $\lambda\leq\mu$ and
    $\psi_\lambda: B\to R$ such that 
    \begin{deflist}
    \item[(1)] $(B_\lambda,\sigma_{\lambda\mu})$ is a directed system with
      $B_0=A$.
    \item[(2)] For all $\lambda$ the map $\psi_\lambda$ induces a
      $\kappa$-isomorphism $B_\lambda/\fm B_\lambda\isom\kappa_\lambda$.
    \item[(3)] For $\lambda\leq\mu$, $B_\mu$ is 
      $B_\lambda$-flat. 
    \end{deflist}
    Let $\xi\leq\gamma$, and suppose that $B_\lambda$,
    $\sigma_{\mu\lambda}$, and $\psi_\lambda$ have been constructed
    satisfying~(1)--(3) for $\lambda\leq\mu < \xi$.  If 
    $\xi=\mu+1$ is a successor then $k_\xi$ is generated over 
    $k_\mu$ by a single element, so we can construct $B_\xi$ and
    $\psi_\xi$ as in the previous two paragraphs.  If 
    $\xi$ is not a successor then we set $B_\xi = \dlimm_{\mu<\xi} B_\mu$,
    and we let $\psi_\xi:B_\xi\to R$ be the natural map.  Then $B_\xi$
    satisfies~(1)--(3) by~[EGA0${}_{\text{III}}$, Lemma~10.3.1.3].\qed
    \nocite{egaIII_1}

\begin{prop} Suppose \label{prop:descend.BTn}
  that the fraction field $K$ of $R$ is algebraically closed.
  Let $G$ be a $\BT_n$ over $R$ for some $n\geq 1$.  
  There exists a complete local noetherian ring
  $R'$, a local homomorphism $R'\to R$ inducing an isomorphism on residue
  fields, and a level-$n$ truncated $p$-divisible group $G'$ over $R'$
  such that $G'\tensor_{R'} R \cong G$.  
\end{prop}

\pf
  First we show that there exists a local noetherian subring $R_0$ of $R$
  with local structure map $R_0\to R$ and a
  finite flat $R_0$-group scheme $G_0$ such that 
  $G_0\tensor_{R_0} R\cong G$.  We will use the standard techniques of
  noetherian approximation from~\cite[\S8]{egaIV_3}.

\subparagraph Let $G = \Spec(A)$, and let
  $m:A\tensor_R A\to A$ and $\mu:A\to A\tensor_R A$ denote the
  multiplication and comultiplication, respectively.  Let
  $\iota:A\to R$ be the coidentity, and let $I = \ker(\iota)$, so
  $A = R\dsum I$.  Choosing an $R$-basis
  for $I$ we obtain isomorphisms $I\cong R^M$ and $A\cong R^{M+1}$, and
  hence matrix representations $m = (m_{ij})$ and $\mu = (\mu_{ij})$.  The
  augmentation ideal $I$ and the structure coefficients $m_{ij},\mu_{ij}$
  determine the Hopf algebra structure on $A$.  Let $R_0$ be a subring of
  $R$ that is finitely generated over $\Z$ and contains the $m_{ij}$ and
  $\mu_{ij}$, let $I_0 = R_0^M$ and let  
  $A_0 = R_0\dsum I_0$.  We define maps
  $m_0:A_0\tensor_{R_0}A_0\to A_0$ and $\mu_0:A_0\to A_0\tensor_{R_0}A_0$ using the
  matrices $(m_{ij})$ and $(\mu_{ij})$, respectively, and we let
  $R_0\to A_0 = R_0\dsum I_0$ (resp. $\iota_0: A_0\to R_0$) be inclusion into
  (resp. projection onto) the first factor.  It is not hard to see that
  these data endow $A_0$ with the structure of Hopf algebra over $R_0$.  If
  $G_0=\Spec(A_0)$ then $G_0\tensor_{R_0} R\cong G$ by construction, and
  $G_0$ is commutative since $R_0\to R$ is an injection.  Replacing $R_0$
  with its localization at $R_0\cap\fm_R$, we may assume that $R_0$ is local.

\subparagraph Let \label{par:directed.system.R}
  $\{R_\alpha\}_{\alpha\in\cI}$ be the directed system of local noetherian
  subrings of $R$ containing $R_0$ and having local structure map
  $R_\alpha\to R$, and for $\alpha\in\cI$ let 
  $G_\alpha = G_0\tensor_{R_0} R_\alpha$.  Suppose that $n\geq 2$, and let 
  $0\leq i\leq n$.  Then $[p^i]_{G_\alpha}: G_\alpha\to G_\alpha[p^{n-i}]$
  is faithfully flat for large enough $\alpha$
  by~\cite[Theorems~8.10.5 and~11.2.6]{egaIV_3}, so $G_\alpha$ is a
  $\BT_n$.  A similar argument shows that $G_\alpha$ is a $\BT_1$ for
  large enough $\alpha$ when $n=1$.  

\subparagraph
  Fix a large $\alpha$ as in~\parref{par:directed.system.R} and let 
  $R_1 = R_\alpha$.  This is a local
  noetherian ring with maximal ideal $\fm_1$ and local structure map 
  $R_1\to R$.  Let $R_2$ be a local noetherian ring with maximal ideal
  $\fm_2$, equipped with a local homomorphism $R_1\to R_2$ such that
  $\fm_2=\fm_1 R_2$ and a local $R_1$-algebra homomorphism
  $R_2\to R$ inducing an isomorphism on residue fields, as in
  Proposition~\ref{prop:extend.resfield}.  For 
  $r\in\R_{\geq 0}$ let $\fa_r = \{ x\in R_2~:~\ord(x)\geq r\}$, so
  $\fa_1\subset\fa_{1/2}\subset\fa_{1/3}\subset\cdots$ with
  $\bigcup_{i=1}^\infty \fa_{1/i} = \fm$.  
  Since $R_2$ is noetherian, there is some $\pi\in R$ with nonzero
  valuation such that $\ord(x) \leq \ord(\pi)$ for all $x\in\fm_2$.  It
  follows that $\ord(x)\leq\ord(\pi^n)$ for all $x\in\fm_2^n$, so the
  $\fm_2$-adic completion $R'$ of $R_2$ maps into $R$.  The
  ring $R'$ and the truncated $p$-divisible group
  $G' = G_\alpha\tensor_{R_\alpha} R'$ satisfy the properties of
  Proposition~\ref{prop:descend.BTn}. 
\qed

\paragraph\pf[of Proposition~\ref{prop:comes.from.display}]
Let $G$ be a connected $\BT_N$ over $R$ for some $N\geq 1$.
Since $K$ is algebraically closed, the residue field $k$ of $R$ is
perfect.  Let $R'$ be a complete 
local noetherian ring with a local homomorphism $R'\to R$ inducing an
isomorphism of residue fields and such that there exists a connected
truncated $p$-divisible group $G'$ of level $N$ over $R'$ with
$G'\tensor_{R'} R\cong G$, as in Proposition~\ref{prop:descend.BTn}.  
By~\cite[Theorem~4.4(e)]{illusie:pdivgp.deformation}, there is a connected
$p$-divisible group $H'$ over $R'$ such that $H'[p^N]\cong G'$.  Let
$H = H'\tensor_{R'} R$.  Then
\[ G\cong G'\tensor_{R'} R\cong(H'\tensor_{R'} R)[p^N]\cong H[p^N]. \]
\qed

\bibliographystyle{thesis}
\bibliography{thesis}

\end{document}